\documentclass[12pt]{article}
\usepackage[left=2cm,top=2cm,right=2cm,bottom=2cm]{geometry}
\usepackage{latexsym,amssymb,amsmath,color,morefloats}
\usepackage{subfig}
\usepackage[pagewise]{lineno}
\usepackage{natbib}
\usepackage{caption,graphicx,hyperref,theorem}
\usepackage{epstopdf,chngcntr}

\pdfoutput=1

\ifpdf
  \DeclareGraphicsExtensions{-eps-converted-to.pdf,.pdf,.png,.jpg}
\else
  \DeclareGraphicsExtensions{-eps-converted-to.pdf}
\fi

\newtheorem{theorem}{Theorem}

\newtheorem{proposition}[theorem]{Proposition}

{\theorembodyfont{\rmfamily}

}

\numberwithin{equation}{section}

\hypersetup{
    colorlinks = true,
    linkcolor = {blue},
    citecolor = {blue},
    urlcolor  = {blue},}

\definecolor{manchester}{rgb}{.42,.17,.58}

\newcommand{\RE}[1]{\textup{Re}\left\{#1\right\}}
\newcommand{\IM}[1]{\textup{Im}\left\{#1\right\}}

\newcommand{\NN}{{\mathbb{N}}}

\newcommand{\pii}{\pi i}

\newcommand{\Der}[2]{\frac{\text{d}#1}{\text{d}#2}}
\newcommand{\ParDer}[2]{\frac{\partial#1}{\partial#2}}
\newcommand{\tmop}[1]{\ensuremath{\operatorname{#1}}}
\newcommand{\mathd}{\mathrm{d}}

\newcommand{\nin}{\not\in}
\newcommand{\tmem}[1]{{\em #1\/}}

\newcommand{\RED}{} 
\newcommand{\mylinenum}{} 

\newcommand{\tmdummy}{$\mbox{}$}
\newcommand{\tmmathbf}[1]{\ensuremath{\boldsymbol{#1}}}

\begingroup 
\makeatletter
\gdef\eqna@origamp{&} 
\catcode`\&\active 
\gdef\eqna@newamp{%
  \ifx\@currenvir\eqna@currenvir 
    \eqna@onlyfirstamp\let\eqna@onlyfirstamp\@empty 
  \else 
    \eqna@origamp 
  \fi
}
\gdef\eqna@hook{
  \let\eqna@currenvir\@currenvir 
  \catcode`\&\active 
  \let&\eqna@newamp 
  \let\eqna@onlyfirstamp\eqna@origamp 
  }
\catcode`\*11 
\gdef\eqnarray{\eqna@hook\align} 
\gdef\eqnarray*{\eqna@hook\align*} 
\global\let\endeqnarray\endalign
\global\let\endeqnarray*\endalign*
\endgroup

\title{Analytical methods for perfect wedge diffraction: a review}

\author{Matthew A. Nethercote$^{*}$, Raphael C. Assier$^{*}$ and I. David Abrahams$^{\dagger}$\\
\footnotesize{$^{*}$ School of Mathematics, University of Manchester, Oxford Road, Manchester, {\rm M13 9PL}, UK}\\
\footnotesize{$^{\dagger}$ Isaac Newton Institute, University of Cambridge, 20 Clarkson Road, Cambridge CB3 0EH, UK}
}

\begin{document}

\maketitle

\begin{abstract}{
The subject of diffraction of waves by sharp boundaries has been studied intensively for well over a century, initiated by groundbreaking mathematicians and physicists including Sommerfeld, Macdonald and Poincar\'e. The significance of such canonical diffraction models, and their analytical solutions, was recognised much more broadly thanks to Keller, who introduced a geometrical theory of diffraction (GTD) in the middle of the last century\RED{, and other important mathematicians such as Fock and Babich}. This has led to a very wide variety of approaches to be developed in order to tackle such two and three dimensional diffraction problems, with the purpose of obtaining elegant and compact analytic solutions capable of easy numerical evaluation. 

The purpose of this review article is to showcase the disparate mathematical techniques that have been proposed. For ease of exposition, mathematical brevity, and for the broadest interest to the reader, all approaches are aimed at one canonical model, namely diffraction of a monochromatic scalar plane wave by a two-dimensional wedge with perfect Dirichlet or Neumann boundaries. The first three approaches offered are those most commonly used today in diffraction theory, although not necessarily in the context of wedge diffraction. These are the Sommerfeld-Malyuzhinets method, the Wiener-Hopf technique, and the Kontorovich-Lebedev transform approach. Then follows three less well-known and somewhat novel methods, which would be of interest even to specialists in the field, \RED{namely the embedding method, a random walk approach, and the technique of functionally-invariant solutions.}

Having offered the exact solution of this problem in a variety of forms, a numerical comparison between the exact solution and several powerful approximations such as GTD is performed and critically assessed.
}\end{abstract}

\mylinenum{
\section{Introduction and formulation}
At the close of the $19^{\text{th}}$ century, wedge diffraction became a core problem in mathematical physics when renowned mathematicians Poincar\'{e} and Sommerfeld studied the diffraction of wave fields in angular domains \citep{Poincare1892part1,Poincare1892part2,Sommerfeld1896,Sommerfeld1901}. Sommerfeld made the first breakthrough when he solved his famous half-plane problem \citep{Sommerfeld1896}, during which he introduced the contour integral representation that we now know as the Sommerfeld \RED{integral. This work has now been translated to English in \citet{Sommerfeld2003}, with additional insightful comments. Sommerfeld} would later be the first to solve problems of wedge diffraction \citep{Sommerfeld1901} where the wedge has  an interior angle equal to $m\pi/n$ ($m<n\in\NN$). 

For wedges with arbitrary interior angles, the solution was first obtained by \citet{Macdonald1902}. He did this by considering a line source incident wave and used separation of variables to get a series solution. The solution was rewritten in Sommerfeld integral form and he then provided the solution for an incident plane wave. \RED{We discuss this line source approach and provide an alternative way to obtain the plane wave solution in \ref{Mac}.}

In the 1950s, Malyuzhinets released a series of papers that culminated in the solution to the problem with impedance boundary conditions, \citep{Malyuzhinets1955-1-russian,Malyuzhinets1955-2-russian,Malyuzhinets1958-1,Malyuzhinets1958-2,Malyuzhinets1958-3}. This result created the first method that we discuss here in Section \ref{SMT}, the Sommerfeld-Malyuzhinets technique \RED{(S-M)}. Other authors who solved the impedance wedge problem independently were \citet{Senior1959} and \citet{Williams1959} but for more details on Malyuzhinets' method, see the review paper \citep{NorrisOsipov1999} or the books \citep{Budaev1995,SMtechnique2007,LyalinovZhu2013}.

One of the most popular methods in diffraction theory is the Wiener-Hopf (W-H) technique, invented by \citet{WienerHopf1931} as a means to solve a special type of integral equation. It was soon discovered to be a useful method for diffraction problems and has appeared in a number of classic articles such as \citep{Copson1946} and \citep{DSJones1952}. Since then, applications of the technique have appeared in a wide array of research areas including diffraction, waveguides and flow problems.

The well-known textbook \citep{Noble1958} provides an excellent tutorial for various aspects and extensions of the W-H technique. In 2007, the Journal of Engineering Mathematics published a W-H special issue led by a historical overview \citep{LawrieAbrahams2007} along with a collection of articles applying the W-H technique to various problems. For wedge problems, the technique was thought to be ineffective due to the two boundaries not being parallel, however (see e.g. \citep{AVShanin1996,Daniele2003}) this can be overcome as discussed later in Section \ref{WHT}.

Another key method is based on the Kontorovich-Lebedev (K-L) transform. First introduced by \citet{KL1939}, this transform is an effective tool when dealing with a radial coordinate. This makes it useful for wedge diffraction problems as evidenced by \citet{Abrahams1986,Abrahams1987}, since obtaining the general solution in that way is a very natural process. We shall discuss this further in Section \ref{MKLT}, but for more details on the transform and its applications, see \citep{Lebedev1965}, \citep{DSJones1964,DSJones1980,DSJones1986} and \citep{Felsen1994}.

In Section \ref{SolAnalysis}, we will focus on the asymptotic technique created by the classic paper \citep{Keller1962}, called the Geometrical Theory of Diffraction \RED{(GTD), see also the book by \citet{BorovikovKinber1994}}. We will also follow the \RED{uniform} GTD extension detailed in literature such as \citep{KP1974,GLJames1986,Pistorius1990} and \citep{SMtechnique2007}. 

\RED{Even though we will not use it in this review, for completeness, it is important to mention an alternative asymptotic technique applied to diffraction problems that is the Physical Theory of Diffraction (PTD) \citep{Ufimtsev1971}. This development was made possible in part thanks to Macdonald's work on Kirchhoff's approximation (see e.g. \citet{Ufimtsev2014}).} A useful paper that compares the GTD and PTD asymptotic techniques as well as the exact solution in series and integral form is \citep{HSU2011}. \RED{Similar methods, describing creeping waves in diffraction by smooth obstacle, have also been developed in \citet{Fock1965} for example.}

Section \ref{Alt} contains a number of alternative methods that are effective but less well-known for wedge diffraction. The \RED{first} of these alternative methods is based on the very powerful concept of embedding formula. This reasonably recent approach consists in expressing the diffraction coefficient (which depends on both the incident and observer angles) of the diffracted field resulting from an incident plane wave in terms of the directivities (depending on one angle only) of simpler problems. These simpler problems are directly related to edge Green's functions. These are Green's functions for which the source is sent towards the geometric singularities of the obstacle. This method was primarily used for planar cracks and slits, and parallel combinations of these (see e.g. \citep{williams,Gautesen1983,Martin1983,Biggs2001,Biggs2002}). \RED{In} \citep{ShaninCraster2005} it was shown that the method can be successfully adapted to wedges, as we will discuss later.

\RED{The second of these alternative methods is the so-called random walk approach.} It is based on the known link between deterministic PDEs and stochastic differential equations (SDEs) given by the Feynman-Kac theorem. It allows to express the solution of a diffraction problem as the mean of a set of solutions to \RED{given SDEs} with carefully chosen initial and final conditions. The method was developed through a series of papers by \citet{BudaevBogy2001,BudaevBogy2002a,BudaevBogy2002b,BudaevBogy2003}, the latter being dedicated to wedge diffraction.

The \RED{last} of these is the method of functionally-invariant solutions also known as the Sobolev-Smirnov method which has been used for a number of plane wave diffraction problems from half-planes \citep{Sobolev1935,Smirnov1964} to wedges \citep{Filippov1964,KMM2015,Babich2015}. A very similar method that develops Busemann's ``conical flow method'' \citep{Busemann1947} was also considered in \citet{KellerBlank1951} \RED{ and \citet{Miles1952-1} for example}.

In this review \RED{(apart from Macdonald's approach discussed in \ref{Mac})}, we will focus primarily on plane wave incidence \RED{rather than line sources. It has to be noted, however, that a broad range of work \citep{Bromwich1915,Oberhettinger1954,Rawlins1987,Rawlins1989}} has also been carried out for both acoustic and electromagnetic sources. \RED{For other reviews of some of the methods used for various types of incident waves (plane, cylindrical, spherical, dipole and pulse), see \citet{Oberhettinger1958} and \citet{Bowman1987}.}

\RED{The elastic wedge problem has equally received a lot of attention. \citet{Knopoff1969} wrote an interesting review of possible approaches to tackle this (still unsolved) problem, it includes attempts using the method of images, the W-H technique, the K-L transform and the conical flow method
. More recent approaches by \citet{CroisilleLebeau1999} or \citet{BudaevBogy1995,BudaevBogy1996,BudaevBogy1998} are also worth mentioning.}


\RED{In Section} \ref{WHM-mapping-section}, in the spirit of \citet{Wegert2012}, we will visualise complex functions using phase \RED{portraits} in order to show domains of analyticity, locations of singularities and orientations of branch cuts. \RED{These portraits assign the argument of a complex function to a HSV colour model.} For example, Figure \ref{wedgegeometry} (left) shows the \RED{phase portrait of $f(z)=z$} which we use as a colour reference. 


Wedge diffraction has a number of physical applications such as the scattering of acoustic pressure fields or electromagnetic fields by sharp structures, scattering of the Sun's radiation by cloud ice crystals and seismology. We shall study this from a more mathematical perspective. From here onward, we will assume that the problem is time-harmonic with time factor $e^{-i\omega t}$, which is therefore suppressed, and we will consider solutions to the homogeneous Helmholtz equation,}
\begin{align}
\label{Intro-Helmholtz}\nabla^2\Phi+k^2\Phi=0.
\end{align}
\mylinenum{inside a wedge-shaped region described in polar coordinates by $\{0<r<\infty,\ -\theta_\text{w}<\theta<\theta_\text{w}\}$ (see Figure \ref{wedgegeometry} (right)). The complementary region, $\{0<r<\infty,\ |\theta|>\theta_\text{w}\}$ is considered to be the wedge scatterer. Defining $\bar{\theta}_\text{w}=\pi-\theta_\text{w}$, the interior angle of the wedge scatterer is $2\bar{\theta}_\text{w}$. Throughout this paper, we will consider two cases of homogeneous boundary conditions \RED{(BCs)}, Dirichlet and Neumann, on \RED{both faces} of the wedge,}
\begin{align}
\label{Intro-DBC}\hbox{Dirichlet \RED{BCs:}}&\ \ \Phi(\theta=\pm\theta_{\text{w}})=0,\\
\label{Intro-NBC}\hbox{Neumann \RED{BCs:}}&\ \ \frac{1}{r}\ParDer{\Phi}{\theta}(\theta=\pm\theta_{\text{w}})=0.
\end{align}
\mylinenum{For acoustics, \eqref{Intro-DBC} and \eqref{Intro-NBC} are called sound-soft and sound-hard \RED{BCs} respectively. For electromagnetic scattering \RED{by an electric (resp. magnetic) polarized plane wave, the solution to the perfect electric conducting (PEC) problem can be expressed in terms of a potential satisfying \eqref{Intro-DBC} (resp. \eqref{Intro-NBC}).} 

\RED{We} define the incident plane wave as, $\Phi_\text{I}=e^{-ikr\cos(\theta-\theta_\text{I})}$ with wavenumber $k>0$ and incident angle $\theta_\text{I}$. Due to the symmetry of the problem, we restrict the incident angle to $\theta_\text{I}\in[0,\theta_{\text{w}}]$. Figure \ref{wedgegeometry} (right) illustrates the geometry of the problem.
\begin{figure}[ht]\centering
\includegraphics[width=0.44\textwidth]{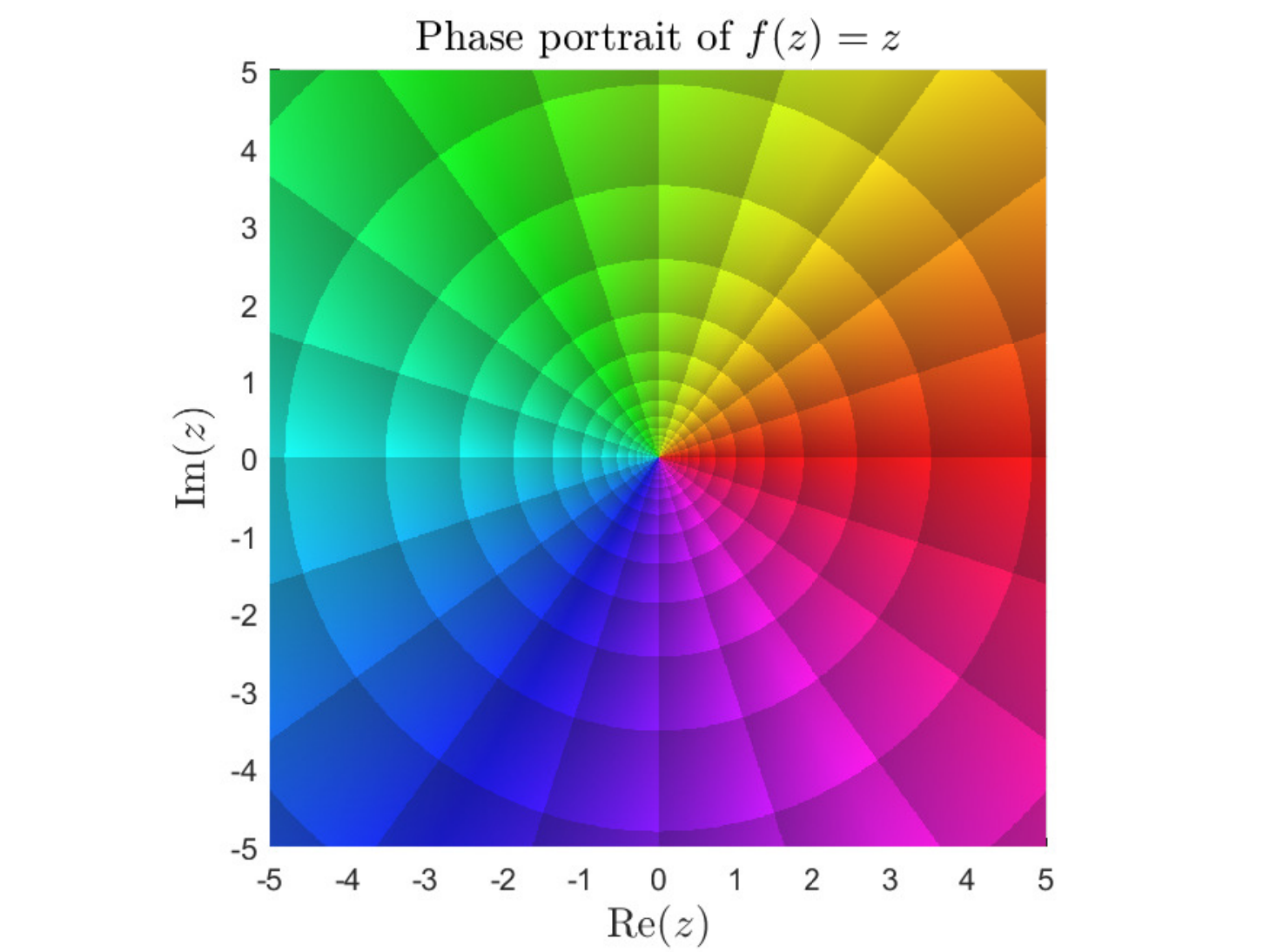}
\includegraphics[width=0.48\textwidth]{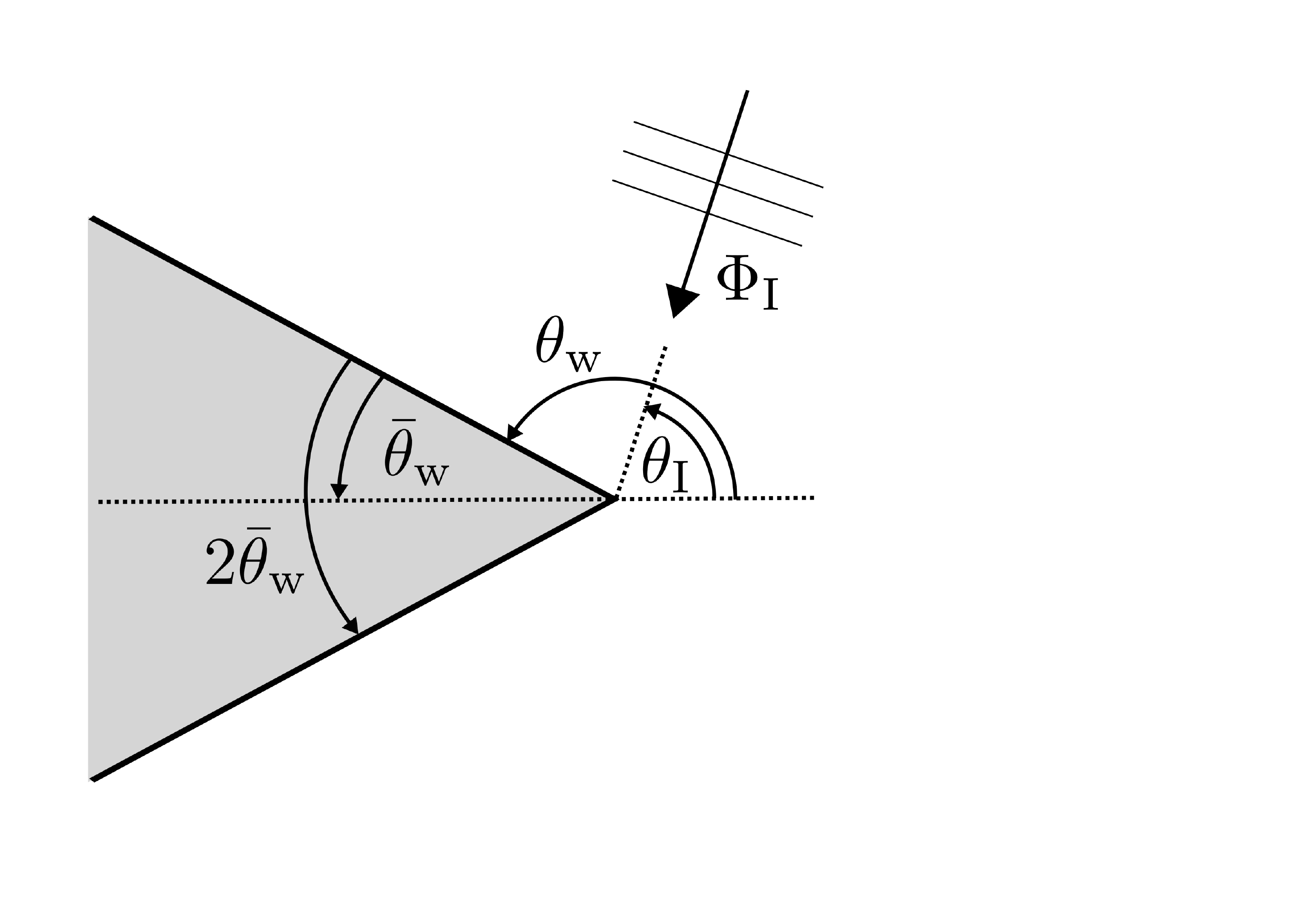}
\caption{\RED{The left figure is a phase portrait of $f(z)=z$ which assigns colours to the complex argument of $f$ (left). For example arg$(f)=0$ is indicated as red and arg$(f)=\pi$ is indicated as cyan. The right figure is the geometry of the problem.}}
\label{wedgegeometry}
\end{figure}

An initial approximation is found using classic Geometrical Optics (GO). The GO part of the solution consists of the incident wave and any reflected waves produced. The rest of the solution is considered to be the diffracted field, $\Phi_\text{Diff}$ which satisfies \RED{a two-dimensional radiation condition (see \citep{SCHOT1992} for a good review on this), written in integral form:}}
\begin{align}\label{Intro-2DSRcondition}
\RED{\lim_{r\rightarrow\infty}
\int_{-\theta_{\text{w}}}^{\theta_{\text{w}}} \left|\ParDer{\Phi_\text{Diff}}{r}-ik\Phi_\text{Diff}\right|^2 r\, \text{d}\theta=0.}
\end{align}

\mylinenum{Lastly, there will be an edge (or Meixner) condition as $r$ becomes small.}
\begin{align}
\label{Intro-Edge}\Phi\sim \RED{\mathcal{A}}+O\left(r^{\RED{\min(\delta,2)}}\right),\ \ \text{where}\ \ \delta=\frac{\pi}{2\theta_\text{w}},
\end{align}
\mylinenum{and $\RED{\mathcal{A}}=0$ for the case of Dirichlet \RED{BCs}. Typically, the edge conditions can be derived using the Frobenius method \RED{\citep{BenderOrszag1999}} while ensuring that the energy remains finite in any neighbourhood of the wedge edge.


\section{The Sommerfeld-Malyuzhinets technique}\label{SMT}  
The first method to be reviewed is the Sommerfeld-Malyuzhinets \RED{(S-M)} technique. Here we will \RED{show briefly how to get} the solution to the perfect wedge problem (\ref{Intro-Helmholtz})-(\ref{Intro-Edge}). For a more thorough explanation, consult Sections 1-4 of \citep{SMtechnique2007}. This technique is based on the general solution of diffraction problems in angular domains being represented as the Sommerfeld integral}
\begin{align}
\Phi(r,\theta)&=\frac{1}{2\pi i}\int_{\gamma_+}e^{-ikr\cos(z)}\left[s(\theta+z)-s(\theta-z)\right]\text{d}z
=\frac{1}{2\pi i}\int_{\gamma_++\gamma_-}\!\!\!\!\!\!\!\!e^{-ikr\cos(z)}s(\theta+z)\text{d}z,\label{SMTSommInt}
\end{align}
\mylinenum{where $\gamma_{\pm}$ are contours defined in Figure \ref{SMTSommcon} (left) and $s(z)$ is an unknown function to be determined. \RED{This representation ensures that the Helmholtz equation is automatically satisfied\footnote{\RED{This is proven using integration by parts and noting that $e^{-ikr\cos(z)}$ satisfies the Helmholtz equation with polar coordinates $(r,z)$.}}. The form of the spectral part, $s(\theta+z)-s(\theta-z)$, is necessary for the radiation conditions to be satisfied. This can be proven by using the method of steepest descent to approximate the integral as $kr\rightarrow\infty$ (see section 3.7 of \citet{SMtechnique2007}).}
The function $s(z)$, referred to as the \RED{spectral function}, is assumed to be meromorphic in the domain}
\begin{align}\label{SMTwidehalfstrip}
\nonumber&\left\{-\pi-\theta_{\text{w}}-\epsilon_1<\RE{z}<\theta_{\text{w}}+\epsilon_1,\ \ \IM{z}>-\epsilon_2\right\}\cup\\
&\left\{-\theta_{\text{w}}-\epsilon_1<\RE{z}<\pi+\theta_{\text{w}}+\epsilon_1,\ \ \IM{z}<\epsilon_2\right\},
\end{align}
\mylinenum{for some $\epsilon_{1,2}>0$ and analytic in the same domain with $\epsilon_{1,2}=0^+$, see Figure \ref{SMTSommcon} (right). \RED{The poles of $s(z)$ will be seen to correspond to the geometrical optics part of the wave field}. The structure of \eqref{SMTSommInt} means that $s(z)$ is defined up to an additive constant. \RED{Using the edge conditions, we can assume} that $s(z)$ has the following behaviour as $|\IM{z}|\rightarrow\infty$,}
\begin{align}\label{SMTsasym}
s(z)&=\pm A+O(e^{-\delta|\IM{z}|}),
\end{align}
\mylinenum{where $A=0$ for the Dirichlet case.\footnote{\RED{We assume \eqref{SMTsasym} for convenience later. Say we assumed $s(z)\rightarrow A_\pm$ as $\IM{z}\rightarrow\pm\infty$ where $A_++A_-\neq0$ instead. Then a later step will require us to redefine $s(z)$ using the additive constant property such that $A_++A_-=0$ implying \eqref{SMTsasym}.}} It is also important to note that $e^{-ikr\cos(z)}$ is an entire $2\pi$-periodic function of $z$, which decays rapidly as $|\IM{z}|\rightarrow\infty$ only in the set of half-strips,}
\begin{align}\label{SMThalfstrips}
\nonumber&\{z:\ (2m-1)\pi<\RE{z}<2m\pi,\ \ \IM{z}<0,\ \ m\in\mathbb{Z}\}\cup\\
&\{z:\ 2m\pi<\RE{z}<(2m+1)\pi,\ \ \IM{z}>0,\ \ m\in\mathbb{Z}\},
\end{align}
\mylinenum{displayed in Figure \ref{SMTSommcon} (left), and grows rapidly in the complementary set. The Sommerfeld contours $\gamma_{\pm}$ are defined so that the integrand is analytic and decays as $|\IM{z}|\rightarrow\infty$ along these contours.\footnote{Note that the Sommerfeld contours are contained in the domain \eqref{SMTwidehalfstrip} with $\epsilon_1=\epsilon_2=0^+$.}
\begin{figure}[ht]\centering
\includegraphics[width=0.45\textwidth]{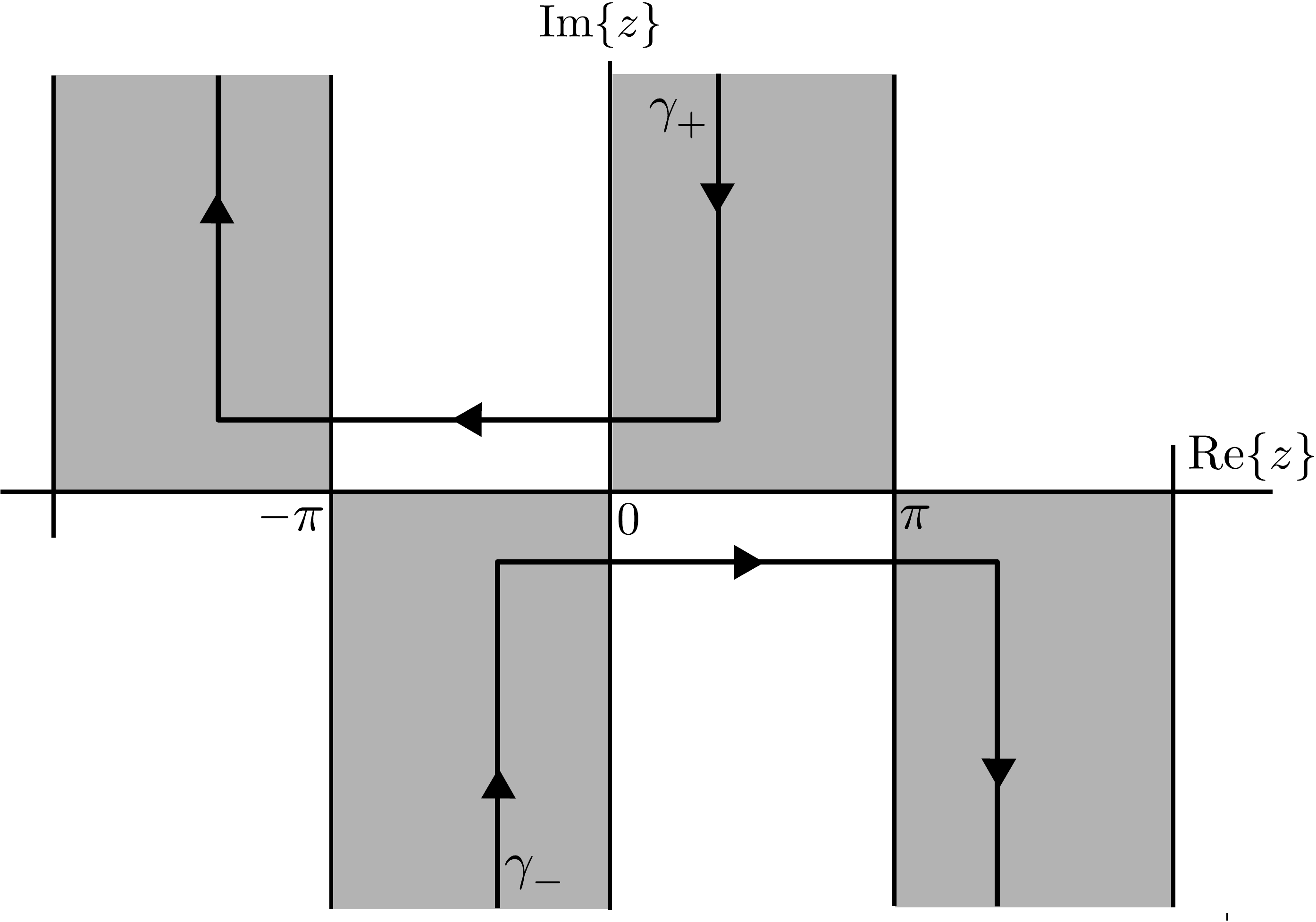}
\includegraphics[width=0.35\textwidth]{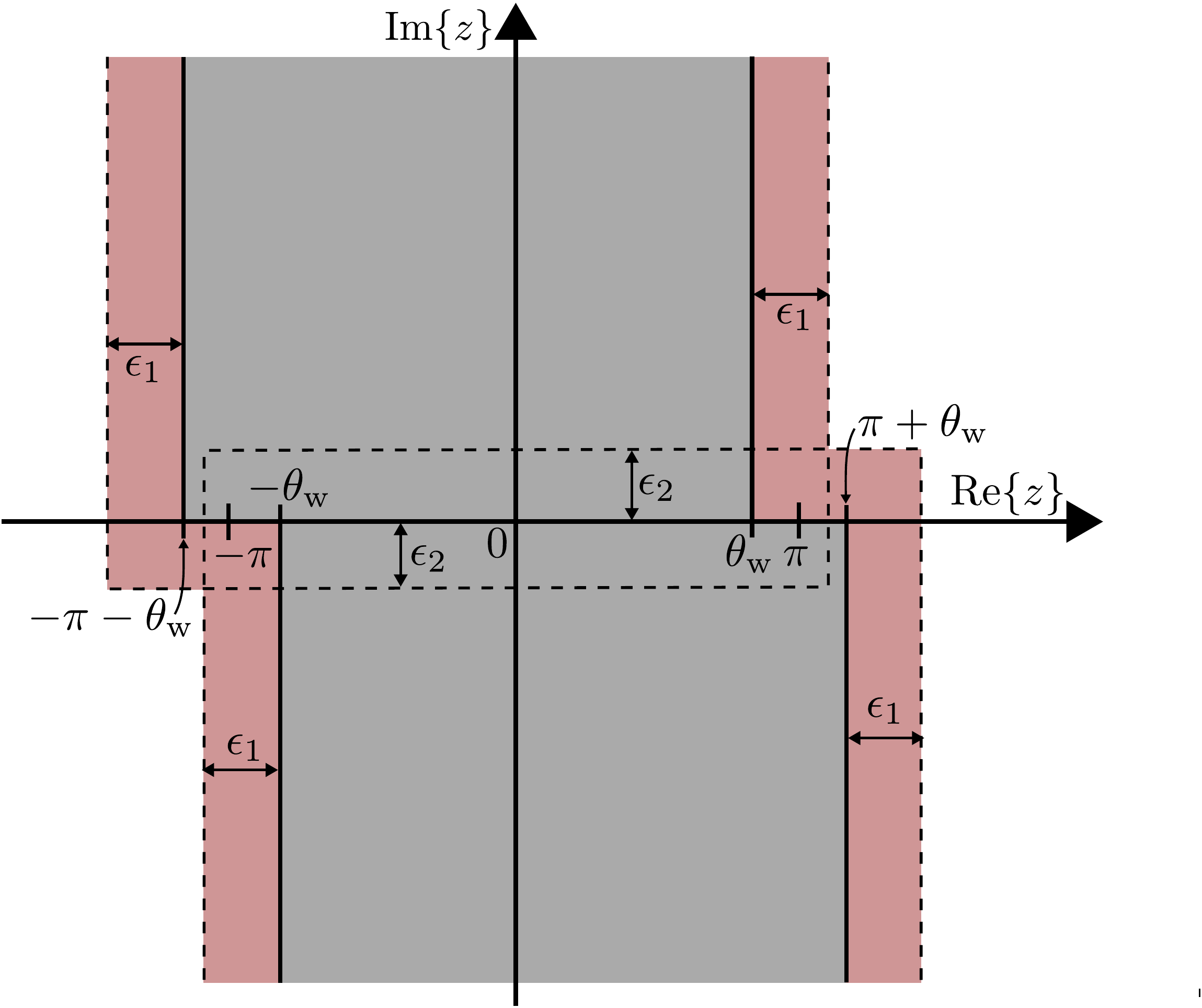}
\caption{The Sommerfeld contours $\gamma_\pm$ and the grey half-strips of exponential decay (left), and regions where $s(z)$ is analytic or meromorphic (right)}
\label{SMTSommcon}
\end{figure}

A crucial part to the \RED{S-M} technique is Malyuzhinets' Theorem or the Sommerfeld Nullification Theorem. This is an important theorem because it allows to obtain a functional equation satisfied by the \RED{spectral function}. The theorem and its proof are presented by Malyuzhinets in \citep{Malyuzhinets1958-1} and more recently in Section 3.4 of \citep{SMtechnique2007}.}

\begin{theorem}[Malyuzhinets' Theorem or Sommerfeld Nullification Theorem] \label{th:MalThm}
\mylinenum{Let the function $\Upsilon(z)$ be analytic and single-valued inside the half-strip,}
\begin{align}\label{MalThmStrip}
\{z:-\pi-\epsilon_1\leq\textup{Re}\{z\}\leq\epsilon_1,\ \ \textup{Im}\{z\}\geq\epsilon_2>0\},\ \ \epsilon_1,\epsilon_2>0.
\end{align}
\mylinenum{If for some constant $D$, the function has the following behaviour as $\IM{z}\rightarrow\infty$ in this half-strip,}
\begin{align}\label{MalThmAsym}
|\Upsilon(z)|\leq\mathrm{constant}\  e^{D\IM{z}},
\end{align}
\mylinenum{and for any $R>0$,}
\begin{align}\label{MalThmInt}
\frac{1}{2\pi i}\int_{\gamma_+}e^{-iR\cos(z)}\Upsilon(z)\text{d}z=0,
\end{align}
\mylinenum{then,}
\begin{align}
\label{MalThmpart1}\Upsilon(z)\equiv0&\ \ \hbox{if}\ \ D<1,\\
\label{MalThmpart2}\Upsilon(z)=\sin(z)\left[\sum_{j=0}^{d-1}c_j(\cos(z))^j\right]&\ \ \hbox{if}\ \ D\geq1,
\end{align}
\mylinenum{where $d$ is the integer part of $D$ and the coefficients $c_j$ are constants.}
\end{theorem}

\mylinenum{With Malyuzhinets' theorem, we have all the tools required to determine the \RED{spectral function} $s(z)$. Recall that it was assumed to be meromorphic in the domain \eqref{SMTwidehalfstrip}, \RED{moreover it has only one simple pole with unit residue at $z=\theta_\text{I}$ within the strip $|\RE{z}|<\theta_{\text{w}}$ correponding to the incident wave $\Phi_\text{I}$.} 

\subsection{Dirichlet boundary condition}\label{SMT-DBC}
Applying the Dirichlet \RED{BCs} \eqref{Intro-DBC} to the general solution \eqref{SMTSommInt} implies that}
\begin{align}\label{SMTSommIntDBC}
\frac{1}{2\pi i}\int_{\gamma_+}e^{-ikr\cos(z)}\left[s(\pm\theta_{\text{w}}+z)-s(\pm\theta_{\text{w}}-z)\right]\text{d}z=0.
\end{align}
\mylinenum{Due to \eqref{SMTsasym}, we can apply Malyuzhinets' theorem to \eqref{SMTSommIntDBC} to produce a pair of functional equations for $s(z)$,}
\begin{align}\label{SMTDBCequ}
s(\pm\theta_{\text{w}}+z)-s(\pm\theta_{\text{w}}-z)&=0.
\end{align}
\mylinenum{These equations imply that $s(z)$ is symmetric about $z=\pm\theta_{\text{w}}$ and as a consequence is $4\theta_{\text{w}}$ periodic. Because the pole at $z=\theta_\text{I}$ produces the incident wave, its residue should be 1, i.~e.}
$$\lim_{z\rightarrow\theta_\text{I}}\left[(z-\theta_\text{I})s(z)\right]=1.$$
\mylinenum{Because this pole is the only singularity in the strip $|\RE{z}|\leq\theta_{\text{w}}$, then by the determined symmetry, we also have poles at $z=2\theta_{\text{w}}-\theta_\text{I}$ and $-2\theta_{\text{w}}-\theta_\text{I}$ with residue $-1$. These two poles correspond to the top and bottom reflected waves respectively. The periodicity implies that each of the poles are repeated every $4\theta_{\text{w}}$ with the same residue. We can therefore express $s(z)$ as a sum of poles,}
\begin{align}\label{SMTSommFuncDpre}
s(z)=\sum_{n=-\infty}^{\infty}\frac{1}{z-\theta_\text{I}-4\theta_{\text{w}}n}-\sum_{n=-\infty}^{\infty}\frac{1}{z+\theta_\text{I}-2\theta_{\text{w}}-4\theta_{\text{w}}n}.
\end{align}
\mylinenum{Using the definition $\delta=\frac{\pi}{2\theta_{\text{w}}}$ and pole expansion of $\cot(z)$,}
$$\cot(z)=\sum_{n=-\infty}^{\infty}\frac{1}{z-n\pi},$$
\mylinenum{we rewrite $s(z)$ as follows,}
\begin{align}\label{SMTSommFuncDirichlet}
s(z)&=\frac{\delta}{2}\left[\cot\left(\frac{(z-\theta_\text{I})\delta}{2}\right)-\cot\left(\frac{(z-2\theta_\text{w}+\theta_\text{I})\delta}{2}\right)\right]=\frac{\delta\cos(\delta\theta_\text{I})}{\sin(\delta z)-\sin(\delta\theta_\text{I})}.
\end{align}

\mylinenum{It is easy to check that \eqref{SMTSommFuncDirichlet} satisfies the functional equations \eqref{SMTDBCequ}, satisfies the estimate $O(e^{-\delta|\IM{z}|})$ as $|\IM{z}|\rightarrow\infty$, and has a single pole with unit residue in the strip $|\RE{z}|\leq\theta_{\text{w}}$. This means that the solution to the exterior wedge problem with Dirichlet \RED{BCs} is,}
\begin{align}\label{SMTSommIntDirichlet}
\Phi(r,\theta)=\frac{1}{2\pi i}\int_{\gamma_++\gamma_-}e^{-ikr\cos(z)}\frac{\delta\cos(\delta\theta_\text{I})}{\sin(\delta(\theta+z))-\sin(\delta\theta_\text{I})}\,\text{d}z.
\end{align}
\mylinenum{

\subsection{Neumann boundary condition}\label{SMT-NBC}
Solving for the Neumann case is done in a very similar way. Applying the Neumann \RED{BCs} \eqref{Intro-NBC} to the general solution \eqref{SMTSommInt} implies that}
\begin{align}\label{SMTSommIntNBC}
\frac{1}{2\pi i}\int_{\gamma_+}e^{-ikr\cos(z)}\left[s'(\pm\theta_{\text{w}}+z)-s'(\pm\theta_{\text{w}}-z)\right]\text{d}z=0.
\end{align}
\mylinenum{We integrate by parts and apply the Malyuzhinets' theorem to \eqref{SMTSommIntNBC} to obtain the following pair of functional equations for $s(z)$,}
\begin{align}\label{SMTNBCequpre}
s(\pm\theta_{\text{w}}+z)+s(\pm\theta_{\text{w}}-z)&=c_\pm.
\end{align}
\mylinenum{Applying \eqref{SMTsasym} determines that $c_\pm=0$. Then the functional equations become,}
\begin{align}\label{SMTNBCequ}
s(\pm\theta_{\text{w}}+z)+s(\pm\theta_{\text{w}}-z)=0.
\end{align}
\mylinenum{These equations imply that $s(z)$ is antisymmetric about $z=\pm\theta_{\text{w}}$. This symmetry also shows that $s(z)$ is a $4\theta_{\text{w}}$ periodic function. \RED{Because} $\theta_\text{I}$ is the only pole of $s(z)$ in the strip $|\RE{z}|\leq\theta_{\text{w}}$ and has unit residue, then due to the antisymmetry and periodicity of $s(z)$, there are more poles located at $z=\theta_\text{I}+4\theta_{\text{w}}n$ and $z=2\theta_{\text{w}}-\theta_\text{I}+4\theta_{\text{w}}n$ for $n\in\mathbb{Z}$, all with unit residue. We can therefore express $s(z)$ as a sum of poles,}
\begin{align}
\nonumber s(z)&=\sum_{n=-\infty}^{\infty}\frac{1}{z-\theta_\text{I}-4\theta_{\text{w}}n}+\sum_{n=-\infty}^{\infty}\frac{1}{z+\theta_\text{I}-2\theta_{\text{w}}-4\theta_{\text{w}}n}\\
&=\frac{\delta}{2}\left[\cot\left(\frac{(z-\theta_\text{I})\delta}{2}\right)+\cot\left(\frac{(z-2\theta_{\text{w}}+\theta_\text{I})\delta}{2}\right)\right]=\frac{\delta\cos(\delta z)}{\sin(\delta z)-\sin(\delta\theta_\text{I})} .
\label{SMTSommFuncNeumann}
\end{align}

\mylinenum{It is easy to check that \eqref{SMTSommFuncNeumann} satisfies the functional equations \eqref{SMTNBCequ}, satisfies the estimate \eqref{SMTsasym} as $|\IM{z}|\rightarrow\infty$, and has one simple pole with unit residue in the strip $|\RE{z}|\leq\theta_{\text{w}}$. Finally, the solution to the exterior wedge problem with Neumann \RED{BCs} is,}
\begin{align}\label{SMTSommIntNeumann}
\Phi(r,\theta)=\frac{1}{2\pi i}\int_{\gamma_++\gamma_-}e^{-ikr\cos(z)}\frac{\delta\cos(\delta(\theta+z))}{\sin(\delta(\theta+z))-\sin(\delta\theta_\text{I})}\,\text{d}z.
\end{align}

\mylinenum{\RED{
\paragraph{Critical analysis} The main advantage of the S-M technique is its ease of implementation. Once in the right form, a natural deformation to the steepest descent contour can transform it into a very simple and fast converging integral (see Section \ref{SolAnalysis}). This aspect clearly makes it the gold standard representation of the general solution to the perfect wedge problem. Moreover, because of the form of the integrand, one can somehow think of the solution as a weighted superposition of plane waves.

An other advantage of such method is its flexibility. It is indeed possible to use it for more complicated problems such as the impedance wedge \citep{Malyuzhinets1958-3} or with various types of incident waves \citep{Bowman1987}. It has also been used in the quarter-plane problem to infer some results about the far-field structure \citep{Lyalinov2013}.

A disadvantage of the method is that it is not particularly constructive. In fact, it often starts with a form of the solution being posed, which, may somehow appear as some kind of ``black magic''. The \ref{app:app1} addresses this issue by showing that the two contour representation of the solution comes naturally from Green's identity. Moreover, once the integral form of the solution is written down, it is not straightforward to prove that the radiation condition is satisfied.
}


\section{The Wiener-Hopf technique}\label{WHT} 
The second method to be reviewed is the Wiener-Hopf (W-H) technique. Before authors such as Shanin and Daniele used the W-H technique for wedge problems, it was mostly used for waveguide problems or more \RED{complicated} half-plane problems. In the two papers \citep{AVShanin1996,AVShanin1998}, Shanin looks at solving wedge problems with inhomogeneous impedance \RED{BCs} via the W-H technique. In a series of papers and internal reports \citep{Daniele2000,Daniele2001,Daniele2003-2,Daniele2003}, Daniele develops several aspects of his method for impenetrable wedge problems.

In this section, we will combine elements from both \citep{AVShanin1996} and \citep{Daniele2003} to rederive \eqref{SMTSommIntDirichlet} and \eqref{SMTSommIntNeumann} using the W-H technique. The idea is to Laplace transform $\Phi$ and its $\theta$ derivative on the two wedge faces, $\theta=\pm\theta_{\text{w}}$, and the line of symmetry $\theta=0$. These transforms are used to produce the W-H equations. After the \RED{BCs} are considered, a mapping to a new complex plane is introduced so that the W-H technique can be applied to produce a solution that will match \eqref{SMTSommIntDirichlet} and \eqref{SMTSommIntNeumann}. We will start by studying this mapping.

\subsection{A \RED{useful} mapping}\label{WHM-mapping-section}
As we will see later, when the W-H equations are derived, they cannot be easily factorised using standard methods. To counter this issue, we will need to map these equations onto a new complex plane so that they can be reduced to classical W-H equations like those in \citep{Noble1958}. In order to do that, Shanin and Daniele use slightly different mappings,}
\begin{align}
\nonumber\hbox{Shanin's mapping:}\quad\eta(\alpha)=&\cos\left(\frac{2\theta_\text{w}}{\pi}\cos^{-1}(\sqrt{\alpha})\right)=\cos\left(\frac{\theta_\text{w}}{\pi}\cos^{-1}(2\alpha-1)\right),&\\
\nonumber\hbox{Daniele's mapping:\footnotemark}\quad\eta(\alpha)=&\ k\cos\left(\frac{\theta_\text{w}}{\pi}\cos^{-1}\left(\frac{\alpha}{k}\right)\right),
\end{align}
\mylinenum{\footnotetext{In Daniele's papers, the mapping has $-k$ in place of $k$ because he uses $e^{i\omega t}$ as time factor.}where $\eta$ and $\alpha$ are the old and new complex variables. Though the two mappings are conceptually equivalent, we will study Daniele's mapping in what follows. Note though that in his work, Daniele has assumed that $k$ has a small imaginary part in order to have a strip of analyticity for the W-H equations. However this is not strictly necessary, and here we will consider $k\in\mathbb{R^+}$, essentially reducing the W-H problem to a Riemann-Hilbert problem (see e.g. \citep{Kisil2015}). This means that $k$ does not need to appear explicitly in the mapping, and we can simply use}
\begin{align}\label{WHM-mapping}
\eta(\alpha)=\cos\left(\frac{\theta_\text{w}}{\pi}\cos^{-1}(\alpha)\right),
\end{align}
\mylinenum{with the corresponding inverse,}
\begin{align}\label{WHM-invmapping}
\alpha(\eta)=\cos\left(\frac{\pi}{\theta_\text{w}}\cos^{-1}(\eta)\right).
\end{align}
\mylinenum{We also consider the intermediate mapping and corresponding inverse,}
\begin{align}\label{WHM-zmapping}
z(\alpha)=\frac{\theta_\text{w}}{\pi}\cos^{-1}(\alpha),\quad\alpha(z)=\cos\left(\frac{\pi}{\theta_\text{w}}z\right).
\end{align}
\mylinenum{The mapping (\ref{WHM-mapping}) has a single branch cut along the real line segment $\alpha\in(-\infty,-1]$ where the local argument of the chosen branch is $(-\pi,\pi]$. This is done by choosing the branch of the inverse cosine such that $\pi-\cos^{-1}(x)=\cos^{-1}(-x)$\RED{, which is standard for programs such as Mathematica and MATLAB.} Note that the intermediate mapping limits $z$ to belong to the strip $\RE{z}\in[0,\theta_\text{w}]$.

One of the most important features of the mapping (\ref{WHM-mapping}) is that the upper half plane (UHP) $\IM{\alpha}\geq0$ is mapped to a subset of the UHP $\IM{\eta}\geq0$ as shown in \RED{Figures \ref{WHM-char-eta-f2-a} and \ref{WHM-char-eta-f2-b}.} This implies in particular that if a function $f(\eta)$ is analytic in the $\eta$ UHP, then the function $f(\eta(\alpha))$ is analytic in the $\alpha$ UHP. Another \RED{noteworthy} property is that if a function \RED{$g(\alpha)$} has no branch point at $\alpha=1$ (or $-1$), then the function \RED{$g(\alpha(z))$} will be symmetric with respect to $z=0$ (or $z=\theta_\text{w}$) in the $z$-plane. 

\RED{The mapping (\ref{WHM-mapping}) is designed specifically for manipulation of the following functions,}}
\RED{\begin{align}
\nonumber f_1(\eta)=\sqrt{1-\eta^2},\quad &f_2(\eta)=\eta\cos(\theta_\text{w})+\sqrt{1-\eta^2}\sin(\theta_\text{w}),\\ 
\label{WHM-Functions}f_3(\eta)=&\ \eta\sin(\theta_\text{w})-\sqrt{1-\eta^2}\cos(\theta_\text{w}),
\end{align}}
\mylinenum{\RED{where the branch for the square root is chosen such that $f_1(0)=1$.} For context, $f_2$ is used to identify domains of analyticity whereas, $f_1$ and $f_3$ are kernel functions that need to be factorised. \RED{Noting that $\eta$ and $z$ have the relation $\eta=\cos(z)$, we} map $f_1$ to the $z$ and $\alpha$ planes.}
\begin{align}\label{WHM-f1-mapped}
f_1(\cos(z))=\sin(z),\quad f_1(\eta(\alpha))=\sin\left(\frac{\theta_\text{w}}{\pi}\cos^{-1}(\alpha)\right).
\end{align}
\mylinenum{In the $\alpha$ plane, $f_1(\eta(\alpha))$ has branch cuts on the segments $\alpha\in(-\infty,-1]$ and $\alpha\in[1,\infty)$ \RED{(see Figure \ref{WHM-fig-f1-mapped} which is a phase plot of $f_1$ on the $\alpha$-plane).} 

Similarly we can map $f_2$ in the $z$ and $\alpha$ planes.}
\begin{align}\label{WHM-f2-mapped}
f_2(\cos(z))=\cos(\theta_\text{w}-z),\quad f_2(\eta(\alpha))=\cos\left(\frac{\theta_\text{w}}{\pi}\cos^{-1}(-\alpha)\right)=\eta(-\alpha).
\end{align}
\mylinenum{In the $\alpha$ plane, $f_2(\eta(\alpha))$ has a single branch cut on the segment $\alpha\in[1,\infty)$ \RED{and is related to $\eta(\alpha)$ by the identity $f_2(\eta(\alpha))=\eta(-\alpha)$, as shown by comparing the phase portraits in Figures \ref{WHM-eta-mapped} and \ref{WHM-fig-f2-mapped}}. The consequence is that the lower half plane (LHP) $\IM{\alpha}\leq0$ is mapped to a subset of the region $\IM{f_2(\eta)}\geq0$ as can be seen from comparing \RED{Figures \ref{WHM-char-eta-f2-c} and \ref{WHM-char-eta-f2-d}}. This also implies, in particular, that if a function $f(f_2)$ is analytic in the $f_2$ UHP, then the function $f(f_2(\eta(\alpha)))$ is analytic in the $\alpha$ LHP.

Lastly, we study the effect of the mapping on $f_3$,}
\begin{align}\label{WHM-f3-mapped}
f_3(\cos(z))=\sin(\theta_\text{w}-z),&\quad f_3(\eta(\alpha))=\sin\left(\frac{\theta_\text{w}}{\pi}\cos^{-1}(-\alpha)\right)=f_1(\eta(-\alpha)),
\end{align}
\mylinenum{showing that $f_1$ and $f_3$ are closely related in the $\alpha$-plane, in the sense that $f_3(\eta(\alpha))$ and $f_1(\eta(\alpha))$ have the same branch cuts but their phase plots are the symmetric images of each other about $z=0$, as illustrated in \RED{Figures \ref{WHM-fig-f1-mapped} and \ref{WHM-fig-f3-mapped}}. This relationship means that factorising $f_3$ in the $\alpha$ plane is analogous to factorising $f_1$. 

\begin{figure}[h!]\centering
	\subfloat[\label{WHM-char-eta-f2-a}]{
	\includegraphics[width=0.24\textwidth]{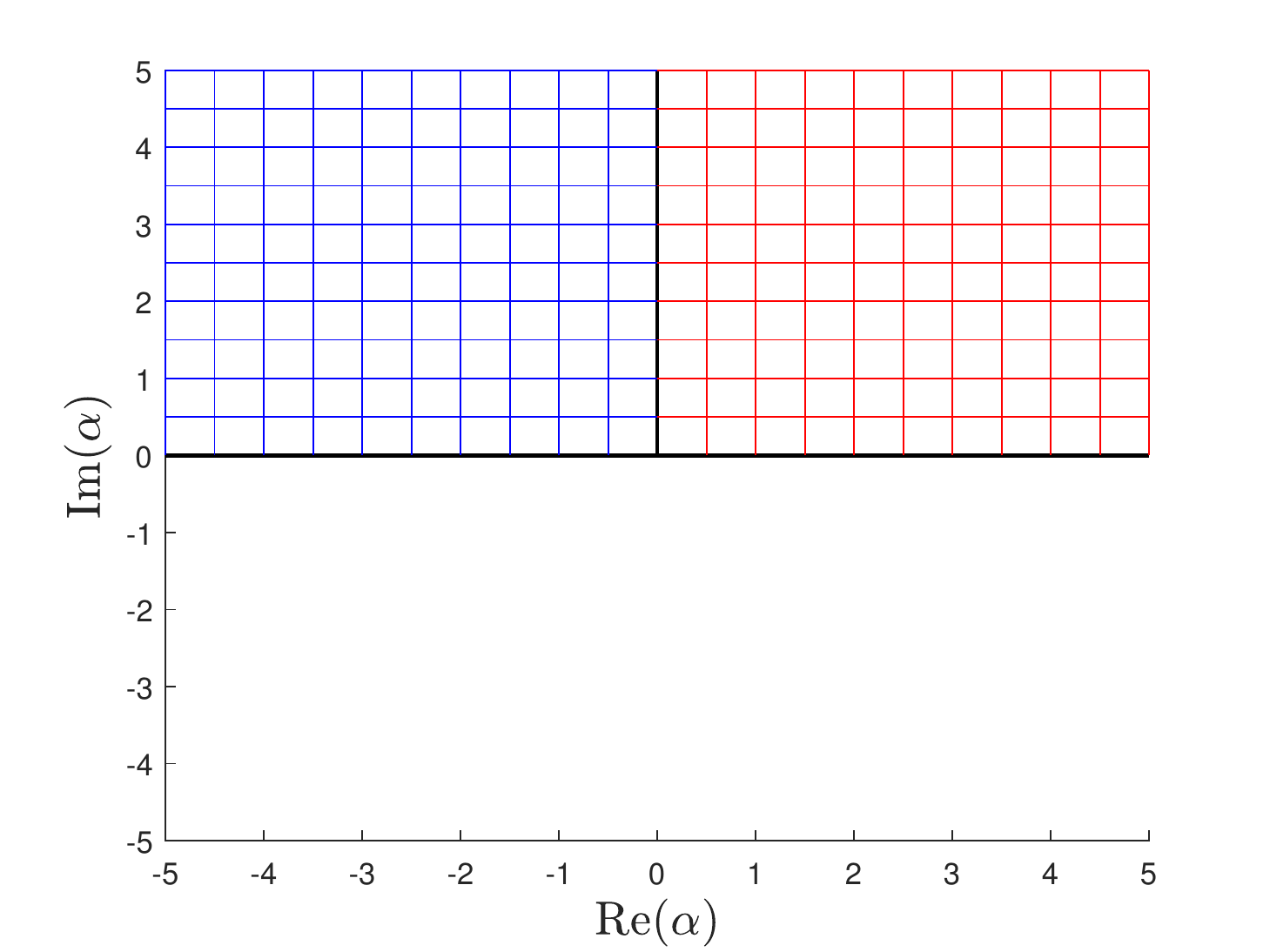}}
	\subfloat[\label{WHM-char-eta-f2-b}]{
	\includegraphics[width=0.24\textwidth]{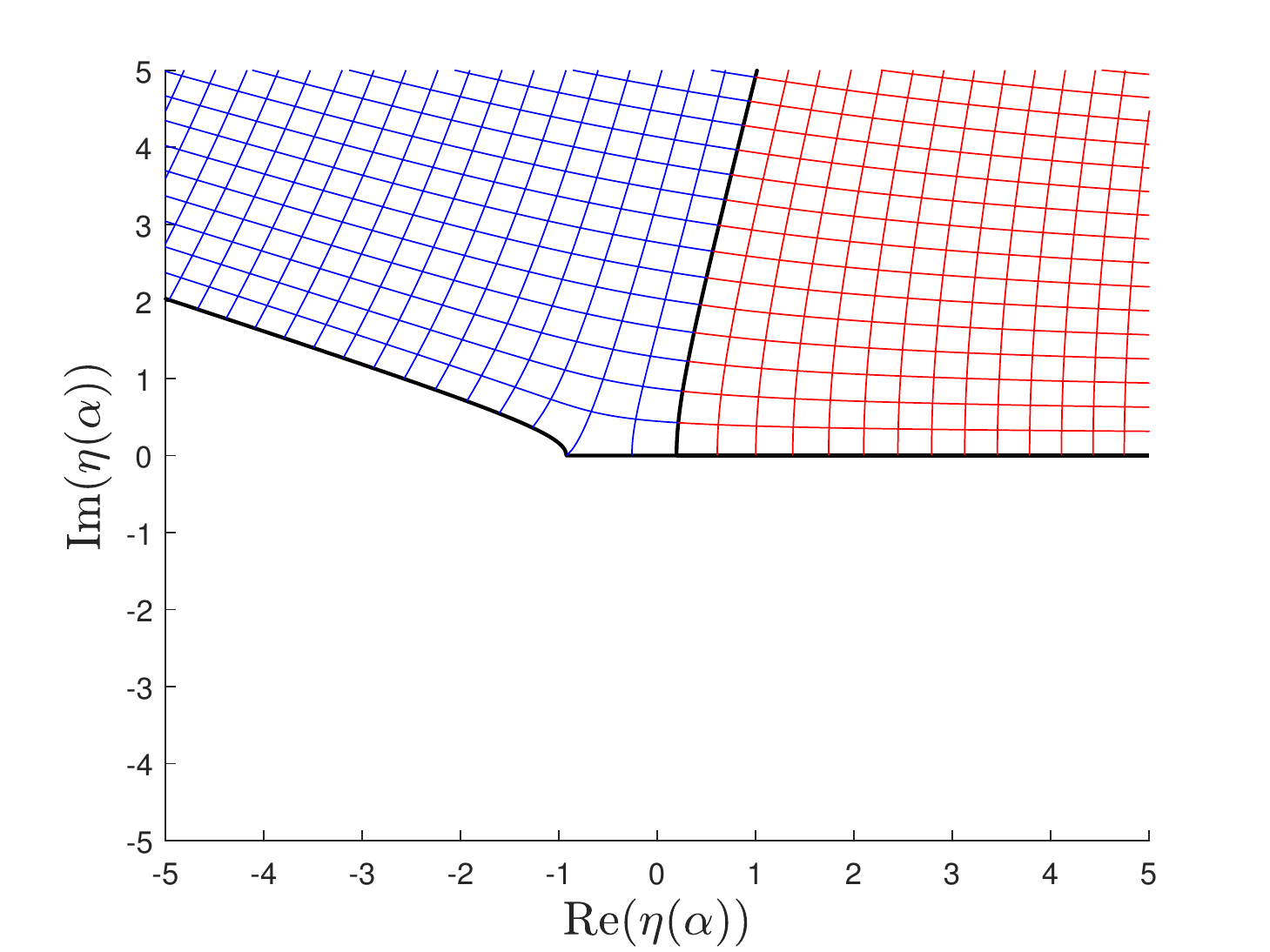}}
	\subfloat[\label{WHM-char-eta-f2-c}]{
	\includegraphics[width=0.24\textwidth]{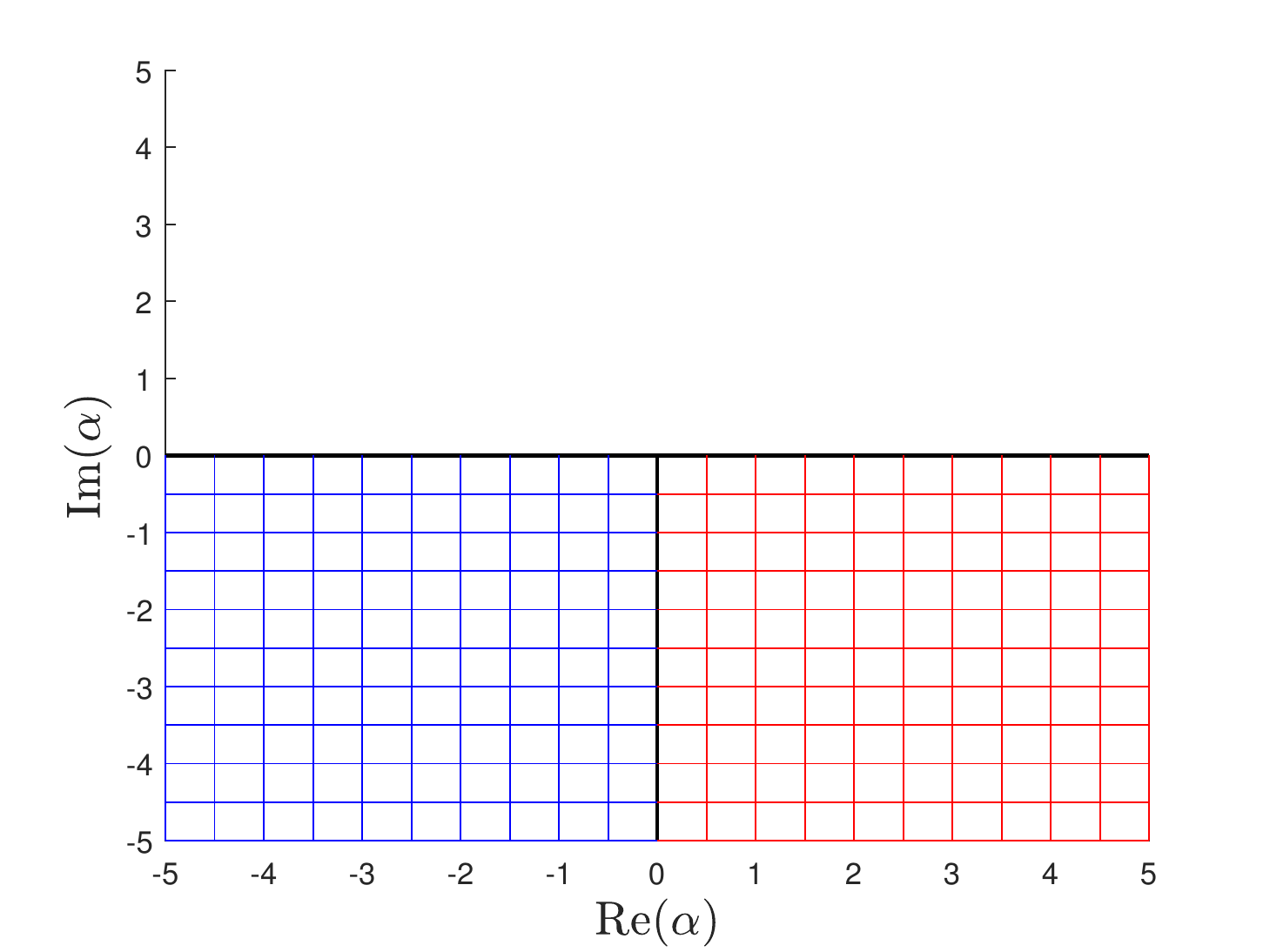}}
	\subfloat[\label{WHM-char-eta-f2-d}]{
	\includegraphics[width=0.24\textwidth]{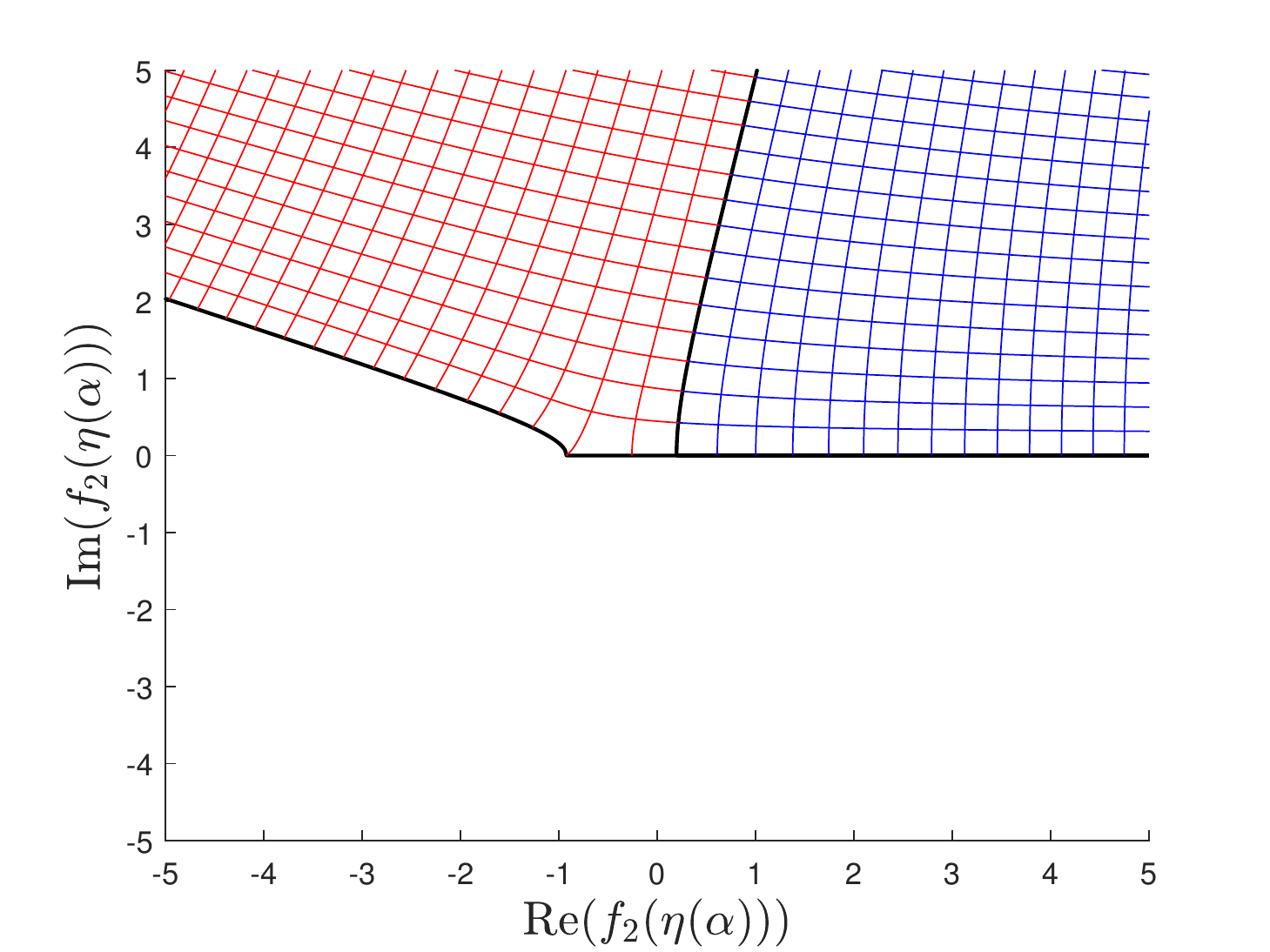}}
	\caption{\RED{With $\theta_{\text{w}}=7\pi/8$, these are pictures of characteristic domains showing that an upper half $\alpha$-plane (a) is mapped onto an upper half $\eta$-plane (b) and a lower half $\alpha$-plane (c) is mapped onto an upper half $f_2$-plane (d).}}
	\label{WHM-char-eta-f2}
\end{figure}

\begin{figure}[h!]\centering
	\subfloat[\label{WHM-eta-mapped}$\eta(\alpha)$]{
	\includegraphics[width=0.242\textwidth]{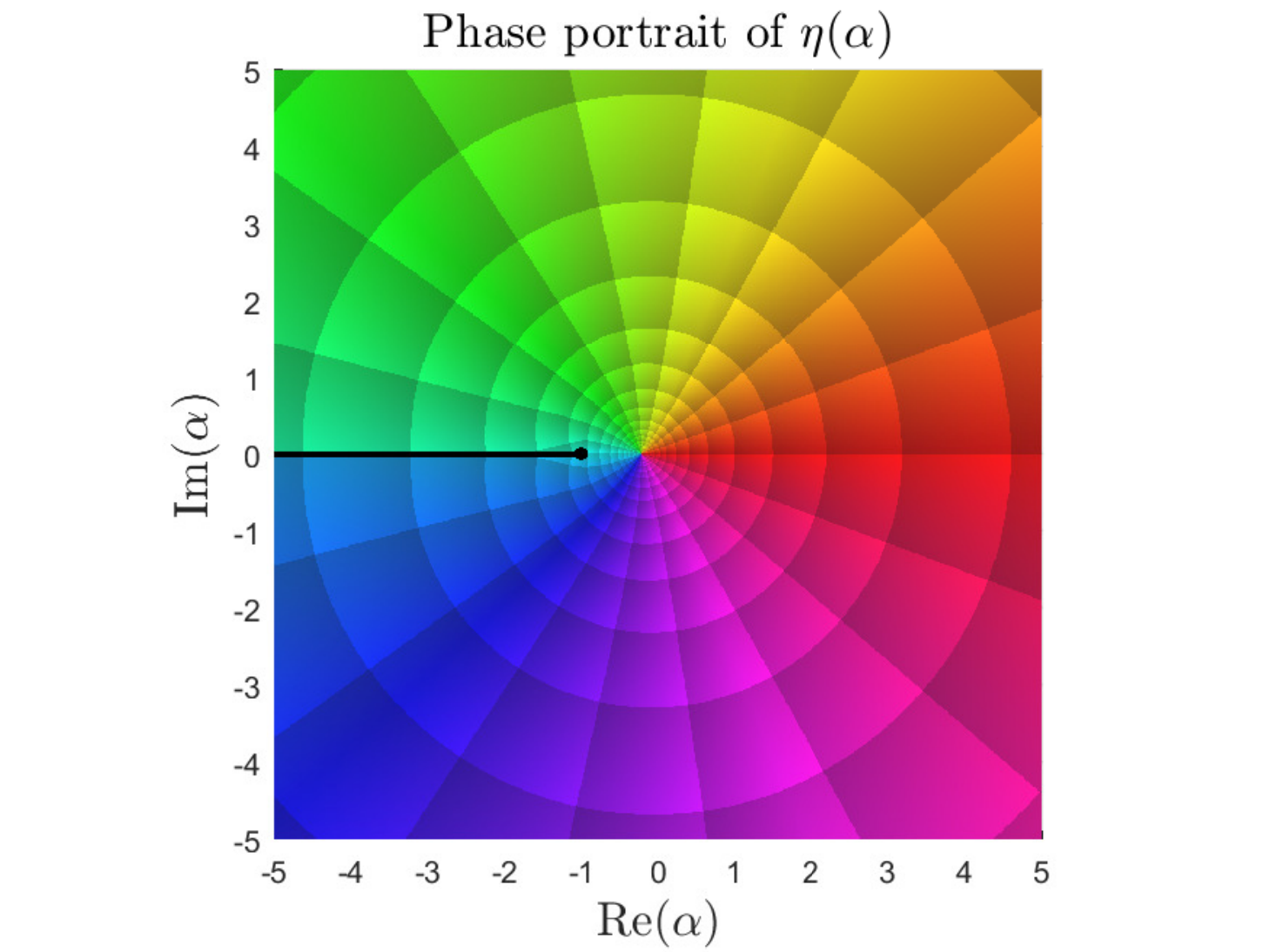}}
	\subfloat[\label{WHM-fig-f1-mapped}$f_1(\eta(\alpha))$]{
	\includegraphics[width=0.242\textwidth]{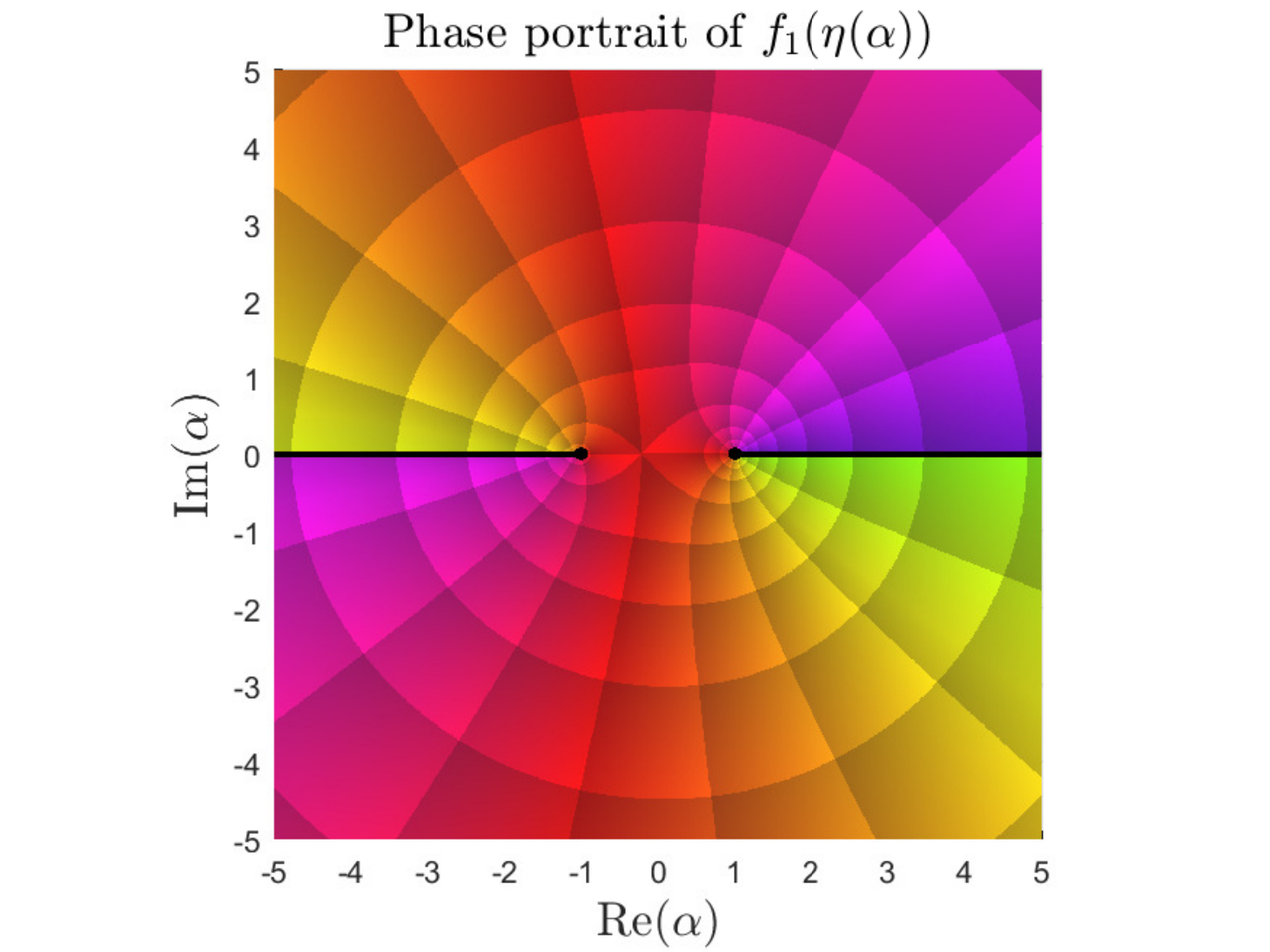}}
	\subfloat[\label{WHM-fig-f2-mapped}$f_2(\eta(\alpha))$]{
	\includegraphics[width=0.242\textwidth]{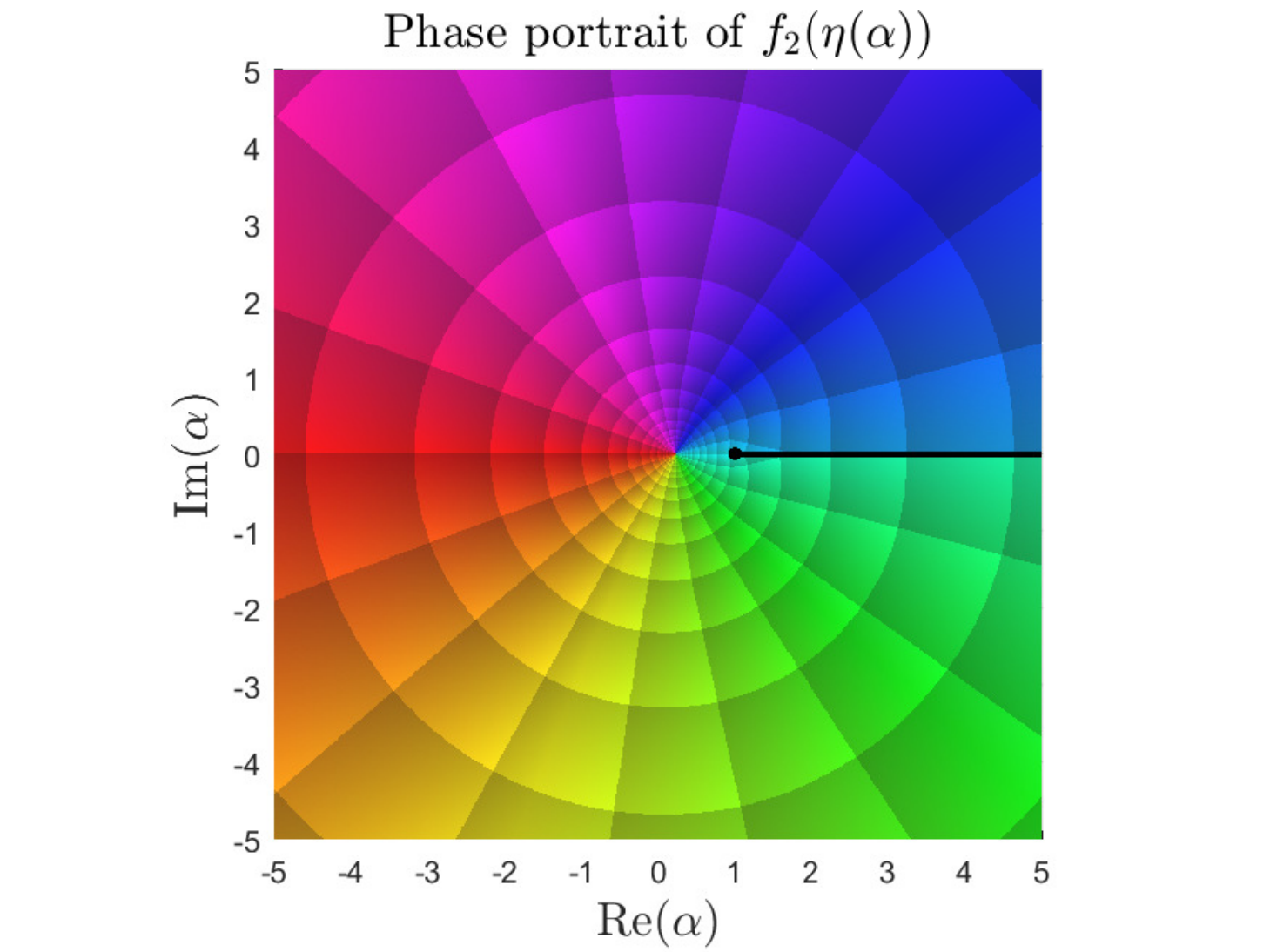}}
	\subfloat[\label{WHM-fig-f3-mapped}$f_3(\eta(\alpha))$]{
	\includegraphics[width=0.242\textwidth]{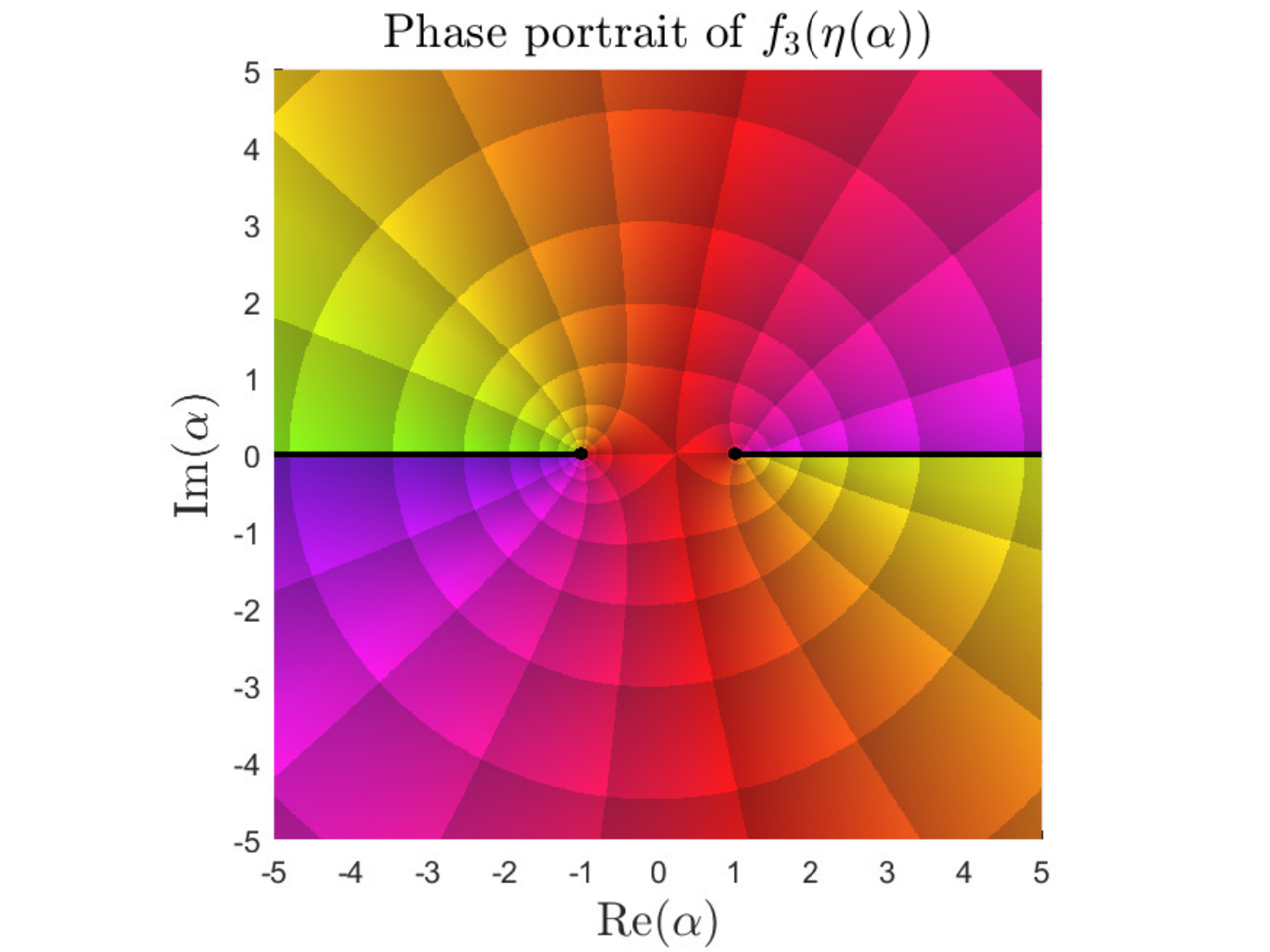}}
	\caption{\RED{Phase plots of $\eta$, $f_1$, $f_2$ and $f_3$ after being mapped onto the $\alpha$-plane. These illustrate any branch cuts and show the relationships between the four functions with $\theta_{\text{w}}=7\pi/8$}}
	\label{WHM-fig-mappedfunctions}
\end{figure}

Before performing such factorisation, let us define the domains $\mathcal{R}_\pm$ on which this factorisation will be done:}
\begin{eqnarray}
\mathcal{R}_+ &=& \{\alpha: \IM{\alpha}>0\}\cup\{\alpha: \RE{\alpha}>-1, \IM{\alpha}=0\} \label{WHM-D-UHP} \\
\mathcal{R}_- &=& \{\alpha: \IM{\alpha}<0\}\cup\{\alpha: \RE{\alpha}<1, \IM{\alpha}=0\} \label{WHM-D-LHP}
\end{eqnarray}
\mylinenum{Note that $\mathcal{R}_+\cup \mathcal{R}_-=\mathbb{C}$ and $\mathcal{R}_+\cap\mathcal{R}_-=(-1,1)$. Factorising $f_1$, requires to write}
\begin{align}\label{WHM-fac-anstaz}
f_1(\eta(\alpha))=f_{1+}(\alpha)f_{1-}(\alpha),
\end{align}
\mylinenum{where $f_{1+}$ is analytic in $\mathcal{R}_+$ and $f_{1-}$ is analytic in $\mathcal{R}_-$.
We expect that $f_{1+}$ will have a branch cut starting at $\alpha=-1$ and $f_{1-}$ a branch cut starting at $\alpha=1$. Hence $f_{1+}(\alpha(z))$ and $f_{1-}(\alpha(z))$ will be symmetric about the points $z=0$ and $\theta_\text{w}$ respectively. 
We can realise this factorisation by using the fact that both the leading order behaviour of $f_1$ as $\alpha\rightarrow1$ and the jump across the cut are consistent with the function $\sqrt{1-\alpha}$. This means that we can define $f_{1-}$ and $f_{1+}$ by,}
\begin{align}\label{WHM-f1-factorised-alpha}
f_{1-}(\alpha)=&\sqrt{\frac{1-\alpha}{2}},\quad f_{1+}(\alpha)=\frac{\sin\left(\frac{\theta_\text{w}}{\pi}\cos^{-1}(\alpha)\right)}{\sqrt{\frac{1-\alpha}{2}}}.
\end{align}

\begin{figure}[h!]\centering
	\subfloat[\label{WHM-(f1-)-part}$f_{1-}(\alpha)$]{
	\includegraphics[width=0.3\textwidth]{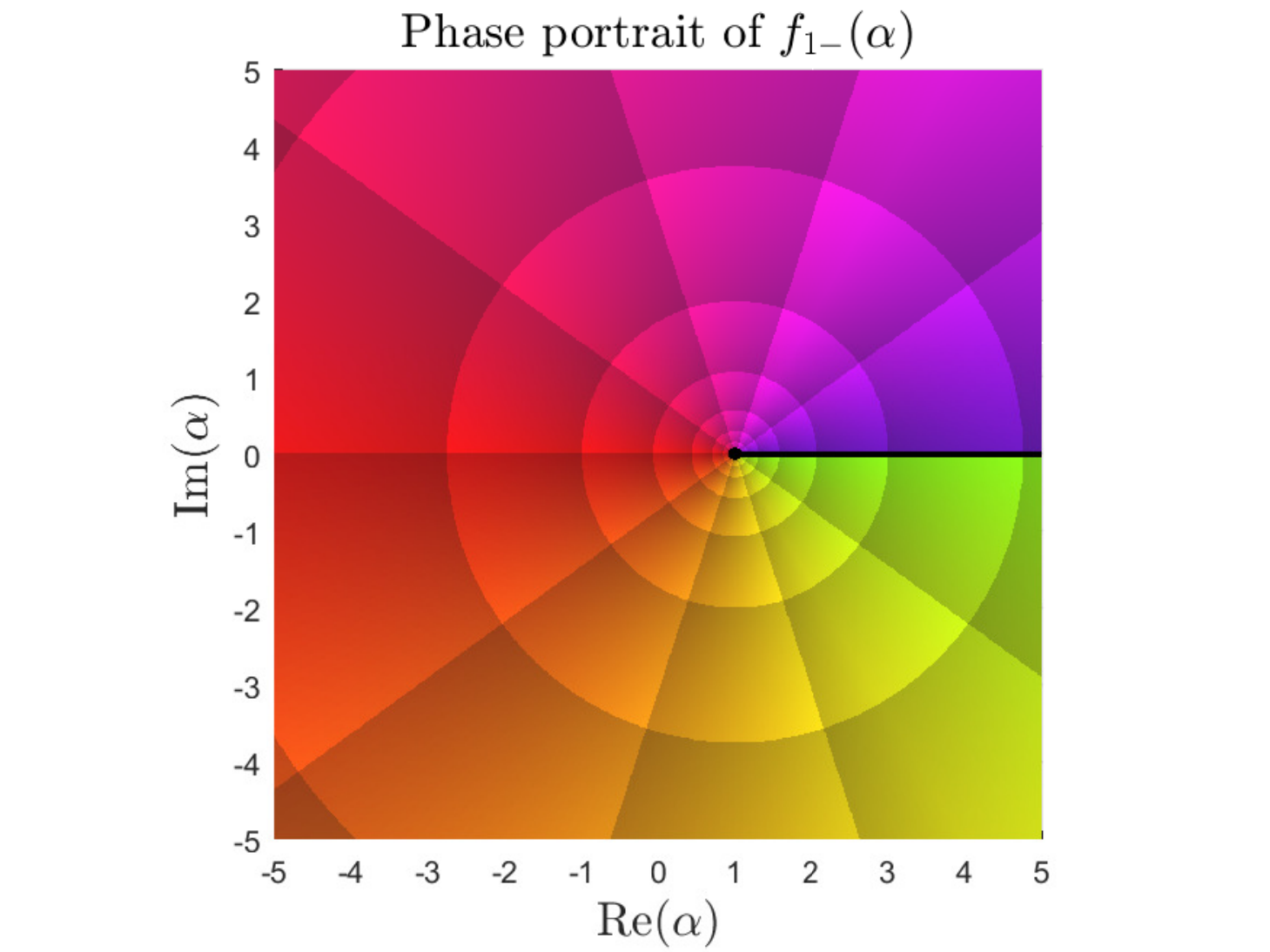}}
	\raisebox{16.7mm}{$\times$}
	\subfloat[\label{WHM-(f1+)-part}$f_{1+}(\alpha)$]{
	\includegraphics[width=0.3\textwidth]{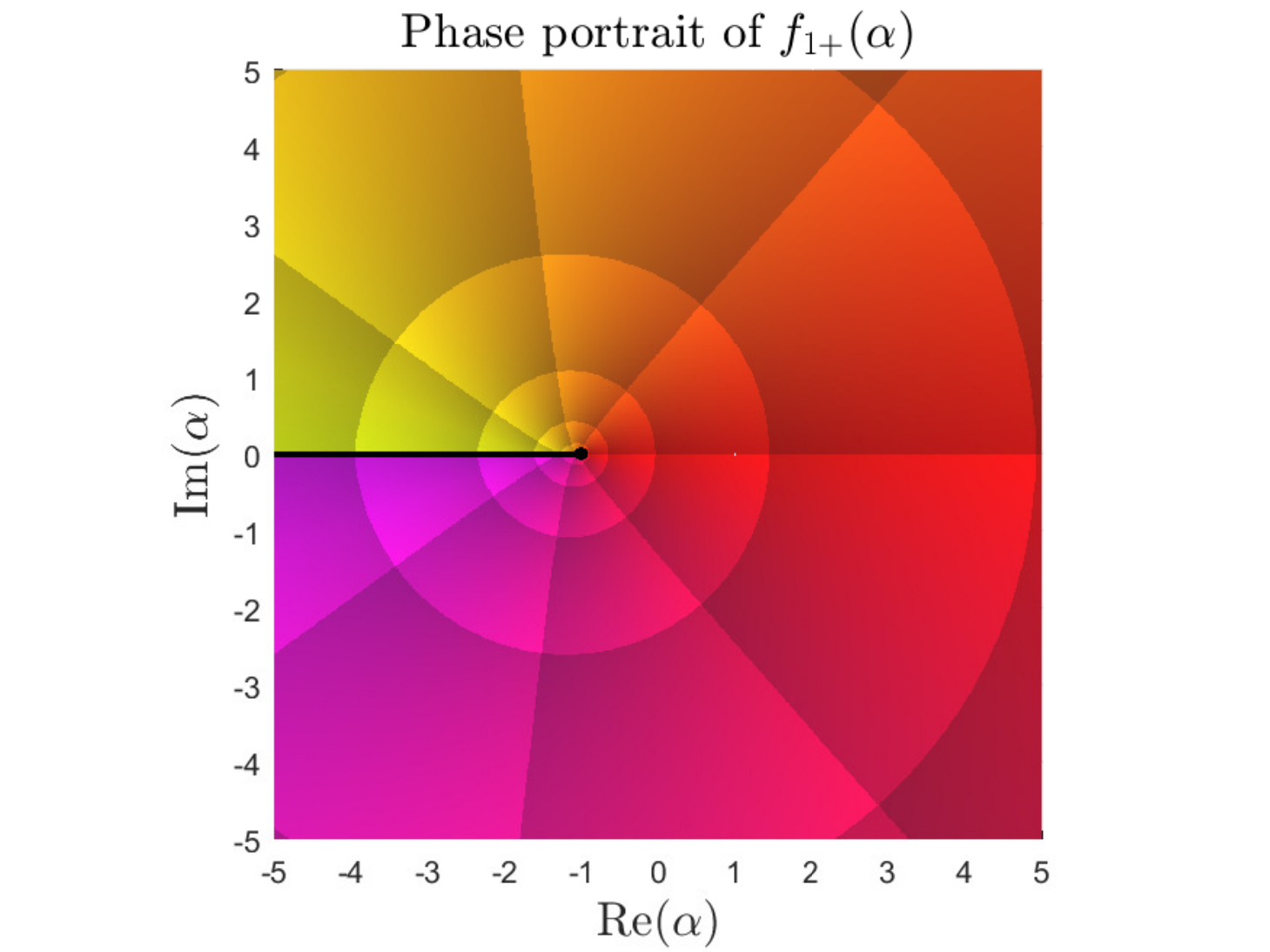}}
	\raisebox{16.7mm}{$=$}
	\subfloat[\label{WHM-f1-whole}$f_1(\eta(\alpha))$]{
	\includegraphics[width=0.3\textwidth]{WH_plots_MATLAB/f1-map-eps-converted-to.pdf}}
	\caption{\RED{Various phase plots helping to illustrate the factorisation \eqref{WHM-fac-anstaz} by displaying each of the parts of $f_1$ in the $\alpha$ plane for $\theta_{\text{w}}=7\pi/8$.}}
	\label{WHM-f1-fac}
\end{figure}

\mylinenum{Clearly $f_{1-}$ has a branch cut on the segment $\alpha\in[1,\infty)$ and is analytic at $\alpha=-1$. While $f_{1+}$ retains the branch cut on the segment $\alpha\in(-\infty,-1]$, and dividing by $\sqrt{1-\alpha}$ has made $\alpha=1$ become a removable singularity and cancelled out the cut discontinuity. Hence we can assign the limiting value to $\alpha=1$ and make $f_{1+}$ analytic at that point. Recalling the definition $\delta=\frac{\pi}{2\theta_{\text{w}}}$, we also map $f_{1-}$ and $f_{1+}$ to the $z$-plane as follows:}
\begin{align}\label{WHM-f1-factorised-z}
f_{1-}(\alpha(z))=&\sin(\delta z),\quad f_{1+}(\alpha(z))=\frac{\sin(z)}{\sin(\delta z)}.
\end{align}
\mylinenum{As anticipated, due to the absence of the branch point at $\alpha=-1$, $f_{1-}(\alpha(z))$ is symmetric about $z=\theta_\text{w}$. Similarly, due to the absence of the branch point at $\alpha=1$, $f_{1+}(\alpha(z))$ is symmetric about $z=0$. \RED{Figure \ref{WHM-f1-fac} illustrates this factorisation by the relationship \eqref{WHM-fac-anstaz}}.}

\mylinenum{\RED{Now we factorise $f_3$ using the established relation $f_3(\eta(\alpha))=f_1(\eta(-\alpha))$ (see \eqref{WHM-f3-mapped}) to find that $f_{3+}(\alpha)=f_{1-}(-\alpha)$ and $f_{3-}(\alpha)=f_{1+}(-\alpha)$. It can be shown that} $f_{3+}$ has a branch cut on the segment $\alpha\in(-\infty,-1]$, is analytic at $\alpha=1$ and $f_{3+}(\alpha(z))$ is symmetric about $z=0$, while $f_{3-}$ has a branch cut on the segment $\alpha\in[1,\infty)$, can be made analytic at $\alpha=-1$ and $f_{3-}(\alpha(z))$ is symmetric about $z=\theta_\text{w}$. 


\subsection{Derivation of the Wiener-Hopf equations}

\RED{Daniele derives the W-H equations by rewriting the Helmholtz equation in terms of an oblique Cartesian coordinate system and uses a Laplace transformation in each of the new coordinates \citep{Daniele2003}.} However this process can be time-consuming and we will show a different method here.


\citep{AVShanin1996} tackles an interior wedge problem with inhomogeneous impedance \RED{BCs}. In this geometry, only one W-H equation is derived, which is obtained via Green's second identity. However we will need to split the exterior wedge region into two halves to obtain two W-H equations. Take Green's second identity for functions $u,v$, twice continuously differentiable on domain $\Omega\in\mathbb{R}^2$ with boundary $\partial \Omega$,}
\begin{align}\label{WHM-Greens2ndID}
\RED{\int_\Omega \left(v\nabla^2u-u\nabla^2v\,\right)\text{d}\Omega=\int_{\partial \Omega}\left(v\ParDer{u}{\boldsymbol{n}}-u\ParDer{v}{\boldsymbol{n}}\,\right)\text{d}S,}
\end{align}
\mylinenum{where $\ParDer{}{\boldsymbol{n}}$ is the normal derivative. Here $u$ is the unknown solution $\Phi$ and we choose a suitable test function for $v$ that satisfies \RED{the} Helmholtz equation \RED{\eqref{Intro-Helmholtz}}. Then the left hand side of \eqref{WHM-Greens2ndID} is automatically zero. We do this for two wedge regions $\theta\in[0,\theta_\text{w}]$ and $\theta\in[-\theta_\text{w},0]$ which require different test functions. The right hand side of \eqref{WHM-Greens2ndID} has two parts, the wedge boundary at $\theta=\pm\theta_\text{w}$ and an imaginary boundary at $\theta=0$. \RED{For the upper region, we choose the test function $v=e^{ikr\cos(\theta-z)}$, leading \eqref{WHM-Greens2ndID} to become}}
\begin{align}
\int_0^\infty&\left[\frac{1}{ikr}\left.\ParDer{\Phi}{\theta}\right|_{\theta=0}\!\!\!\!-\sin(z)\left.\Phi\right|_{\theta=0}\right]e^{ikr\cos(z)}\text{d}r\nonumber\\
=&\int_0^\infty\left[\frac{1}{ikr}\left.\ParDer{\Phi}{\theta}\right|_{\theta=\theta_\text{w}}\!\!\!\!\!\!+\sin(\theta_\text{w}-z)\left.\Phi\right|_{\theta=\theta_\text{w}}\right]e^{ikr\cos(\theta_\text{w}-z)}\text{d}r.
\label{WHM-Greens2ndID1eqn}\end{align}
\mylinenum{\RED{For the lower region, $\theta\in[-\theta_\text{w},0]$, we choose the slightly modified test function $v=e^{ikr\cos(\theta+z)},$ leading \eqref{WHM-Greens2ndID} to become}}
\begin{align}
\int_0^\infty&\left[\frac{1}{ikr}\left.\ParDer{\Phi}{\theta}\right|_{\theta=0}\!\!\!\!+\sin(z)\left.\Phi\right|_{\theta=0}\right]e^{ikr\cos(z)}\text{d}r\nonumber\\
=&\int_0^\infty\left[\frac{1}{ikr}\left.\ParDer{\Phi}{\theta}\right|_{\theta=-\theta_\text{w}}\!\!\!\!\!\!\!\!-\sin(\theta_\text{w}-z)\left.\Phi\right|_{\theta=-\theta_\text{w}}\right]e^{ikr\cos(\theta_\text{w}-z)}\text{d}r.
\label{WHM-Greens2ndID2eqn}\end{align}
\mylinenum{Define the Laplace transform with the following inverse,}
\begin{align}\label{WHM-LapTrans}
F(\eta)=\int_0^\infty f(r)e^{ikr\eta}\text{d}r,\quad f(r)=\frac{k}{2\pi}\int_{-\infty}^{\infty}F(\eta)e^{-ikr\eta}\text{d}\eta,
\end{align}
\mylinenum{where $F(\eta)$ is analytic in the half-plane $\IM{\eta}>0$, then we define the transforms of $\Phi$ and $\frac{1}{ikr}\ParDer{\Phi}{\theta}$,}
\begin{align}\label{WHM-ULapTrans}
U(\eta,\theta)=\int_0^\infty \Phi(r,\theta) e^{ikr\eta}\text{d}r,\quad V(\eta,\theta)=\int_0^\infty \frac{1}{ikr}\ParDer{\Phi}{\theta}(r,\theta)e^{ikr\eta}\text{d}r.
\end{align}
\mylinenum{These transforms are applied to both \eqref{WHM-Greens2ndID1eqn} and \eqref{WHM-Greens2ndID2eqn} which produces the W-H equations,}
\begin{align}
\nonumber V(\cos(z),0)&-\sin(z)U(\cos(z),0)\\
\nonumber&=V(\cos(\theta_\text{w}-z),\theta_\text{w})+\sin(\theta_\text{w}-z)U(\cos(\theta_\text{w}-z),\theta_\text{w}),\\
\nonumber V(\cos(z),0)&+\sin(z)U(\cos(z),0)\\
\label{WHM-preWHsystem}&=V(\cos(\theta_\text{w}-z),-\theta_\text{w})-\sin(\theta_\text{w}-z)U(\cos(\theta_\text{w}-z),-\theta_\text{w}).
\end{align}
\mylinenum{Adding and subtracting these two equations leads to,}
\begin{align}
\nonumber 2V(\cos(z),0)=\ &V(\cos(\theta_\text{w}-z),-\theta_\text{w})+V(\cos(\theta_\text{w}-z),\theta_\text{w})\\
\nonumber&+\sin(\theta_\text{w}-z)\left[U(\cos(\theta_\text{w}-z),\theta_\text{w})-U(\cos(\theta_\text{w}-z),-\theta_\text{w})\right],\\
\nonumber2\sin(z)U(\cos(z),0)=\ &V(\cos(\theta_\text{w}-z),-\theta_\text{w})-V(\cos(\theta_\text{w}-z),\theta_\text{w})\\
\label{WHM-WHsystem}&-\sin(\theta_\text{w}-z)\left[U(\cos(\theta_\text{w}-z),\theta_\text{w})+U(\cos(\theta_\text{w}-z),-\theta_\text{w})\right].
\end{align}
\mylinenum{These are the so-called generalised W-H equations. In this system, the functions \RED{$U(\cdot,\theta)$ and $V(\cdot,\theta)$ are analytic in a region containing the upper half plane} regardless of the value of $\theta$. We solve the system \eqref{WHM-WHsystem} for $U(\cos(z),0)$ and $V(\cos(z),0)$ using the boundary data on the right hand side. 
Noting that $\cos(\mathbb{R})$ is equal to $\gamma_+$ with the opposite orientation, the inverse transform of $U(\cos(z),\theta)$ is,}
\begin{align}
\Phi(r,\theta)=\frac{1}{2\pi i}\int_{\gamma_+}e^{-ikr\cos(z)}U(\cos(z),\theta)\left(ik\sin(z)\right)\text{d}z,
\end{align}
\mylinenum{which is clearly very similar to the Sommerfeld integral. Applying Malyuzhinets' theorem, we find that,}
\begin{align}\label{WHM-Somm}
ik\sin(z)U(\cos(z),\theta)=s(\theta+z)-s(\theta-z).
\end{align}
\mylinenum{\RED{We can derive a second formula by comparing the inverse transform of $V(\cos(z),\theta)$ with the following,}}
\begin{align}
\label{WHM-SommInt-der}\frac{1}{ikr}\ParDer{\Phi}{\theta}(r,\theta)
=&-\frac{1}{2\pii}\int_{\gamma_+}e^{-ikr\cos(z)}\sin(z)\left[s(\theta+z)+s(\theta-z)\right]\text{d}z,
\end{align}
\mylinenum{which is obtained by differentiating \eqref{SMTSommInt} with respect to $\theta$ and integration by parts. Applying Malyuzhinets' theorem, we find that,}
\begin{align}\label{WHM-Somm(2)}
-ikV(\cos(z),\theta)=s(\theta+z)+s(\theta-z).
\end{align}
\mylinenum{\RED{Adding \eqref{WHM-Somm} to \eqref{WHM-Somm(2)} and setting $\theta=0$ implies that for all $z$,}}
\RED{\begin{align}\label{WHM-Somm-Sol}
s(z)=\frac{ik}{2}\left[\sin(z)U(\cos(z),0)-V(\cos(z),0)\right],
\end{align}}
\mylinenum{\RED{establishing a link between the spectral function and the W-H unknowns. We will now} apply the \RED{BCs} to solve the W-H system \eqref{WHM-WHsystem}.

\subsection{Dirichlet boundary condition}
The transformed Dirichlet \RED{BCs} are $U(\RED{\cos(\theta_{\text{w}}-z)},\pm\theta_{\text{w}})=0$ which simplify \eqref{WHM-WHsystem} to,}
\begin{align}
2V(\cos(z),0)&=V(\cos(\theta_\text{w}-z),-\theta_\text{w})+V(\cos(\theta_\text{w}-z),\theta_\text{w}),\nonumber\\
2\sin(z)U(\cos(z),0)&=V(\cos(\theta_\text{w}-z),-\theta_\text{w})-V(\cos(\theta_\text{w}-z),\theta_\text{w}).\label{WHM-D}
\end{align}
\mylinenum{In this form the W-H technique cannot be applied, so \eqref{WHM-zmapping} and \eqref{WHM-f2-mapped} \RED{(discussed in Section \ref{WHM-mapping-section})} are used here,}
\begin{align}
\label{WHM-D(1)}2V(\cos(z(\alpha)),0)&=V(\cos(z(-\alpha)),-\theta_\text{w})+V(\cos(z(-\alpha)),\theta_\text{w}),\\
\label{WHM-D(2)}2\sin(z(\alpha))U(\cos(z(\alpha)),0)&=V(\cos(z(-\alpha)),-\theta_\text{w})-V(\cos(z(-\alpha)),\theta_\text{w}).
\end{align}
\mylinenum{In the $\alpha$-plane, $U(\cos(z(\alpha)),0)$ and $V(\cos(z(\alpha)),0)$ are analytic in $\mathcal{R}_+$, except for some potential poles on the real line segment \RED{$\mathcal{R}_+\cap\mathcal{R}_-$}. Similarly, $V(\cos(z(-\alpha)),\pm\theta_\text{w})$ are analytic in $\mathcal{R}_-$, except for some potential poles on \RED{$\mathcal{R}_+\cap\mathcal{R}_-$}. We have already factorised \RED{$\sin(z(\alpha))=f_1(\eta(\alpha))$} in \eqref{WHM-f1-factorised-alpha}, leading \eqref{WHM-D(2)} to become,}
\begin{align}\label{WHM-D(3)}
2f_{1+}(\alpha)U(\cos(z(\alpha)),0)&=\frac{1}{f_{1-}(\alpha)}\left[V(\cos(z(-\alpha)),-\theta_\text{w})-V(\cos(z(-\alpha)),\theta_\text{w})\right].
\end{align}
\mylinenum{For both equations, \eqref{WHM-D(1)} and \eqref{WHM-D(3)}, the left sides are meromorphic in $\mathcal{R}_+$ and the right sides are meromorphic in $\mathcal{R}_-$, however due to \RED{the} potential poles, these equations cannot be used to create an entire function. To counteract this we remove the poles on the right side using the knowledge of the GO component of the solution. Assuming that $\theta_\text{w}>\pi/2$, the GO components of $V(\cos(z(-\alpha)),\pm\theta_\text{w})$ are,}
\begin{align}
\nonumber V^{(\text{GO})}(\cos(z(-\alpha)),\theta_{\text{w}})&=\frac{2i\sin(\theta_\text{w}-\theta_\text{I})}{k\cos(z(-\alpha))-k\cos(\theta_\text{w}-\theta_\text{I})},\\
\nonumber V^{(\text{GO})}(\cos(z(-\alpha)),-\theta_{\text{w}})&=-\frac{2i\sin(\theta_\text{w}+\theta_\text{I})\mathcal{H}(\pi-\theta_{\text{w}}-\theta_\text{I})}{k\cos(z(-\alpha))-k\cos(\theta_{\text{w}}+\theta_\text{I})},
\end{align}
\mylinenum{where $\mathcal{H}$ is the Heaviside function. The two poles correspond to $z(-\alpha)=\pm(\theta_\text{w}-\theta_\text{I})$, i.e. to $z(\alpha)=\pm \theta_\text{I}$. However there is no $\alpha$ in the chosen branch of inverse cosine that satisfies $z(\alpha)=-\theta_\text{I}$. \RED{This means that the only pole that needs to be removed is that of $V^{(\text{GO})}(\cos(z(-\alpha)),\theta_{\text{w}})$ at $z(\alpha)=\theta_\text{I}$,} corresponding to $\alpha=\alpha_0=\cos\left(2\delta\theta_\text{I}\right)$. The residue at this pole is,}
\begin{align}\label{WHM-pole-resD}
\lim_{\alpha\rightarrow\alpha_0}\left[\left(\alpha-\alpha_0\right)V^{(\text{GO})}(\cos(z(-\alpha)),\theta_{\text{w}})\right]=\frac{4\delta\sin\left(2\delta\theta_\text{I}\right)}{ik}.
\end{align}
\mylinenum{With this residue, we remove the pole from equations \eqref{WHM-D(1)} and \eqref{WHM-D(3)},}
\begin{align}
\nonumber&V(\cos(z(\alpha)),0)-\frac{2\delta\sin\left(2\delta\theta_\text{I}\right)}{ik\left(\alpha-\cos\left(2\delta\theta_\text{I}\right)\right)}\\
\label{WHM-D(4)}&=\frac{1}{2}\left[V(\cos(z(-\alpha)),-\theta_\text{w})+V(\cos(z(-\alpha)),\theta_\text{w})\right]-\frac{2\delta\sin\left(2\delta\theta_\text{I}\right)}{ik\left(\alpha-\cos\left(2\delta\theta_\text{I}\right)\right)},\\
\nonumber &f_{1+}(\alpha)U(\cos(z(\alpha)),0)+\frac{2\delta\sin\left(2\delta\theta_\text{I}\right)}{ikf_{1-}(\alpha_0)\left(\alpha-\cos\left(2\delta\theta_\text{I}\right)\right)}\\
\label{WHM-D(5)}&=\frac{1}{2f_{1-}(\alpha)}\left[V(\cos(z(-\alpha)),-\theta_\text{w})-V(\cos(z(-\alpha)),\theta_\text{w})\right]+\frac{2\delta\sin\left(2\delta\theta_\text{I}\right)}{ikf_{1-}(\alpha_0)\left(\alpha-\cos\left(2\delta\theta_\text{I}\right)\right)}.\end{align}
\mylinenum{In both equations \eqref{WHM-D(4)} and \eqref{WHM-D(5)}, the left sides are now analytic in $\mathcal{R}_+$ and the right sides are analytic in $\mathcal{R}_-$. In order to apply Liouville's theorem, we must determine the behaviour of each part in equations \eqref{WHM-D(4)} and \eqref{WHM-D(5)} as $|\alpha|\rightarrow\infty$. The edge condition \eqref{Intro-Edge} for the Dirichlet case implies that $\Phi=O(r^\delta)$ and $\frac{1}{r}\ParDer{\Phi}{\theta}=O(r^{\delta-1})$. Using the well-known fact that for any function $f(r)$ behaving like $r^\delta$ as $r\rightarrow 0$, its Laplace transform $F(\eta)$, as defined by \eqref{WHM-LapTrans}, behaves like $|\eta|^{-\delta-1}$ as $|\eta|\rightarrow\infty$, and noting that $\eta(\alpha)=O(|\alpha|^\frac{\theta_\text{w}}{\pi})$ as $|\alpha|\rightarrow\infty$, we can show that as $|\alpha|\rightarrow\infty$, \RED{$U(\cos(z(\alpha)),0)=O(|\alpha|^{-\frac{1}{2}-\frac{\theta_\text{w}}{\pi}})$, and $V(\cos(z(\alpha)),0)=O(|\alpha|^{-\frac{1}{2}})$} within $\mathcal{R_+}$, \RED{while $V(\cos(z(-\alpha)),\pm\theta_\text{w})=O(|\alpha|^{-\frac{1}{2}})$} within $\mathcal{R_-}$. We can also determine that $f_{1-}(\alpha)=O(|\alpha|^{\frac{1}{2}})$ and $f_{1+}(\alpha)=O(|\alpha|^{\frac{\theta_\mathrm{w}}{\pi}-\frac{1}{2}})$ as $|\alpha|\rightarrow\infty$. This means that all parts of equations \eqref{WHM-D(4)} and \eqref{WHM-D(5)} are decaying as $|\alpha|\rightarrow\infty$ in the \RED{appropriate} half plane. Construct the two functions,}
\begin{align}\RED{
\Psi_1(\alpha)=\begin{cases}
\hbox{LHS}\eqref{WHM-D(4)}&\text{ in }\mathcal{R}_+,\\
\hbox{RHS}\eqref{WHM-D(4)}&\text{ in }\mathcal{R}_-,\\
\eqref{WHM-D(4)}&\text{ in }\mathcal{R}_+\cap\mathcal{R}_-,\end{cases}\ \Psi_2(\alpha)=\begin{cases}
\hbox{LHS}\eqref{WHM-D(5)}&\text{ in }\mathcal{R}_+,\\
\hbox{RHS}\eqref{WHM-D(5)}&\text{ in }\mathcal{R}_-,\\
\eqref{WHM-D(5)}&\text{ in }\mathcal{R}_+\cap\mathcal{R}_-.\end{cases}}
\end{align}
\mylinenum{Both $\Psi_1$ and $\Psi_2$ are entire and decaying at infinity, therefore Liouville's theorem can be applied to show that \RED{$\Psi_1,\Psi_2\equiv0$}. It is hence possible to determine $V(\cos(z),0)$ and $U(\cos(z),0)$,}
\begin{align}\label{WHM-D-Vsol}
V(\cos(z),0)&=\frac{2\delta\sin(2\delta\theta_\text{I})}{ik(\cos(2\delta z)-\cos(2\delta\theta_\text{I}))},\\
\label{WHM-D-Usol}\sin(z)U(\cos(z),0)&=-\frac{4\delta\cos(\delta\theta_\text{I})\sin(\delta z)}{ik(\cos(2\delta z)-\cos(2\delta\theta_\text{I}))}.
\end{align}
\mylinenum{Equations \eqref{WHM-D-Vsol} and \eqref{WHM-D-Usol} can be substituted into \eqref{WHM-Somm-Sol} to obtain,}
\begin{align}\label{WHMSommFuncDirichlet}
s(z)=\frac{ik}{2}\left[\sin(z)U(\cos(z),0)-V(\cos(z),0)\right]=\frac{\delta\cos(\delta\theta_\text{I})}{\sin\left(\delta z\right)-\sin(\delta\theta_\text{I})},
\end{align}
\mylinenum{which is the exact \RED{spectral function} \eqref{SMTSommFuncDirichlet} obtained using the \RED{S-M} technique.

\subsection{Neumann boundary condition}
The Neumann problem is solved in a similar way to the Dirichlet problem. The transformed Neumann \RED{BCs} are $V(\RED{\cos(\theta_{\text{w}}-z)},\pm\theta_{\text{w}})=0$, \RED{leading to a simplification of \eqref{WHM-WHsystem}. Again, in the resulting form, the W-H technique cannot be applied directly and we need the \textit{useful} mapping of Section \ref{WHM-mapping-section} together with the factorisation of $f_1(\eta(\alpha))=\sin(z(\alpha))$ and $f_3(\eta(\alpha))=\sin(z(-\alpha))$ given in the same section. This leads to}}
\begin{align}
\nonumber\frac{1}{f_{3+}(\alpha)}V(\cos(z(\alpha)),0)&=\frac{f_{3-}(\alpha)}{2}\left[U(\cos(z(-\alpha)),\theta_\text{w})-U(\cos(z(-\alpha)),-\theta_\text{w})\right],\\
\label{WHM-N(2)}\frac{f_{1+}(\alpha)}{f_{3+}(\alpha)}U(\cos(z(\alpha)),0)&=-\frac{f_{3-}(\alpha)}{2f_{1-}(\alpha)}\left[U(\cos(z(-\alpha)),\theta_\text{w})+U(\cos(z(-\alpha)),-\theta_\text{w})\right].
\end{align}
\mylinenum{\RED{The left (resp. right) sides of \eqref{WHM-N(2)} are meromorphic in $\mathcal{R}_+$ (resp. $\mathcal{R}_-$) but as for the Dirichlet case, there are potential poles on $\mathcal{R}_+\cap\mathcal{R}_-$.} To counteract this we remove the poles on the right side using the knowledge of the GO component of the solution. Assuming that $\theta_\text{w}>\pi/2$, the GO components of $U(\cos(z(-\alpha)),\pm\theta_\text{w})$ are,}
\begin{align}
\nonumber U^{(\text{GO})}(\cos(z(-\alpha)),\theta_{\text{w}})&=\frac{2i}{k(\cos(z(-\alpha))-\cos(\theta_{\text{w}}-\theta_\text{I}))},\\
\nonumber U^{(\text{GO})}(\cos(z(-\alpha)),-\theta_{\text{w}})&=\frac{2i\mathcal{H}(\pi-\theta_{\text{w}}-\theta_\text{I})}{k(\cos(z(-\alpha))-\cos(\theta_{\text{w}}+\theta_\text{I}))}.
\end{align}
\mylinenum{\RED{As in the Dirichlet case, we only need to remove the pole of $U^{(\text{GO})}(\cos(z(-\alpha)),\theta_{\text{w}})$ at $\alpha=\alpha_0=\cos\left(2\delta\theta_\text{I}\right)$ with residue,}}
\begin{align}\label{WHM-pole-resN}
\lim_{\alpha\rightarrow\alpha_0}\left[\left(\alpha-\alpha_0\right)U^{(\text{GO})}(\cos(z(-\alpha)),\theta_{\text{w}})\right]=\frac{4\delta\sin(2\delta\theta_\text{I})}{ik\sin(\theta_\text{w}-\theta_\text{I})}.
\end{align}
\mylinenum{Using this residue, and the fact that $f_{1-}(\alpha_0)=\sin(\delta\theta_\text{I})$ and $f_{3-}(\alpha_0)=\frac{\sin(\theta_\text{w}-\theta_\text{I})}{\cos(\delta\theta_\text{I})}$, we can remove this pole from the W-H system \eqref{WHM-N(2)} to get}
\begin{align}
\nonumber\frac{1}{f_{3+}(\alpha)}&V(\cos(z(\alpha)),0)-\frac{4\delta\sin(\delta\theta_\text{I})}{ik(\alpha-\cos(2\delta\theta_\text{I}))}\\
\label{WHM-N(4)}&=\frac{f_{3-}(\alpha)}{2}\left[U(\cos(z(-\alpha)),\theta_\text{w})-U(\cos(z(-\alpha)),-\theta_\text{w})\right]-\frac{4\delta\sin(\delta\theta_\text{I})}{ik(\alpha-\cos(2\delta\theta_\text{I}))},\\
\nonumber\frac{f_{1+}(\alpha)}{f_{3+}(\alpha)}&U(\cos(z(\alpha)),0)+\frac{4\delta}{ik(\alpha-\cos(2\delta\theta_\text{I})}\\
\label{WHM-N(5)}&=-\frac{f_{3-}(\alpha)}{2f_{1-}(\alpha)}\left[U(\cos(z(-\alpha)),\theta_\text{w})+U(\cos(z(-\alpha)),-\theta_\text{w})\right]+\frac{4\delta}{ik(\alpha-\cos(2\delta\theta_\text{I})}.
\end{align}

\mylinenum{In both equations \eqref{WHM-N(4)} and \eqref{WHM-N(5)}, \RED{the left (resp. right) sides are now analytic in $\mathcal{R}_+$ (resp. $\mathcal{R}_-$).} As before, in order to apply Liouville's theorem, we must determine the behaviour of each part in equations \eqref{WHM-N(4)} and \eqref{WHM-N(5)} as $|\alpha|\rightarrow\infty$. Using the edge conditions \eqref{Intro-Edge} for the Neumann case, and the reasoning developed in the Dirichlet case, we can show that, as $|\alpha|\rightarrow\infty$,
\RED{$U(\cos(z(\alpha)),0)=O(|\alpha|^{-\frac{\theta_\text{w}}{\pi}})$ and $V(\cos(z(\alpha)),0)=O(|\alpha|^{-\frac{1}{2}})$ within $\mathcal{R_+}$, while $U(\cos(z(-\alpha)),\pm\theta_\text{w})=O(|\alpha|^{-\frac{\theta_\text{w}}{\pi}})$ within $\mathcal{R_-}$.} We can also determine that $f_{1-}(\alpha),f_{3+}(\alpha)=O(|\alpha|^{\frac{1}{2}})$ and $f_{1+}(\alpha),f_{3-}(\alpha)=O(|\alpha|^{\frac{\theta_\text{w}}{\pi}-\frac{1}{2}})$ as $|\alpha|\rightarrow\infty$. \RED{As before we can hence construct two decaying entire functions and apply Liouville's theorem to obtain}}
\begin{align}\label{WHM-N-Vsol}
V(\cos(z),0)&=\frac{4\delta\sin(\delta\theta_\text{I})\cos(\delta z)}{ik(\cos(2\delta z)-\cos(2\delta\theta_\text{I}))},\\
\label{WHM-N-Usol}\sin(z)U(\cos(z),0)&=-\frac{2\delta\sin(2\delta z)}{ik(\cos(2\delta z)-\cos(2\delta\theta_\text{I})}.
\end{align}
\mylinenum{Equations \eqref{WHM-N-Vsol} and \eqref{WHM-N-Usol} can be substituted into \eqref{WHM-Somm-Sol} to get,}
\begin{align}\label{WHMSommFuncNeumann}
s(z)=\frac{ik}{2}\left[\sin(z)U(\cos(z),0)-V(\cos(z),0)\right]=\frac{\delta\cos(\delta z)}{\sin\left(\delta z\right)-\sin(\delta\theta_\text{I})},
\end{align}
\mylinenum{which is the exact \RED{spectral function} \eqref{SMTSommFuncNeumann} obtained using the \RED{S-M} technique.

\RED{ 
\paragraph{Critical analysis} 
The main disadvantage of this method is that, as our derivation shows, it is not naturally designed to tackle the wedge problem. As a result, we do not directly obtain a usual W-H equation, and have to rely on a sophisticated mapping in order to get back to the usual framework.

The advantage of this section, however, is to show the flexibility of the W-H method, and that it can work, even in non-flat/parallel geometries. As for the S-M\\ technique it is possible to adapt such method to more complicated cases such as inhomogeneous impedance \citep{AVShanin1998}, skew incidence \citep{DanieleLombardi2006} or even the penetrable wedge \citep{DanieleLombardi2011}.

Moreover, in usual flat geometries, it is known that the W-H technique can be adapted to handle finite structures. It generally results in matrix W-H problems. This is encouraging in our case since there is a chance of tackling geometries such as the truncated wedge (tip removed) with such method. It is perhaps surprising that this problem can be recast in an analytical continuation problem of the W-H problem type. This may give insight to the solution of a broader class of diffraction problems \citep{DanieleZich2014}.
}


\section{The Kontorovich-Lebedev transform method}\label{MKLT} 
The third method to be reviewed relies on the Kontorovich-Lebedev (K-L) transform. Introduced in \citep{KL1939}, this transform is useful because the resulting transformed Helmholtz equation is easy to solve. However, the inverse transform is known to have convergence issues, but there are alternative versions involving a convergence factor that can help with this (see e.g. \citep{DSJones1980}). For any function $f(r)$, define the K-L transform and its inverse (which can have many variations from numerous sources such as \citep{Lebedev1965}, \citep{Abrahams1986} and \citep{DSJones1986}) as}
\begin{align}\label{MKLT-KL-transform}
F(\nu)=\int_0^\infty \frac{f(r)}{r}H^{(1)}_\nu(kr)\text{d}r,\quad f(r)=\frac{1}{2}\int_{-i\infty}^{i\infty}\nu J_\nu(kr)F(\nu)\text{d}\nu,
\end{align}
\mylinenum{where $J_\nu$ and $H^{(1)}_\nu$ are the Bessel and Hankel functions of the first kind. The transform is valid if,}
\begin{align}\label{MKLT-KLconditions}
\left|\int_c^\infty \frac{f(r)e^{-ikr}}{r^{\frac{3}{2}}}\text{d}r\right|<\infty,\ \ \text{and}\ \ \int_0^\epsilon \left|\frac{f(r)\ln(kr)}{r}\right|\text{d}r<\infty,
\end{align}
\mylinenum{for all $c>0$ and $0<\epsilon\ll1$. Alternatively, if the second integral condition is not satisfied because $f(r)$ tends to a constant as $r\rightarrow0$, then $F(\nu)$ contains a pole at $\nu=0$ on the integration contour of the inverse transform. This pole is interpreted as,}
\begin{align}\label{MKLT-origin-pole}
\frac{1}{\nu}=\frac{1}{2}\lim_{\epsilon\rightarrow0}\left[\frac{1}{\nu-\epsilon}+\frac{1}{\nu+\epsilon}\right].
\end{align}
\mylinenum{If the integrand of the inverse transform fails to converge as $\nu\rightarrow\pm i\infty$, an alternative version with a convergence factor (proposed by \citet{DSJones1980}), should be used:}
\begin{align}\label{MKLT-KL-con-fac}
f(r)=\frac{1}{2}\lim_{\epsilon\rightarrow0}\left[\int_{-i\infty}^{i\infty}e^{\epsilon\nu^2}\nu J_\nu(kr)F(\nu)\text{d}\nu\right].
\end{align}

\mylinenum{To adapt this to our problem, we first split the total wave field $\Phi$ into its incident and scattered parts $\Phi(r,\theta)=e^{-ikr\cos(\theta-\theta_\text{I})}+\Phi_\text{S}(r,\theta)$, where the scattered part $\Phi_\text{S}$ satisfies Helmholtz's equation \eqref{Intro-Helmholtz} and two types of \RED{BCs},}
\begin{align}\label{MKLT-DBC}
\textrm{Dirichlet \RED{BCs:}}\quad&\Phi_\text{S}(r,\pm\theta_\text{w})=-e^{-ikr\cos(\theta_\text{w}\mp\theta_\text{I})},\\
\label{MKLT-NBC}\textrm{Neumann \RED{BCs:}}\quad&\frac{1}{r}\ParDer{\Phi_\text{S}}{\theta}(r,\pm\theta_\text{w})=\mp ik\sin(\theta_\text{w}\mp\theta_\text{I})e^{-ikr\cos(\theta_\text{w}\mp\theta_\text{I})}.
\end{align}
\mylinenum{\RED{For our problem, the K-L transform and the associated inverse are given below,}}
\begin{align}
\label{MKLT-myKL-transform}
\Psi(\nu,\theta)=\int_0^\infty \frac{\Phi_\text{S}(r,\theta)}{r}H^{(1)}_\nu(kr)\text{d}r,\ \ \Phi_\text{S}(r,\theta)=\frac{1}{2}\int_{-i\infty}^{i\infty}\nu J_\nu(kr)\Psi(\nu,\theta)\text{d}\nu,
\end{align}
\mylinenum{where the first integral condition \eqref{MKLT-KLconditions} is satisfied due to the radiation condition \eqref{Intro-2DSRcondition}. The edge condition \eqref{Intro-Edge} implies that $\Psi(\nu,\theta)$ will have a pole at $\nu=0$. Using \eqref{MKLT-myKL-transform}, we find the transformed boundary data,}
\begin{align}\label{MKLT-data-Dpre}
\textrm{Dirichlet:}\quad\Psi^\pm(\nu)&=\Psi(\nu,\pm\theta_\text{w})=-\int_0^\infty \frac{1}{r}e^{-ikr\cos(\theta_\text{w}\mp\theta_\text{I})}H^{(1)}_\nu(kr)\text{d}r,\\
\label{MKLT-data-Npre}\textrm{Neumann:}\quad\Psi_\theta^\pm(\nu)&=\ParDer{\Psi}{\theta}(\nu,\pm\theta_\text{w})=\mp ik\sin(\theta_\text{w}\mp\theta_\text{I})\int_0^\infty e^{-ikr\cos(\theta_\text{w}\mp\theta_\text{I})}H^{(1)}_\nu(kr)\text{d}r.
\end{align}

\mylinenum{From equation 6.611.5 of \citep{TablesISP8th}, for $\RE{\nu}\in(-1,1)$, we know that}
\begin{align}\label{MKLT-ISP-66115}
\int_0^\infty e^{-aR}H^{(1)}_\nu(R)\text{d}R=\frac{i\left(\sqrt{a^2+1}+a\right)^{-\nu}}{\sin(\pi\nu)\sqrt{a^2+1}}\left[e^{-i\pi\nu}-\left(\sqrt{a^2+1}+a\right)^{2\nu}\right],
\end{align}
\mylinenum{and integrating \eqref{MKLT-ISP-66115} with respect to $a$, we obtain}
\begin{align}\label{MKLT-ISP-66115-int}
-\int_0^\infty\frac{e^{-aR}}{R}H^{(1)}_\nu(R)\text{d}R=-\frac{i\left(\sqrt{a^2+1}+a\right)^{-\nu}}{\nu\sin(\pi\nu)}\left[e^{-i\pi\nu}+\left(\sqrt{a^2+1}+a\right)^{2\nu}\right].
\end{align}
\mylinenum{Now, let $R=kr$ and $a=i\cos(\theta_\text{w}\mp\theta_\text{I})$, \RED{then use \eqref{MKLT-ISP-66115-int} (resp. \eqref{MKLT-ISP-66115}) to evaluate \eqref{MKLT-data-Dpre} (resp. \eqref{MKLT-data-Npre}) explicitly to get}}
\begin{align}
\label{MKLT-data-D}\textrm{Dirichlet:}\quad\Psi^\pm(\nu)&=\frac{2(-i)^{1+\nu}}{\nu\sin(\pi\nu)}\cos((\theta_\text{w}\mp\theta_\text{I}-\pi)\nu),\\
\label{MKLT-data-N}\textrm{Neumann:}\quad\Psi_\theta^\pm(\nu)&=\mp\frac{2(-i)^{1+\nu}}{\sin(\pi\nu)}\sin((\theta_\text{w}\mp\theta_\text{I}-\pi)\nu).
\end{align}

\mylinenum{The advantage of the K-L transform is that if $\Psi$ satisfies the following governing equation,} 
\begin{align}\label{MKLT-Psi-gov-eq}
\frac{\partial^2\Psi}{\partial\theta^2}+\nu^2\Psi=0,
\end{align}
\mylinenum{then $\Phi_\text{S}$ satisfies Helmholtz's equation. For the Dirichlet case, the solution of \eqref{MKLT-Psi-gov-eq} is,}
\begin{align}\label{MKLT-AfterBCs-D}
\Psi(\nu,\theta)=\frac{1}{\sin(2\theta_\text{w}\nu)}\left[\Psi^-(\nu)\sin((\theta_\text{w}-\theta)\nu)+\Psi^+(\nu)\sin((\theta_\text{w}+\theta)\nu)\right].
\end{align}
\mylinenum{This means that the exact solution with Dirichlet \RED{BCs} is}
\begin{align}
\Phi(r,\theta)=&\ \int_{-i\infty}^{i\infty}\frac{J_\nu(kr)}{i^{1+\nu}\sin(\pi\nu)\sin(2\theta_\text{w}\nu)} \big[\cos((\theta_\text{w}+\theta_\text{I}-\pi)\nu)\sin((\theta_\text{w}-\theta)\nu) \nonumber\\
& +\cos((\theta_\text{w}-\theta_\text{I}-\pi)\nu)\sin((\theta_\text{w}+\theta)\nu)\big]\text{d}\nu +e^{-ikr\cos(\theta-\theta_\text{I})}. \label{MKLT-KLSol-D} 
\end{align}
\mylinenum{For the Neumann case, the solution of \eqref{MKLT-Psi-gov-eq} is,}
\begin{align}\label{MKLT-AfterBCs-N}
\Psi(\nu,\theta)=\frac{1}{\nu\sin(2\theta_\text{w}\nu)}\left[\Psi_\theta^-(\nu)\cos((\theta_\text{w}-\theta)\nu)-\Psi_\theta^+(\nu)\cos((\theta_\text{w}+\theta)\nu)\right].
\end{align}
\mylinenum{This means that the exact solution with Neumann \RED{BCs} is}
\begin{align}
\Phi(r,\theta)=&\ \int_{-i\infty}^{i\infty}\frac{J_\nu(kr)}{i^{1+\nu}\sin(\pi\nu)\sin(2\theta_\text{w}\nu)} \big[\sin((\theta_\text{w}+\theta_\text{I}-\pi)\nu)\cos((\theta_\text{w}-\theta)\nu) \nonumber\\ 
&+\sin((\theta_\text{w}-\theta_\text{I}-\pi)\nu)\cos((\theta_\text{w}+\theta)\nu)\big]\text{d}\nu + e^{-ikr\cos(\theta-\theta_\text{I})}. \label{MKLT-KLSol-N}
\end{align}

\begin{figure}[ht]\centering
\includegraphics[width=0.545\textwidth]{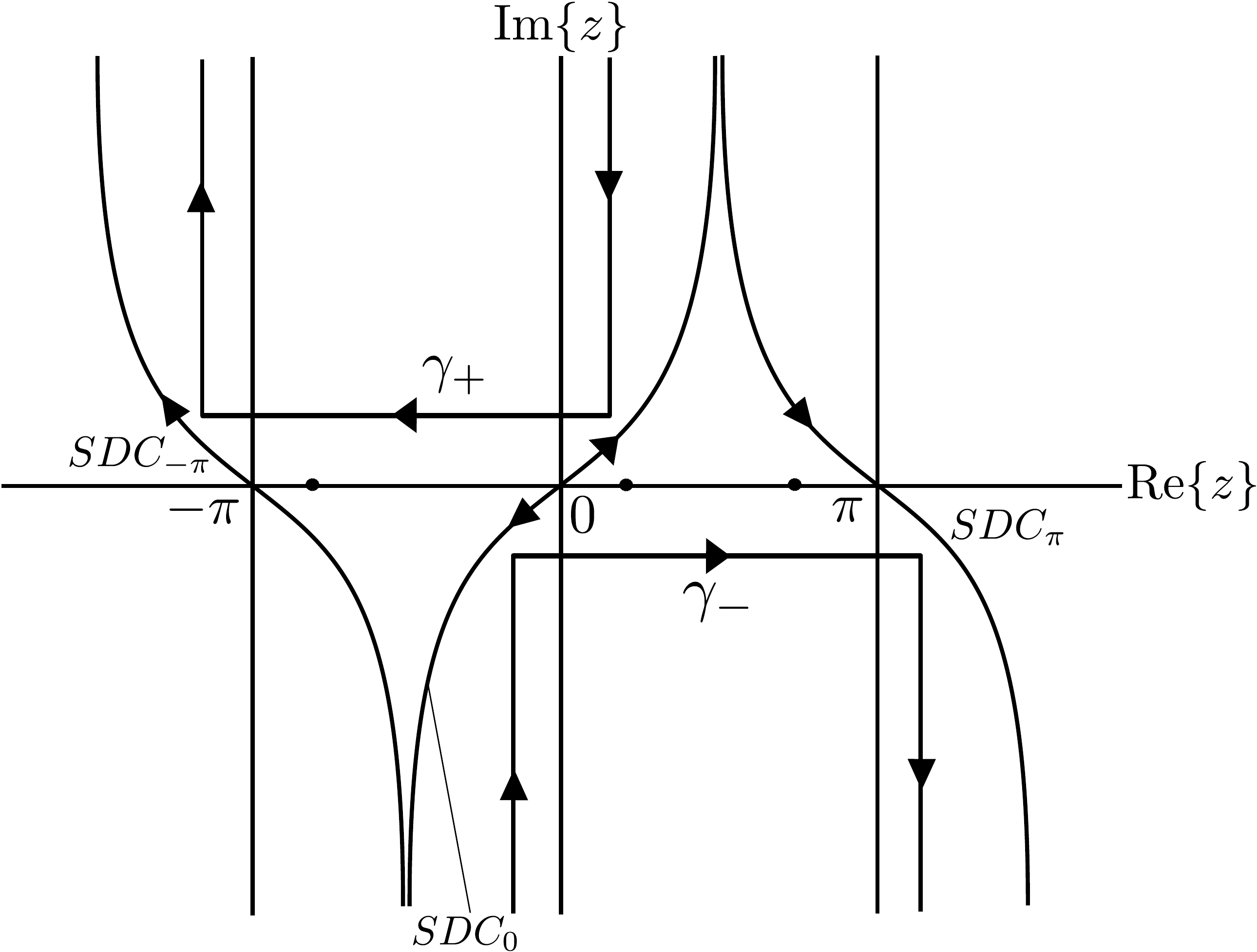}
\includegraphics[width=0.445\textwidth]{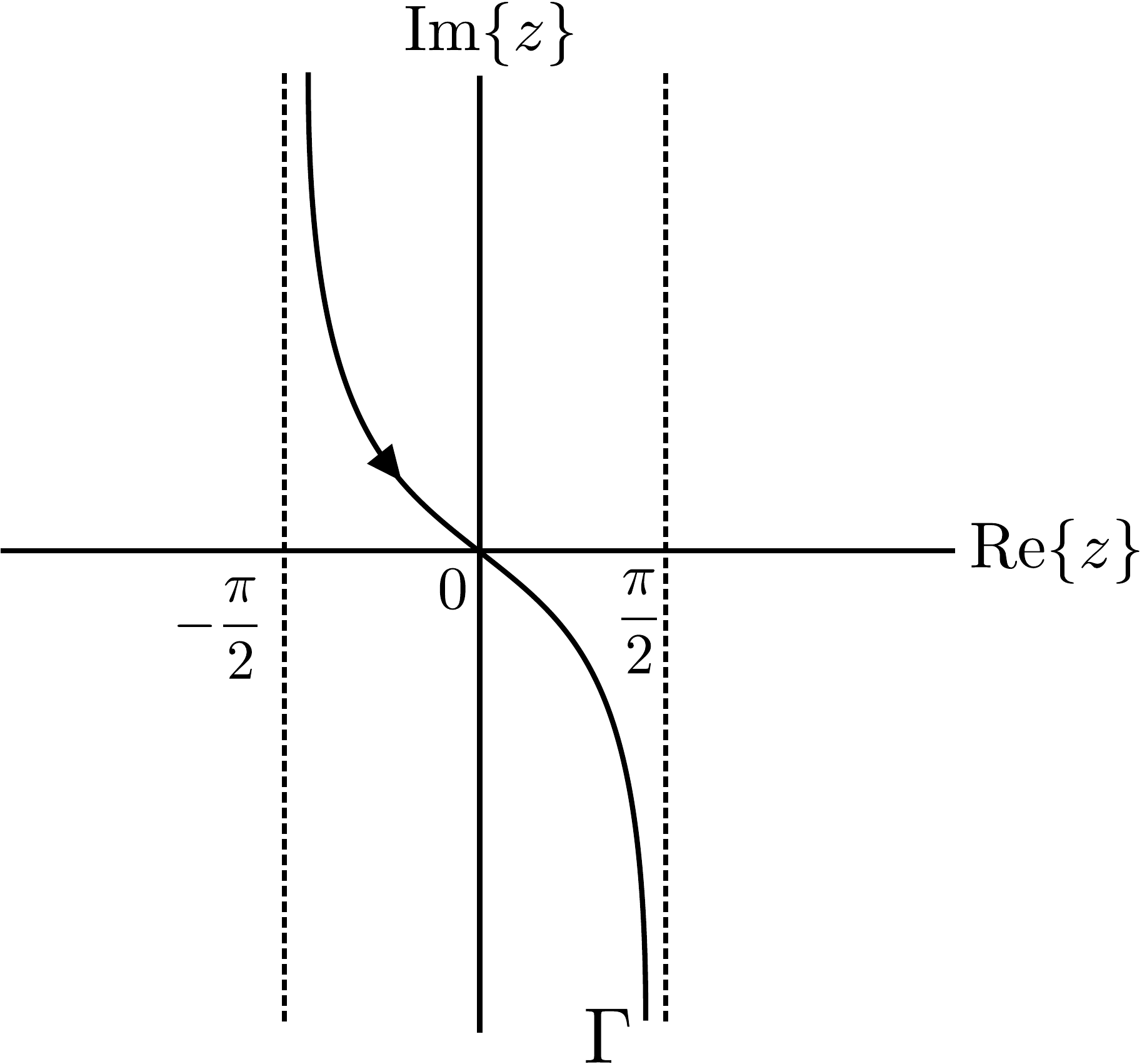}
\caption{The Sommerfeld contours $\gamma_\pm$, the steepest descent contours $SDC_{-\pi}$, $SDC_0$, $SDC_\pi$ and, shown as dots, possible poles on the real line (left) and the $\Gamma$ contour (right) }
\label{SA-SDCs}
\end{figure}

\mylinenum{While it is difficult to see it by inspection, these integral solutions \eqref{MKLT-KLSol-D} and \eqref{MKLT-KLSol-N} are equivalent to the Sommerfeld integral solutions \eqref{SMTSommIntDirichlet} and \eqref{SMTSommIntNeumann}. \RED{As we discuss in \ref{app:app1}, the} connection between the Sommerfeld inverse formula \eqref{SMTSommInt} and the K-L inverse transform \eqref{MKLT-KL-transform} was made in \citep{Malyuzhinets1958-2} by using the Sommerfeld integral form of Bessel functions. Here, we shall show equivalence by first rewriting the integrals \eqref{MKLT-KLSol-D} and \eqref{MKLT-KLSol-N} as series, then convert the result into Sommerfeld integrals.

The solutions \eqref{MKLT-KLSol-D} and \eqref{MKLT-KLSol-N} are evaluated by deforming the contour to the right and summing the residues of the poles crossed in the process. The double pole at $\nu=0$ is interpreted using \eqref{MKLT-origin-pole}. We can simplify the result using the Jacobi-Anger expansion of the incident wave to obtain the following series solutions for Dirichlet and Neumann \RED{BCs} respectively,}
\begin{align}
\label{MKLT-SeriesSol-D}\Phi(r,\theta)&=2\delta\sum_{n=1}^\infty(-i)^{\delta n}J_{\delta n}(kr)\left[\cos((\theta-\theta_\text{I})\delta n)-\cos((\theta-2\theta_\text{w}+\theta_\text{I})\delta n)\right],\\
\label{MKLT-SeriesSol-N}\Phi(r,\theta)&=2\delta J_{0}(kr)+2\delta\sum_{n=1}^\infty(-i)^{\delta n}J_{\delta n}(kr)\left[\cos((\theta-\theta_\text{I})\delta n)+\cos((\theta-2\theta_\text{w}+\theta_\text{I})\delta n)\right],
\end{align}
\mylinenum{where, as before, $\delta=\frac{\pi}{2\theta_\text{w}}$. These series solutions can be matched with classical series solutions obtained by \citet{Macdonald1902} \RED{(see \ref{Mac} for more details)}.

Finally, we need to transform \eqref{MKLT-SeriesSol-D} and \eqref{MKLT-SeriesSol-N} into Sommerfeld integrals. We do this by using the Sommerfeld integral formula for the Bessel function of the first kind,}
\begin{align}\label{MKLT-BesselSI}
J_v(R)=-\frac{1}{2\pi}\int_{\gamma_+}e^{-iR\cos(z)+ivz+iv\frac{\pi}{2}}\text{d}z,
\end{align}
\mylinenum{and equation 1.461.2 from \citep{TablesISP8th}}
\begin{align}\label{MKLT-ISP-14162}
1+2\sum_{n=1}^\infty e^{ina}\cos(nb)=\frac{i\sin(a)}{\cos(b)-\cos(a)},
\end{align}
\mylinenum{which converges if $\IM{a}>0$. This means that \eqref{MKLT-SeriesSol-D} and \eqref{MKLT-SeriesSol-N} can be written in Sommerfeld integral form as}
\begin{align}\label{MKLT-SommInt}
\Phi(r,\theta)&=\frac{1}{2\pi i}\int_{\gamma_+}e^{-ikr\cos(z)}\left[\frac{\delta\sin(\delta z)}{\cos(\delta(\theta-\theta_\text{I}))-\cos(\delta z)}\pm\frac{\delta\sin(\delta z)}{\cos(\delta(\theta+\theta_\text{I}))+\cos(\delta z)}\right]\text{d}z,
\end{align}
\mylinenum{where the plus and minus signs denote the Dirichlet and Neumann solutions respectively. Using standard trigonometric identities, it is trivial to show that the square brackets in \eqref{MKLT-SommInt} are an alternate form of $s(\theta+z)-s(\theta-z)$. Hence the Kontorovich-Lebedev solutions \eqref{MKLT-KLSol-D} and \eqref{MKLT-KLSol-N} match with the Sommerfeld integrals \eqref{SMTSommIntDirichlet} and \eqref{SMTSommIntNeumann} respectively.

\RED{
\paragraph{Critical analysis} The main advantage of such transform is that it is a very natural way to tackle the wedge problem. It hence leads to a constructive proof of the form of the solution in the K-L space, see e.g. \eqref{MKLT-AfterBCs-D}. In addition, it leads easily to a near-field expansion of the solution.

The clear disadvantage of such method lies in the convergence issues of the inverse K-L transform. It does require some regularisation in order to be evaluated numerically, and even in this case, the computation of inverse K-L transform remains cumbersome. To evaluate the far-field it is usually necessary to convert it into a Sommerfeld type integral.
}


\section{Solution analysis and evaluation}\label{SolAnalysis} 
We shall now compare the exact integral and series solutions with some GTD approximations and evaluate them for some representative values of $\theta_{\text{w}}$. Note that the K-L integrals do not need to be plotted since we have already shown that they are equivalent to the Sommerfeld integrals \RED{\eqref{SMTSommIntDirichlet} and \eqref{SMTSommIntNeumann}. Numerical computation of the Sommerfeld} integrals can be slow if $kr\gg1$ because $e^{-ikr\cos(z)}$ will oscillate rapidly along the Sommerfeld contours. 

Another way to evaluate these integrals is to deform the Sommerfeld contours to the steepest descent contours shown on the left side of Figure \ref{SA-SDCs}. During this deformation, all poles on the real line segment $|\RE{z}|\leq\pi$ are crossed. Their contribution, which can be calculated exactly using residues, correspond to the GO component of the field, $\Phi_{\text{GO}}$, leaving behind the diffracted part $\Phi_{\text{Diff}}$.
The steepest descent contour $SDC_0$ is repeated twice in opposite directions so is cancelled out and the other two contours are translations of each other. The exponential term $e^{-ikr\cos(z)}$ does not oscillate along these contours so computation time is significantly \RED{reduced}, even for large $kr$. This means that the \RED{S-M integrals \eqref{SMTSommIntDirichlet}-\eqref{SMTSommIntNeumann}} are equivalent to,}
\begin{align}
\label{SA-SDCInt}\Phi&=\Phi_{\text{GO}}+\Phi_{\text{Diff}}=\Phi_{\text{GO}}+\frac{1}{2\pi i}\int_{SDC_{-\pi}+SDC_\pi}e^{-ikr\cos(z)}s(\theta+z)\text{d}z.
\end{align}

\mylinenum{Since they are translations of each other, we can transform $SDC_{-\pi}$ and $SDC_{\pi}$ onto the $\Gamma$ contour which is illustrated on the right side of Figure \ref{SA-SDCs},}
\begin{align}
\label{SA-GammaInt}\Phi_{\text{Diff}}&=\frac{1}{2\pi i}\int_{\Gamma}e^{ikr\cos(z)}\left[s(\theta+z+\pi)-s(\theta+z-\pi)\right]\text{d}z.
\end{align}
\mylinenum{By the \RED{method of steepest descent}, $\Gamma$ satisfies,}
\begin{align}
\label{SA-Gammaequ}\Gamma(\tau)=\cosh^{-1}\left(\frac{1}{\cos(\tau)}\right),\quad\text{and}\quad\sin(\tau)\sinh(\Gamma(\tau))\leq0,
\end{align}
\mylinenum{where $\RE{z}=\tau\in(-\frac{\pi}{2},\frac{\pi}{2})$. \RED{This is rewritten in terms of the Gudermannian function $\textrm{gd}(x)$,}}
\RED{\begin{align}
\label{SA-Gammaequ2}\Gamma(\tau)=i\textrm{gd}(i\tau)=\ln\left|\sec(\tau)-\tan(\tau)\right|.
\end{align}}
\mylinenum{Using the following parametrisation, $z=\tau+i\Gamma(\tau)$, and noting these identities,}
\begin{align}
\nonumber\sinh(\Gamma(\tau))=-\tan(\tau),\ \ \Der{\Gamma}{\tau}=-\sec(\tau),
\end{align}
\mylinenum{\RED{the diffracted part is written as a simple integral:}}
\begin{align}
\label{SA-GammaInt2}\Phi_{\text{Diff}}&=\frac{e^{ikr}}{2\pi i}\int_{-\frac{\pi}{2}}^{\frac{\pi}{2}}e^{-kr\sin(\tau)\tan(\tau)}\big[s(\theta+\pi+\tau+i\Gamma(\tau))-s(\theta-\pi+\tau+i\Gamma(\tau))\big]\left(1-i\sec(\tau)\right)\text{d}\tau.
\end{align}

\mylinenum{As stated earlier, this integral will be much faster to evaluate numerically than the \RED{S-M integrals \eqref{SMTSommIntDirichlet}-\eqref{SMTSommIntNeumann}.} However, difficulties can arise when $\theta$ is in a small neighbourhood of the GO discontinuities because \RED{one} of the poles will be very close to the contour of integration, which will cause numerical issues.


\subsection{Comparison with simpler problems}\label{SA-comparison}
In this subsection we will show that the solution is consistent with the simple case when the wedge opens up to \RED{form a half-space or closes to make} a half-plane. First, we look at the case where $\theta_{\text{w}}=\pi/2$ to form a \RED{half-space problem. The} solution is easily obtainable via the method of images,}
\begin{align}\label{SA-InfPlaneSol}
\Phi(r,\theta)=\Phi_\text{I}\mp\Phi_\text{R}.
\end{align}
\mylinenum{\RED{Here} the upper and lower signs denote the Dirichlet and Neumann solutions respectively. Obviously \eqref{SA-InfPlaneSol} is equal to the GO component so we need to show that the diffracted part \eqref{SA-GammaInt} is identically zero. Expressing $s(z)$ in terms of the cotangent, we find that,}
\begin{align}
\nonumber s(z+\pi)-s(z-\pi)=\ &\frac{1}{2}\left[\cot\left(\frac{1}{2}(z-\theta_\text{I}+\pi)\right)\mp\cot\left(\frac{1}{2}(z+\theta_\text{I})\right)\right]\\
\label{ISWinfphid2}&-\frac{1}{2}\left[\cot\left(\frac{1}{2}(z-\theta_\text{I}-\pi)\right)\mp\cot\left(\frac{1}{2}(z+\theta_\text{I}-2\pi)\right)\right],
\end{align}
\mylinenum{which is identically zero due to the periodicity of cotangent. This implies that \RED{$\Phi_\mathrm{Diff}\equiv0$}, as required. 

For another comparison, we look at the case where the wedge closes to form a half-plane. Hence we let $\theta_{\text{w}}=\pi$ and match the \RED{S-M integrals \eqref{SMTSommIntDirichlet}-\eqref{SMTSommIntNeumann}} with the known solution to the half-plane problem in terms of Fresnel integrals,}
\begin{align}
\nonumber\Phi(r,\theta)=\ &\Phi_\text{I}\left[\frac{1}{2}+\pi^{-\frac{1}{2}}e^{-\frac{i\pi}{4}}\mathcal{F}\left((2kr)^\frac{1}{2}\cos\left(\frac{1}{2}(\theta-\theta_\text{I})\right)\right)\right]\\
\label{SA-HalfPlaneSol}&\mp\Phi_\text{R}\left[\frac{1}{2}-\pi^{-\frac{1}{2}}e^{-\frac{i\pi}{4}}\mathcal{F}\left((2kr)^\frac{1}{2}\cos\left(\frac{1}{2}(\theta+\theta_\text{I})\right)\right)\right],
\end{align}
\mylinenum{where the upper and lower signs correspond to the Dirichlet and Neumann solution respectively, $\Phi_\text{R}=e^{-ikr\cos(\theta+\theta_\text{I})}$ is the reflected wave and $\mathcal{F}(v)$ is the Fresnel integral defined\footnote{\RED{Fresnel integrals can be written in many different ways, see for example \citep{Handbook,Noble1958,Assier2012b} and references therein.}} by}
\begin{align}
\label{SA-FresnelInt}\mathcal{F}(v)=\int_0^ve^{iu^2}\text{d}u.
\end{align}
\mylinenum{\RED{Having} $\theta_{\text{w}}=\pi$ implies that $\delta=1/2$ and,}
\begin{align}
\label{SA-sfunction}s(z)=\frac{1}{4}\left[\cot\left(\frac{1}{4}(z-\theta_\text{I})\right)\mp\cot\left(\frac{1}{4}(z+\theta_\text{I}-2\pi)\right)\right].
\end{align}
\mylinenum{Hence, we can rewrite the \RED{S-M integrals \eqref{SMTSommIntDirichlet}-\eqref{SMTSommIntNeumann}} in the following form,}
\begin{align}
\label{SA-HalfPlaneSol2}\Phi(r,\theta)=\Phi_\mathcal{F}(r,\theta-\theta_\text{I})\mp\Phi_\mathcal{F}(r,\theta+\theta_\text{I}-2\pi),
\end{align}
\mylinenum{where}
\begin{align}
\label{SA-PhiFint}\Phi_\mathcal{F}(r,\lambda)=\frac{1}{8\pi i}\int_{\gamma_++\gamma_-}e^{-ikr\cos(z)}\cot\left(\frac{1}{4}(z+\lambda)\right)\text{d}z.
\end{align}
\mylinenum{It is possible to express \eqref{SA-PhiFint} in terms of a Fresnel integral (a procedure to do this can be found in section 5.3 in \citep{SMtechnique2007}), leading to}
\begin{align}
\label{SA-PhiF}\Phi_\mathcal{F}(r,\lambda)=e^{-ikr\cos(\lambda)}\left[\frac{1}{2}+\pi^{-\frac{1}{2}}e^{-\frac{i\pi}{4}}\mathcal{F}\left((2kr)^{\frac{1}{2}}\cos\left(\frac{\lambda}{2}\right)\right)\right].
\end{align}
\mylinenum{\RED{Using this, we recover exactly \eqref{SA-HalfPlaneSol} from \eqref{SA-HalfPlaneSol2}, as expected.} Now that the solution matches with that of the \RED{half-space} and half-plane problems, we shall focus on deriving the GTD approximation for non-degenerate wedges.


\subsection{Geometrical Theory of Diffraction (GTD)}\label{SA-secGTD}
\citet{Keller1962} defined the Geometrical Theory of Diffraction to be an extension of classic Geometrical Optics including diffraction terms. The GTD approximation is simply an asymptotic approximation of the total wave field as $kr\rightarrow\infty$, creating a high-frequency or far-field approximation. To derive the GTD approximation of the case presented here, we continue with the \RED{method of steepest descent} applied to \eqref{SA-GammaInt} as $kr\rightarrow\infty$. Equation \eqref{SA-GammaInt} is of the form $\Phi_{\text{Diff}}=\int_{\Gamma}g(z)e^{-kr\psi(z)}\,\mathd z$, where $kr$ is a big parameter, $g(z)=\frac{(s(\theta+z+\pi)-s(\theta+z-\pi))}{2\pi i}$ and $\psi(z)=-i\cos(z)$. The latter has a saddle point at $z = 0$ and is such that $\psi''(0)=i\neq 0$. Since $g(0)$ is also not zero, we can apply the method of steepest descent in its simplest form (see e.g. \citep{Bleistein}) to get}
\begin{align}
\Phi_{\text{Diff}}\underset{kr\rightarrow\infty}{\sim}&\sqrt{\frac{2\pi}{kr\psi''(0)}}g(0)e^{-kr\psi(0)}=\frac{e^{ikr+i\pi/4}}{\sqrt{2\pi kr}}(s(\theta-\pi)-s(\theta+\pi)) 
\end{align}
\mylinenum{Hence we can write}
\begin{align}
\label{SA-GTD} \Phi(r,\theta)\underset{kr \rightarrow \infty}{\sim}\Phi_{\text{GO}}+\frac{e^{ikr+\frac{i\pi}{4}}}{\sqrt{2\pi kr}}\left[s(\theta-\pi)-s(\theta+\pi)\right].
\end{align}

\mylinenum{In this GTD approximation, the term,}
\begin{align}
D(\theta,\theta_\text{I})=\frac{e^{\frac{i\pi}{4}}}{\sqrt{2\pi}}\left[s(\theta-\pi)-s(\theta+\pi)\right],
\label{SA-diffraction-coefficient}
\end{align}
\mylinenum{is known as the diffraction coefficient. Unfortunately, this GTD approximation is singular for certain values of $\theta$, for example in the case where $\theta_{\text{w}}>\pi/2$, the GTD is invalid at $\theta_\text{I}-\pi$, $2\theta_{\text{w}}-\theta_\text{I}-\pi$ and $\pi-2\theta_{\text{w}}-\theta_\text{I}$, which correspond to the GO discontinuities. This is the main issue with GTD: while it is a much more accurate approximation than the Geometrical Optics, it becomes invalid at the GO discontinuities. The pursuit of an approximation that is uniformly valid for all $\theta$ has led to the improved Uniform Geometrical Theory of Diffraction \citep{KP1974}. We follow section 5.5 in \citep{SMtechnique2007} to find the uniform GTD approximation (UTD). 

Restricting ourselves to the specific case where $\theta_\text{w}>\pi/2$ and $|\theta_\text{I}|<\bar{\theta}_\text{w}=\pi-\theta_\text{w}$, there are only two values where the standard GTD is invalid, $\theta=\pi-2\theta_{\text{w}}-\theta_\text{I}$ and $2\theta_{\text{w}}-\pi-\theta_\text{I}$. To produce the uniform approximation, we first construct a function that is a linear combination of $\Phi$ and $\Phi_\mathcal{F}$ defined by \eqref{SA-PhiFint}. The idea is to remove the poles causing the singularities in \eqref{SA-GTD} and then use the \RED{method of steepest descent}. Consider the following,}
\begin{align}
\label{SA-UTD1st-Xi}\RED{\Xi}(r,\theta)=&\ \Phi(r,\theta)\pm\Phi_\mathcal{F}(r,\theta+\theta_\text{I}-2\theta_{\text{w}})\pm\Phi_\mathcal{F}(r,\theta+\theta_\text{I}+2\theta_{\text{w}}),
\end{align}
\mylinenum{where $\Phi_\mathcal{F}$ is defined in \eqref{SA-PhiFint}. The upper and lower signs denote the Dirichlet and Neumann solutions respectively. The combination of $\Phi$ and $\Phi_\mathcal{F}$ has effectively removed the poles at $2\theta_{\text{w}}-\theta_\text{I}-\theta$ and $-2\theta_{\text{w}}-\theta_\text{I}-\theta$, but the pole at $\theta_\text{I}-\theta$ remains for all values of $\theta$. We use the \RED{method of steepest descent} to approximate \RED{$\Xi$},}
\begin{align}
\nonumber\RED{\Xi}(r,\theta)\sim&\ \Phi_\text{I}+\frac{e^{ikr+\frac{i\pi}{4}}}{\sqrt{2\pi kr}}\bigg[s(\theta-\pi)-s(\theta+\pi)\mp\frac{1}{2}\sec\left(\frac{1}{2}(\theta+\theta_\text{I}-2\theta_{\text{w}})\right)\\
\label{SA-UTD1st-MoSD}&\mp\frac{1}{2}\sec\left(\frac{1}{2}(\theta+\theta_\text{I}+2\theta_{\text{w}})\right)\bigg].
\end{align}
\mylinenum{We rearrange \eqref{SA-UTD1st-Xi} and \RED{use} \eqref{SA-PhiF} and \eqref{SA-UTD1st-MoSD} to find the \RED{UTD} approximation.}
\begin{align}
\nonumber\Phi(r,\theta)\sim&\ \Phi_\text{I}\mp\Phi_{\text{R}_1}\!\!\left[\frac{1}{2}+\pi^{-\frac{1}{2}}e^{-\frac{i\pi}{4}}\mathcal{F}\left((2kr)^{\frac{1}{2}}\cos\left(\frac{\theta+\theta_\text{I}-2\theta_{\text{w}}}{2}\right)\right)\right]\\
&\mp\Phi_{\text{R}_2}\!\!\left[\frac{1}{2}+\pi^{-\frac{1}{2}}e^{-\frac{i\pi}{4}}\mathcal{F}\left((2kr)^{\frac{1}{2}}\cos\left(\frac{\theta+\theta_\text{I}+2\theta_{\text{w}}}{2}\right)\right)\right]\label{SA-UTD1st-sol}\\
&+\frac{e^{ikr+\frac{i\pi}{4}}}{\sqrt{2\pi kr}}\bigg[s(\theta-\pi)-s(\theta+\pi) \mp\frac{1}{2}\sec\left(\frac{1}{2}(\theta+\theta_\text{I}-2\theta_{\text{w}})\right)\mp\frac{1}{2}\sec\left(\frac{1}{2}(\theta+\theta_\text{I}+2\theta_{\text{w}})\right)\bigg]\nonumber
\end{align}
\mylinenum{where $\Phi_{\text{R}_1}=e^{-ikr\cos(\theta-2\theta_\text{w}+\theta_\text{I})}$ and $\Phi_{\text{R}_2}=e^{-ikr\cos(\theta+2\theta_\text{w}+\theta_\text{I})}$ are the reflections of the incident wave from the top and bottom face respectively. If we restricted the incident angle to $\theta_\text{I}>\pi-\theta_{\text{w}}$ instead, we would need to use the following function,}
\begin{align}
\label{SA-UTD2nd-Xi}\RED{\Xi}(r,\theta)=\Phi(r,\theta)-\Phi_\mathcal{F}(r,\theta-\theta_\text{I})\pm\Phi_\mathcal{F}(r,\theta+\theta_\text{I}-2\theta_{\text{w}}),
\end{align}
\mylinenum{where the same method as the first case will produce \RED{another UTD} approximation for $\Phi$}
\begin{align}
\nonumber\Phi(r,\theta)\sim\ &\Phi_\text{I}\left[\frac{1}{2}+\pi^{-\frac{1}{2}}e^{-\frac{i\pi}{4}}\mathcal{F}\left((2kr)^{\frac{1}{2}}\cos\left(\frac{\theta-\theta_\text{I}}{2}\right)\right)\right]\\
&\mp\Phi_{\text{R}_1}\!\!\left[\frac{1}{2}+\pi^{-\frac{1}{2}}e^{-\frac{i\pi}{4}}\mathcal{F}\left((2kr)^{\frac{1}{2}}\cos\left(\frac{\theta+\theta_\text{I}-2\theta_{\text{w}}}{2}\right)\right)\right]\label{SA-UTD2nd-sol}\\
&+\frac{e^{ikr+\frac{i\pi}{4}}}{\sqrt{2\pi kr}}\left[s(\theta-\pi)-s(\theta+\pi)+\frac{1}{2}\sec\left(\frac{1}{2}(\theta-\theta_\text{I})\right)\mp\frac{1}{2}\sec\left(\frac{1}{2}(\theta+\theta_\text{I}-2\theta_{\text{w}})\right)\right].\nonumber
\end{align}

\mylinenum{These two approximations are uniformly valid for $-\theta_{\text{w}}<\theta<\theta_{\text{w}}$, however the \RED{BCs} are only satisfied in the limit $kr\rightarrow\infty$. Another potential accuracy issue occurs when $\theta_\text{I}$ approaches $\bar{\theta}_\text{w}$. This situation corresponds to a transition in the GO field, from a case when only one reflected wave is present to a case when two reflected waves occur. Finally, note that using the asymptotic expansions for large argument for the Fresnel integrals will simplify the above formulas to produce the GTD approximation \eqref{SA-GTD} again.


\subsection{Graphical comparison of evaluation methods}\label{plotcomparison}
The exact solution to the perfect wedge problem has been written as a Sommerfeld integral on the usual Sommerfeld contour as in \eqref{SMTSommIntDirichlet}-\eqref{SMTSommIntNeumann} or on its steepest descent contour as in \eqref{SA-GammaInt2}. Both formulations are exact and equivalent, but the latter is much easier to evaluate numerically. We have also presented three different approximations, a truncated infinite series \eqref{MKLT-SeriesSol-D}, a GTD approximation \eqref{SA-GTD} and a UTD approximation  \eqref{SA-UTD1st-sol} or \eqref{SA-UTD2nd-sol}.
In this subsection, we will plot the exact solution and each of the approximations and compare their accuracy and computational speed. For the series solutions we shall truncate at 100 terms, which is enough for the wavenumbers considered here.

In Figure \ref{fig:phicompare0}, we will consider the wedge defined by $2\bar{\theta}_{\text{w}}=\pi/4$ for zero incidence angle, $\theta_\text{I}=0$. This corresponds to a case where the GO part of the field exhibits two reflected waves. In Figure \ref{fig:phicomparehalfpi}, we consider the same wedge, but with an incident angle $\theta_\text{I}=\pi/2$, corresponding to a GO field with a single reflected wave. In both cases, we will plot the real part of the total field $\Phi$ against $\theta$ for different values of $kr$ and different \RED{BCs}. In both figures, the thick plain line represents the exact Sommerfeld solution (SI/SDC), the thick dashed line is the truncated series approximation, the dotted line and the thin line represent the UTD and GTD approximations respectively.
\begin{figure}[h!]\centering
\includegraphics[width=0.45\textwidth]{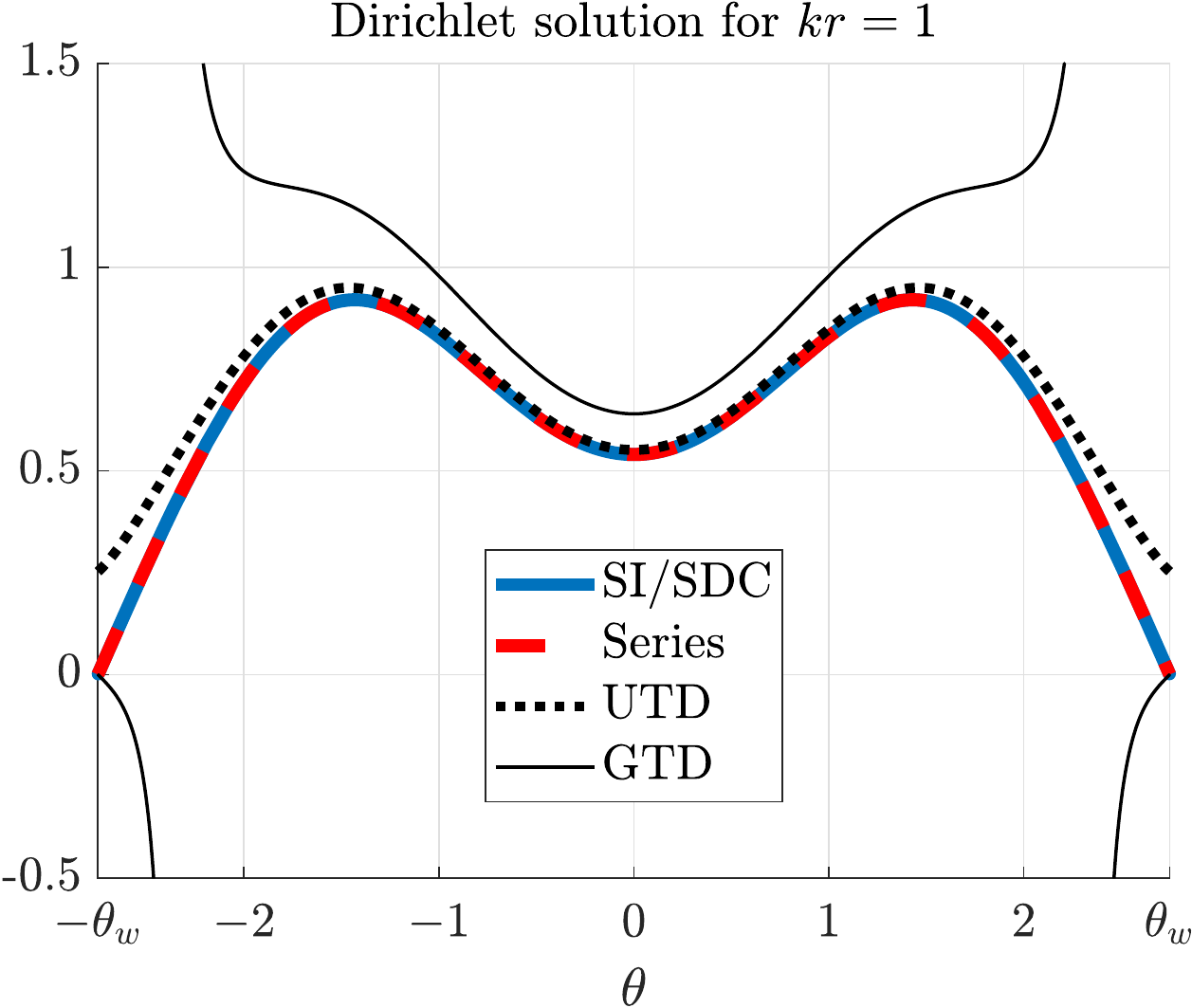}\qquad
\includegraphics[width=0.45\textwidth]{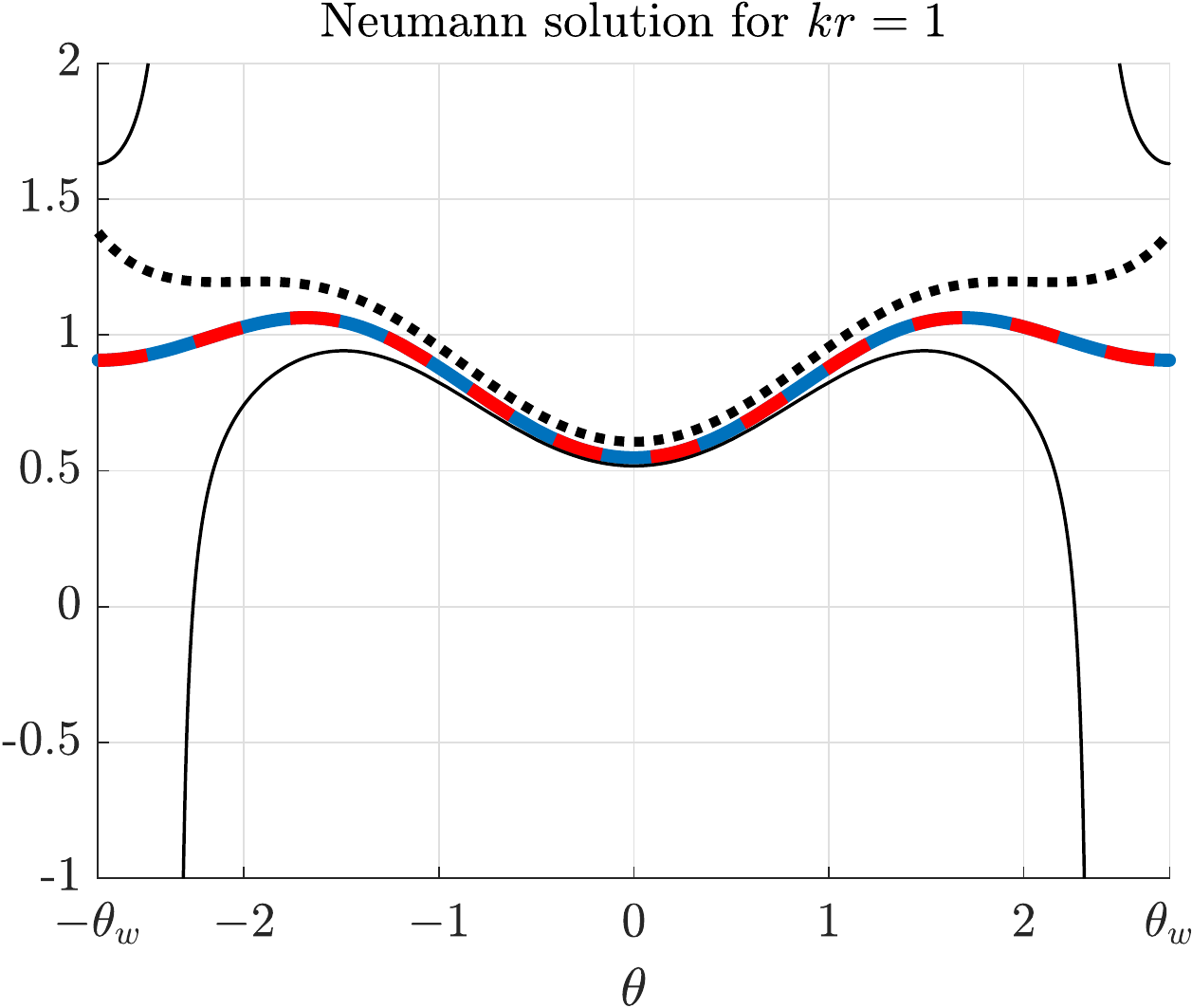}\\
\includegraphics[width=0.45\textwidth]{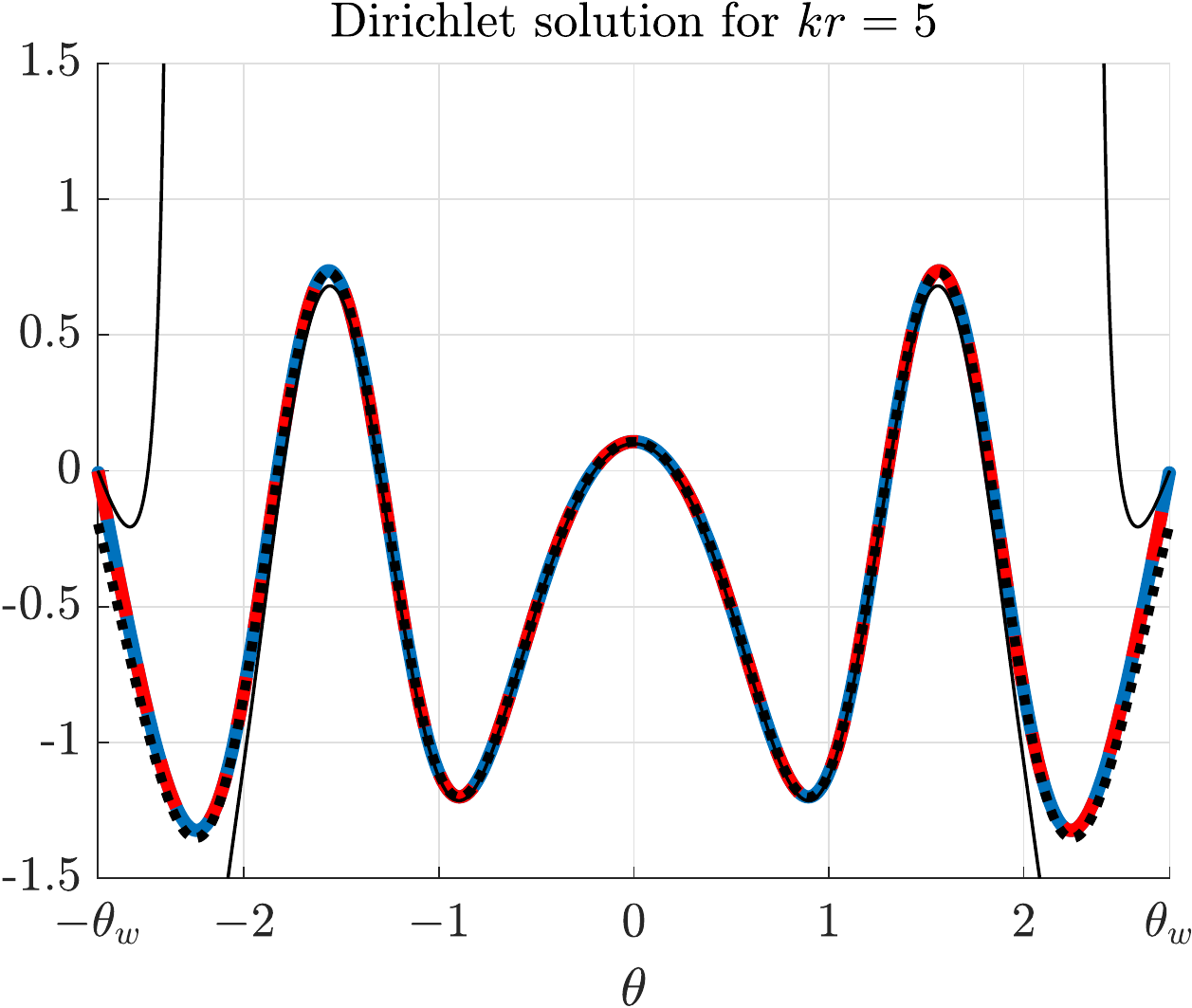}\qquad
\includegraphics[width=0.45\textwidth]{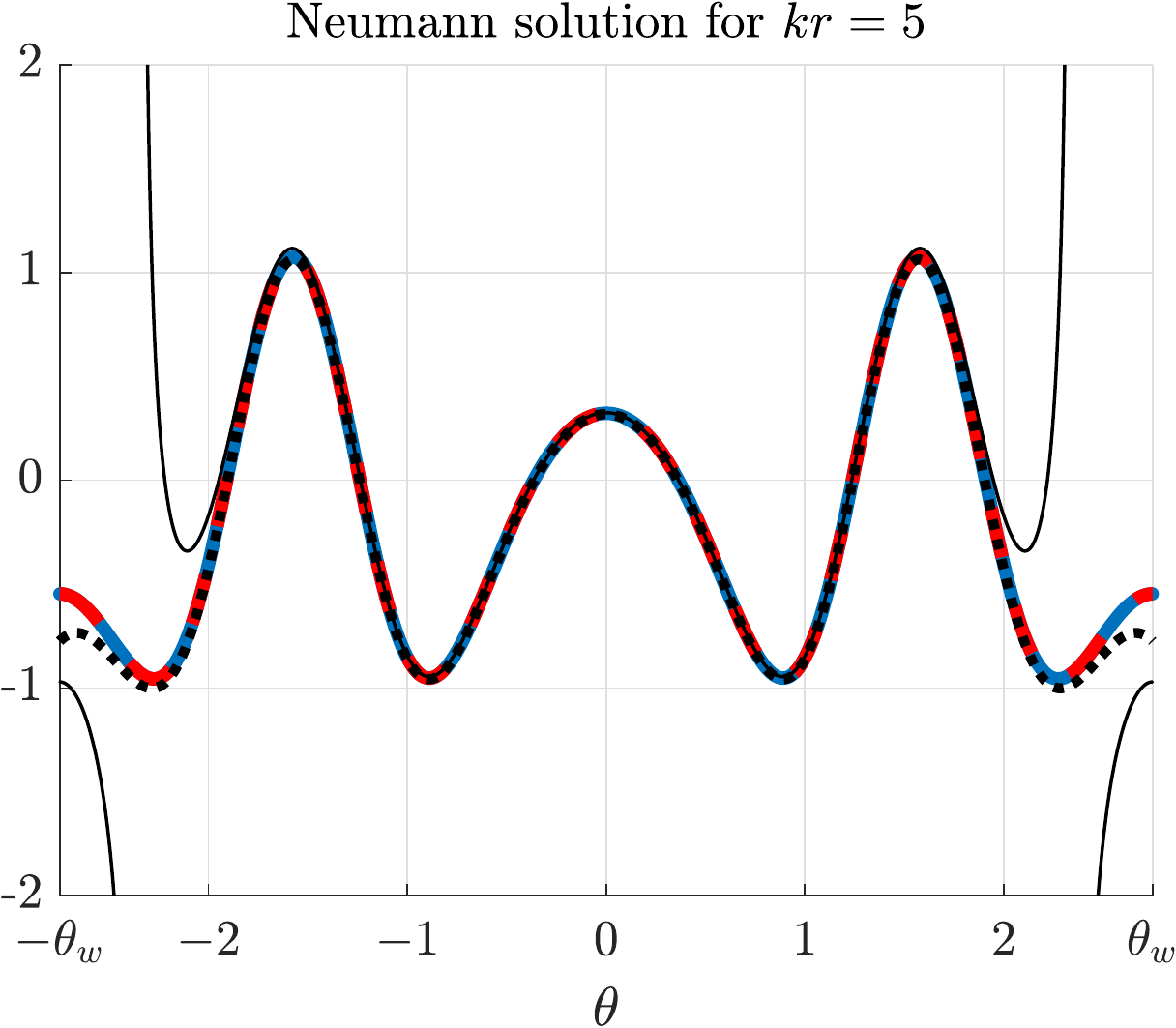}\\
\includegraphics[width=0.45\textwidth]{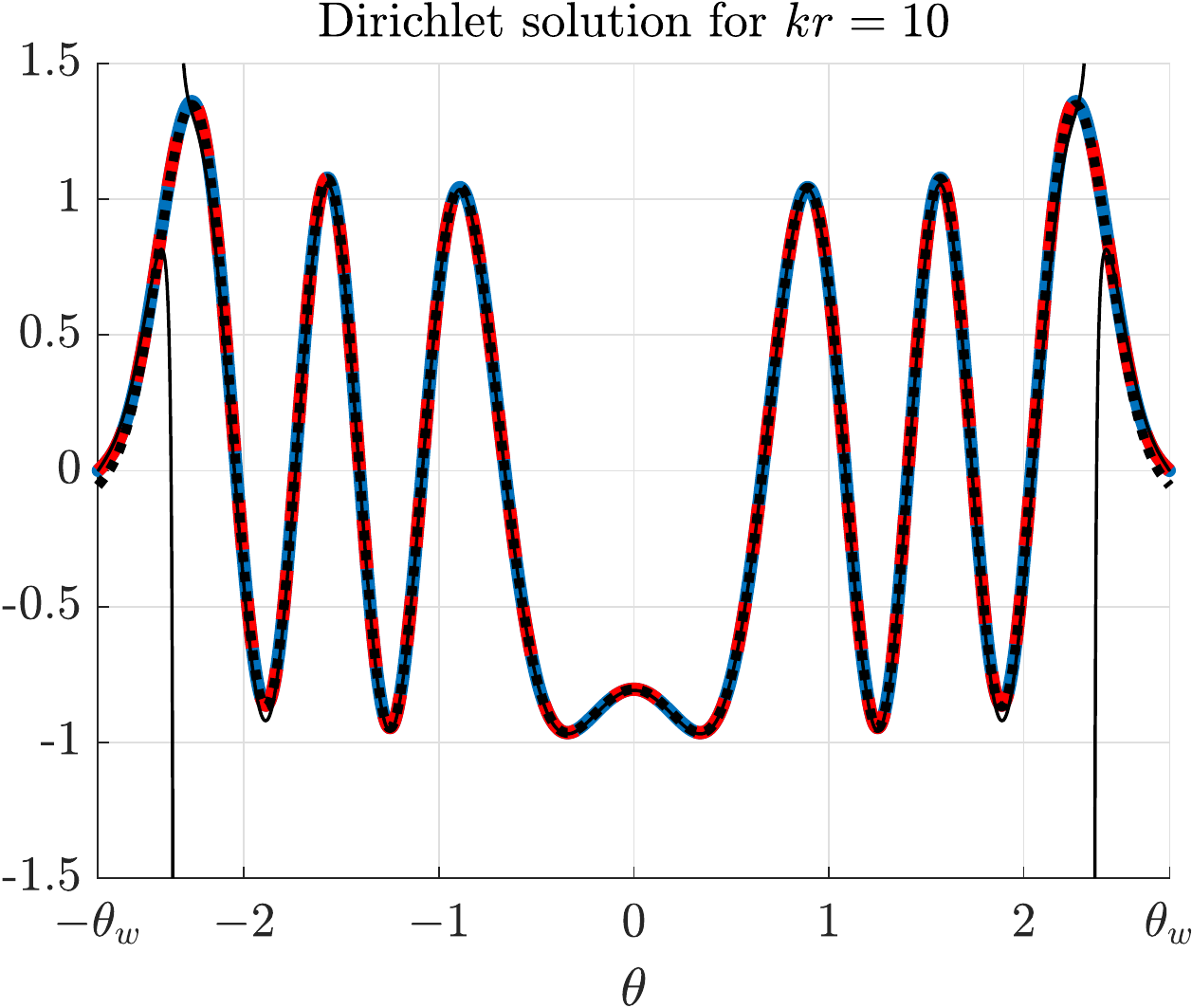}\qquad
\includegraphics[width=0.45\textwidth]{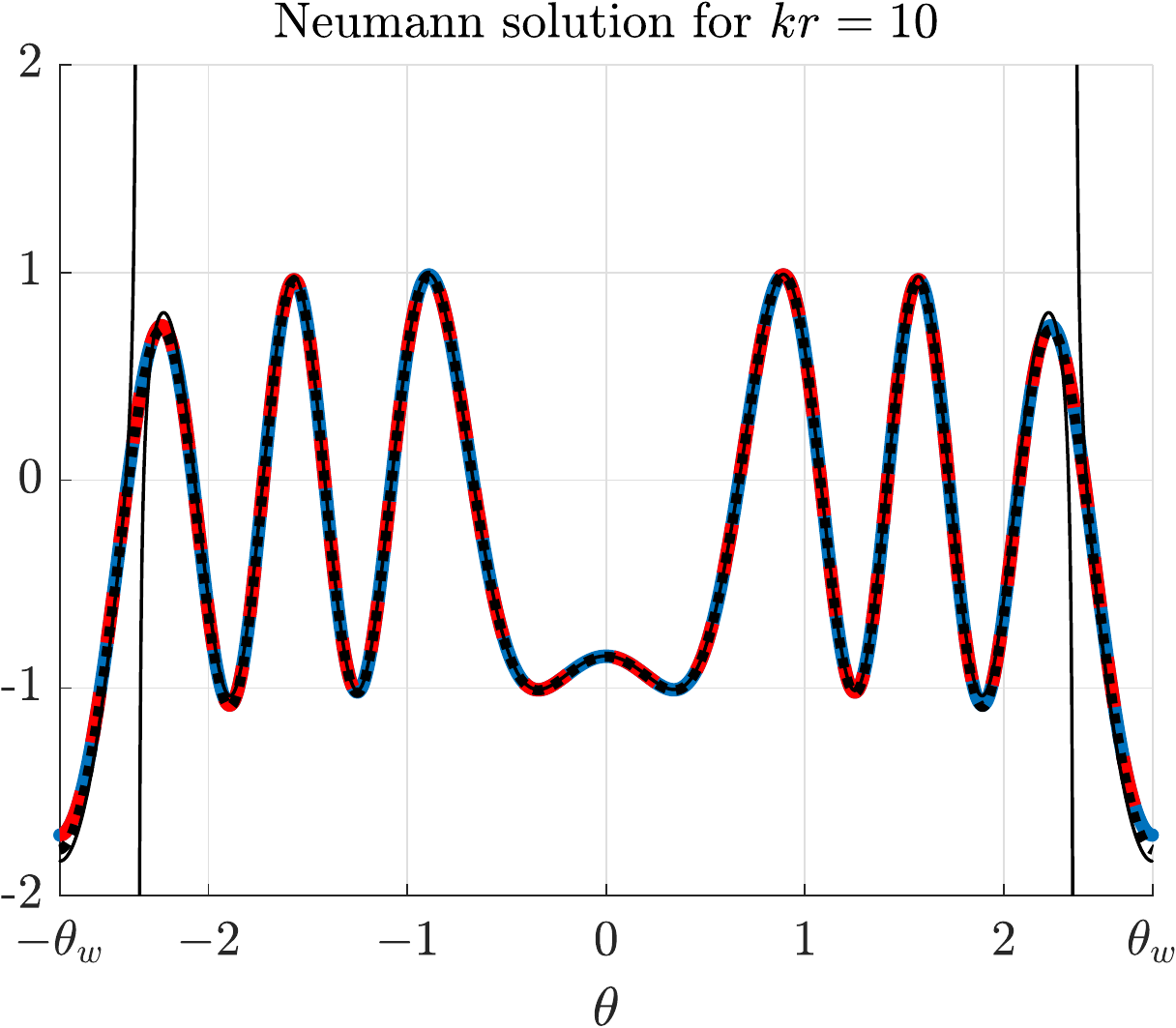}
\caption{Comparison between the real part of the exact solution (SI/SDC) and various approximations for Dirichlet and Neumann \RED{BCs}, for $kr=1,5,10$ and for a wedge characterised by $\theta_\text{w}=7\pi/8$ and an incident angle $\theta_\text{I}=0$.}
\label{fig:phicompare0}
\end{figure}

\begin{figure}[h!]\centering
\includegraphics[width=0.45\textwidth]{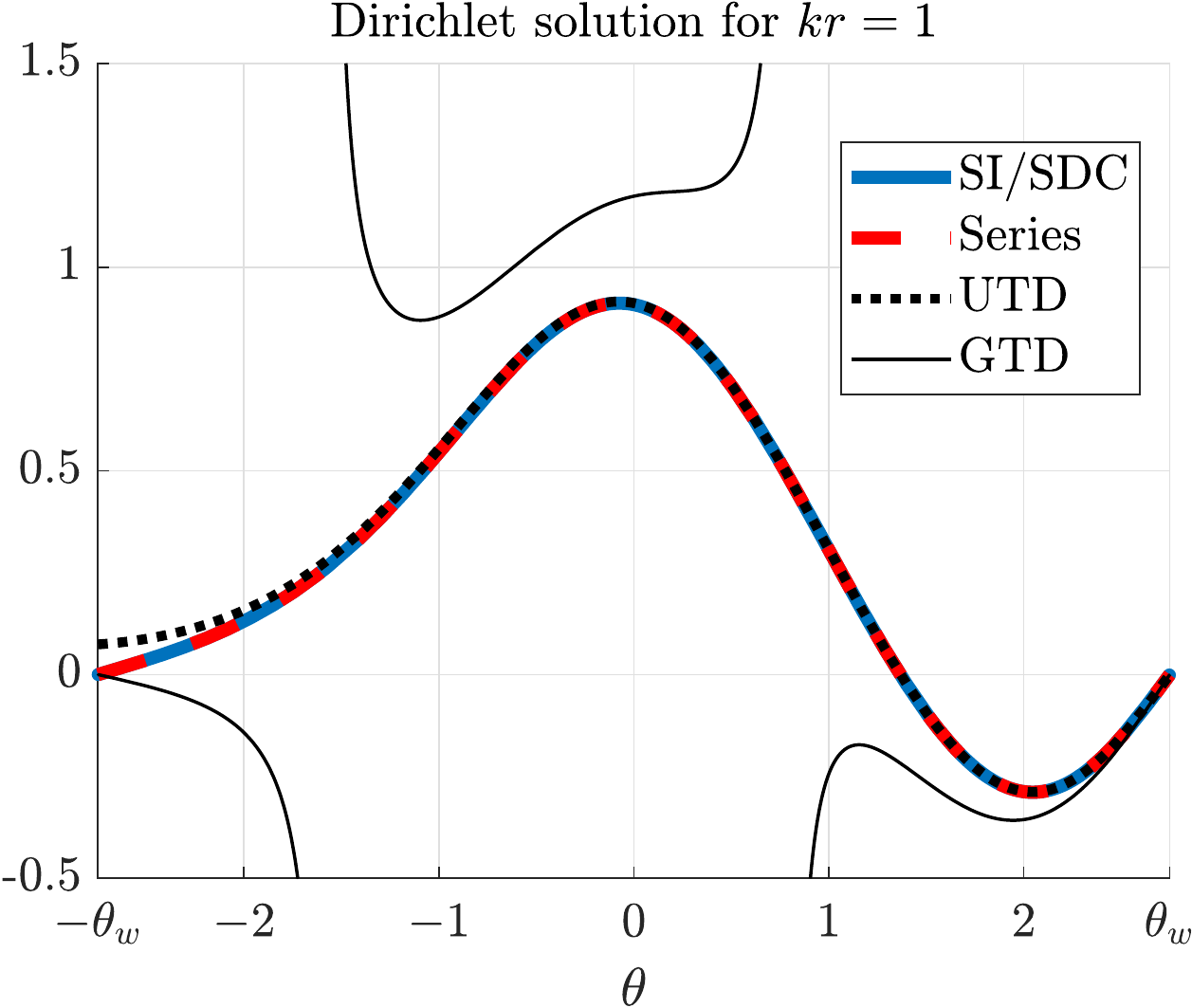}\qquad
\includegraphics[width=0.45\textwidth]{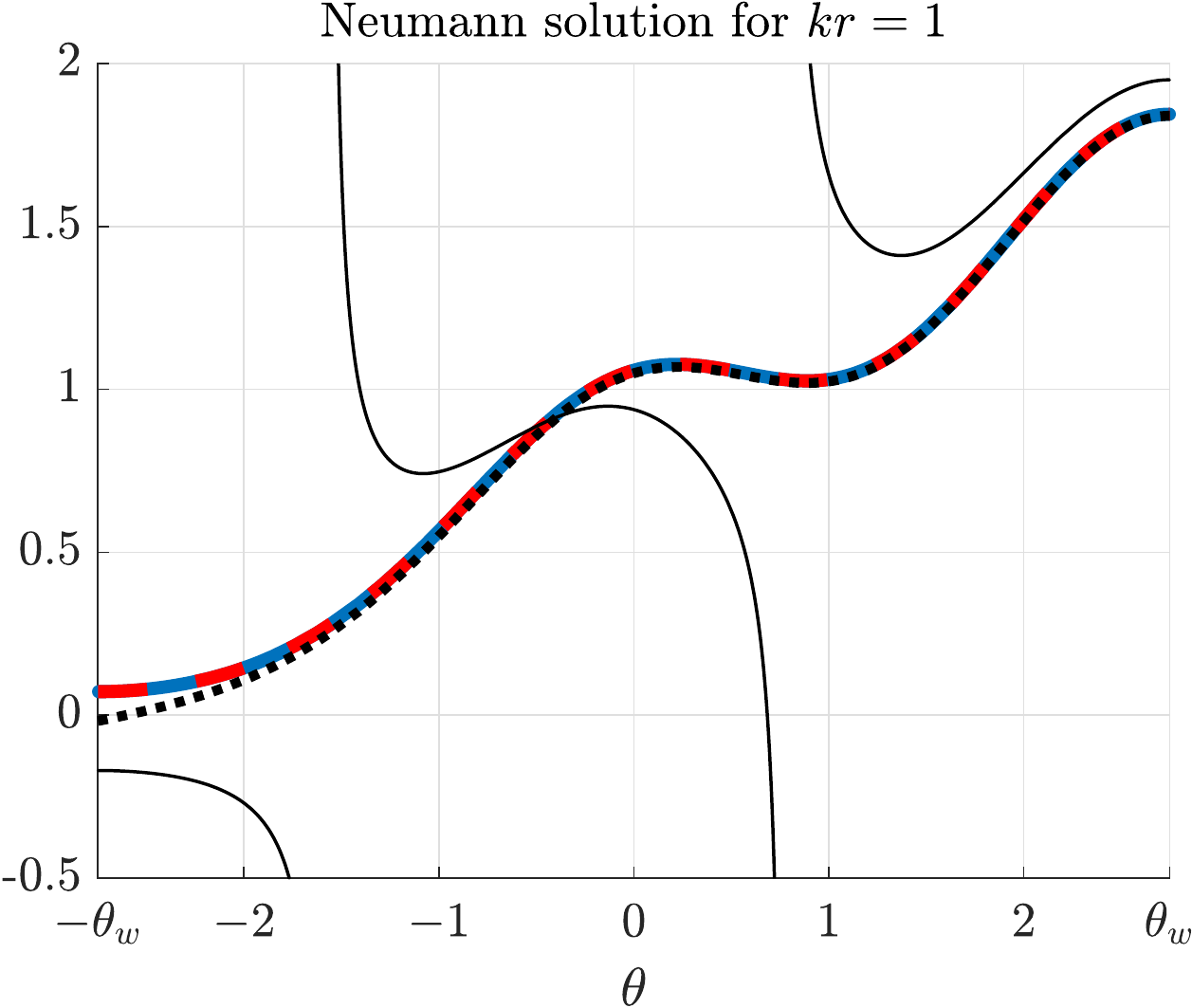}\\
\includegraphics[width=0.45\textwidth]{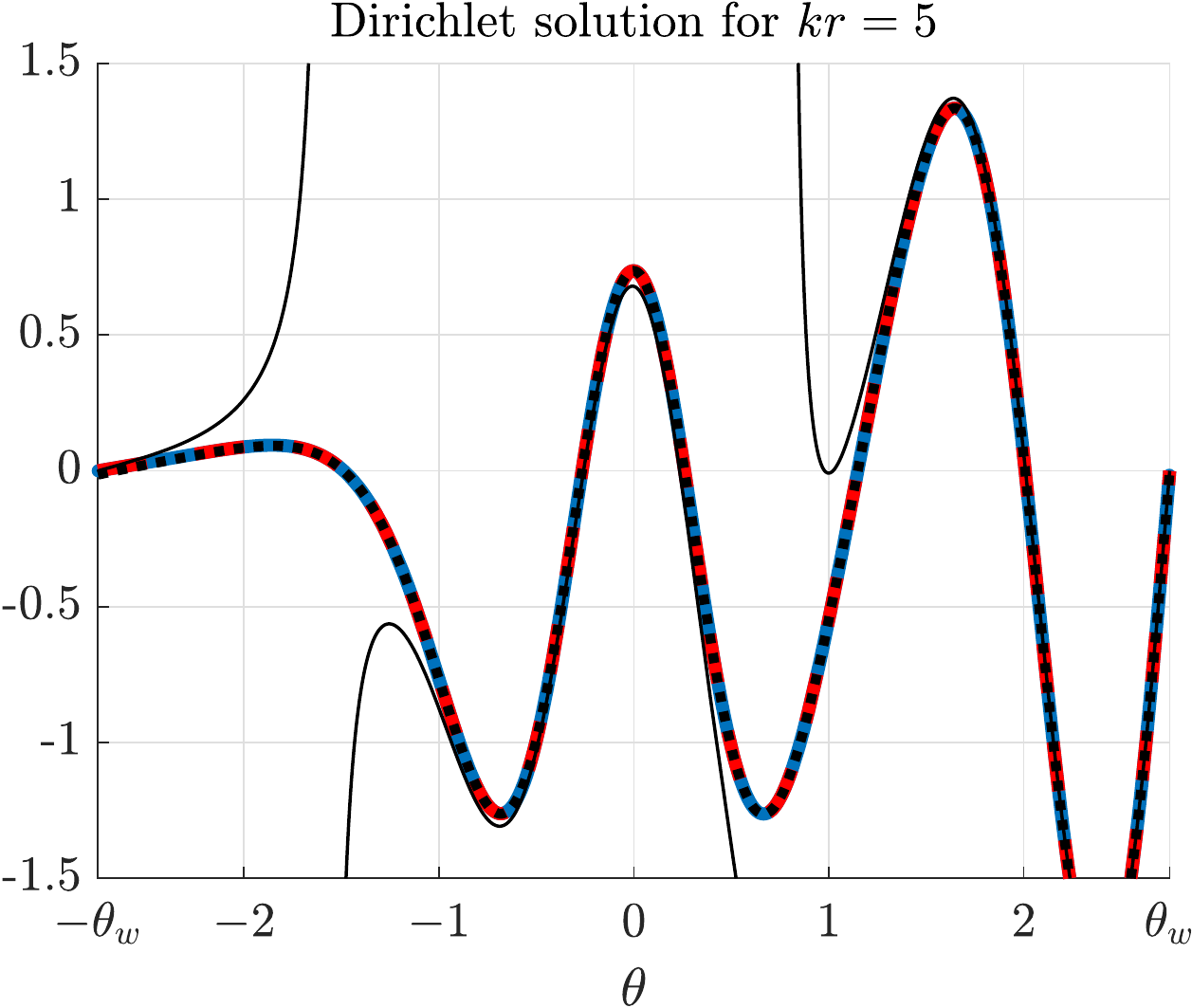}\qquad
\includegraphics[width=0.45\textwidth]{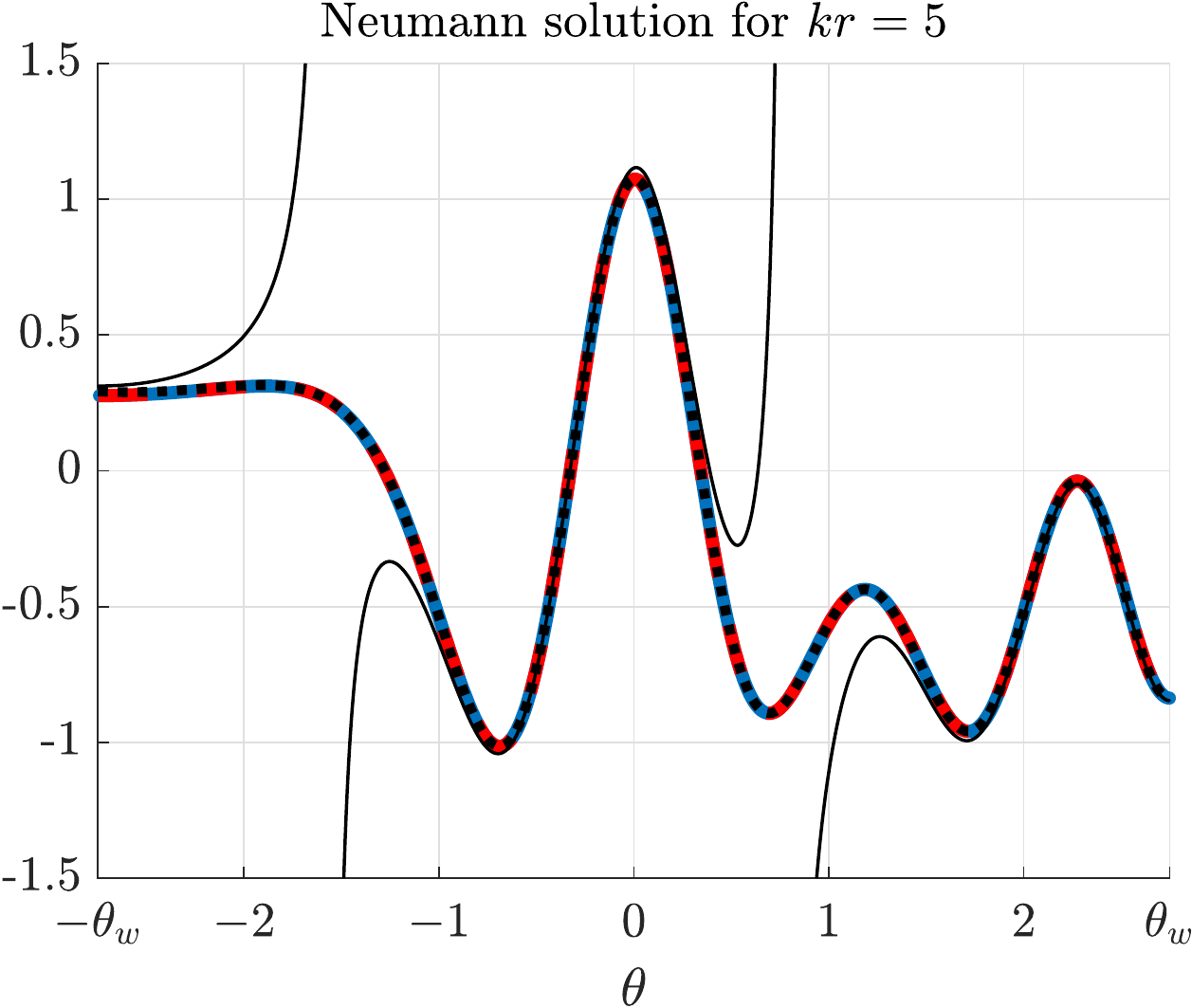}\\
\includegraphics[width=0.45\textwidth]{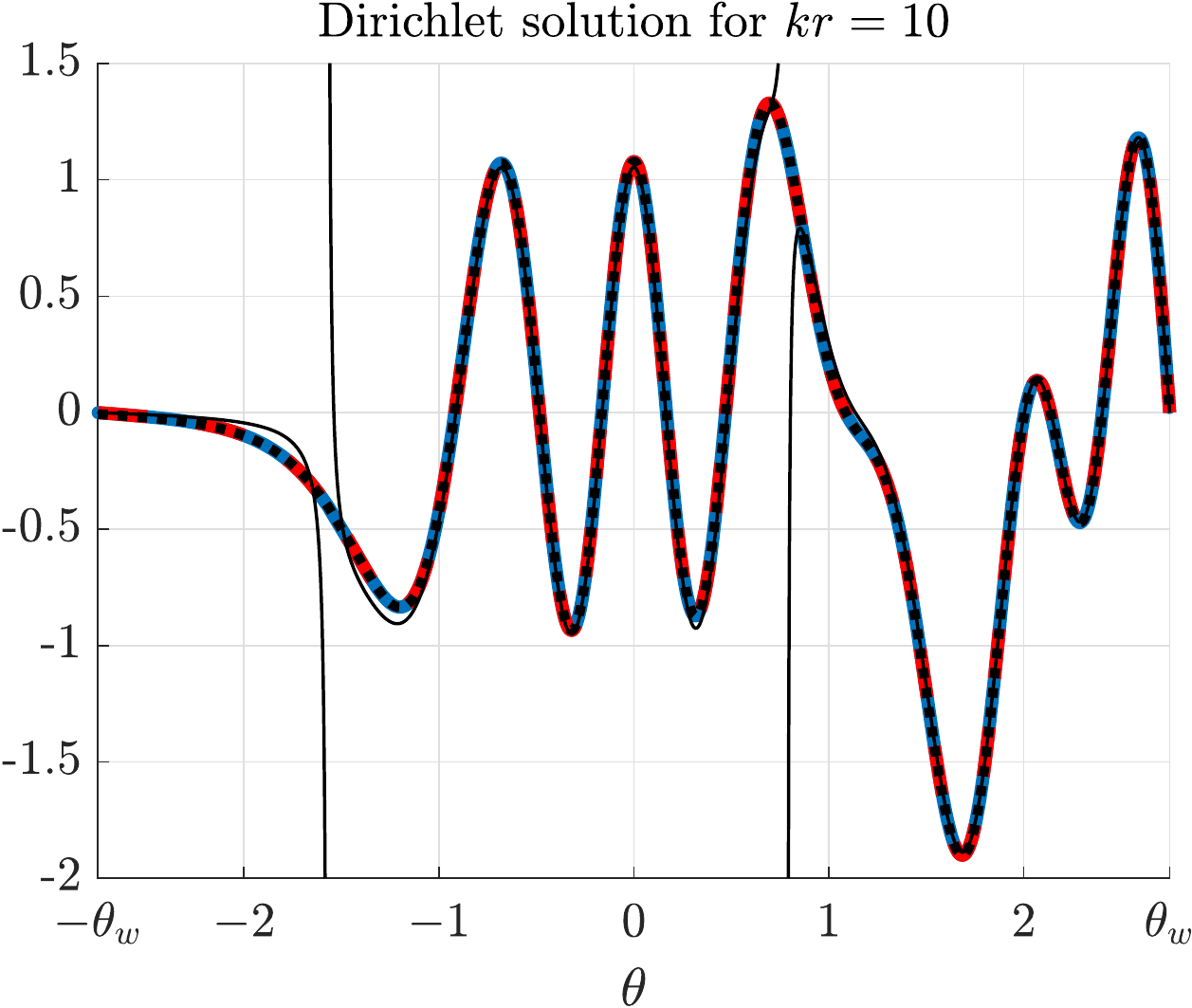}\qquad
\includegraphics[width=0.45\textwidth]{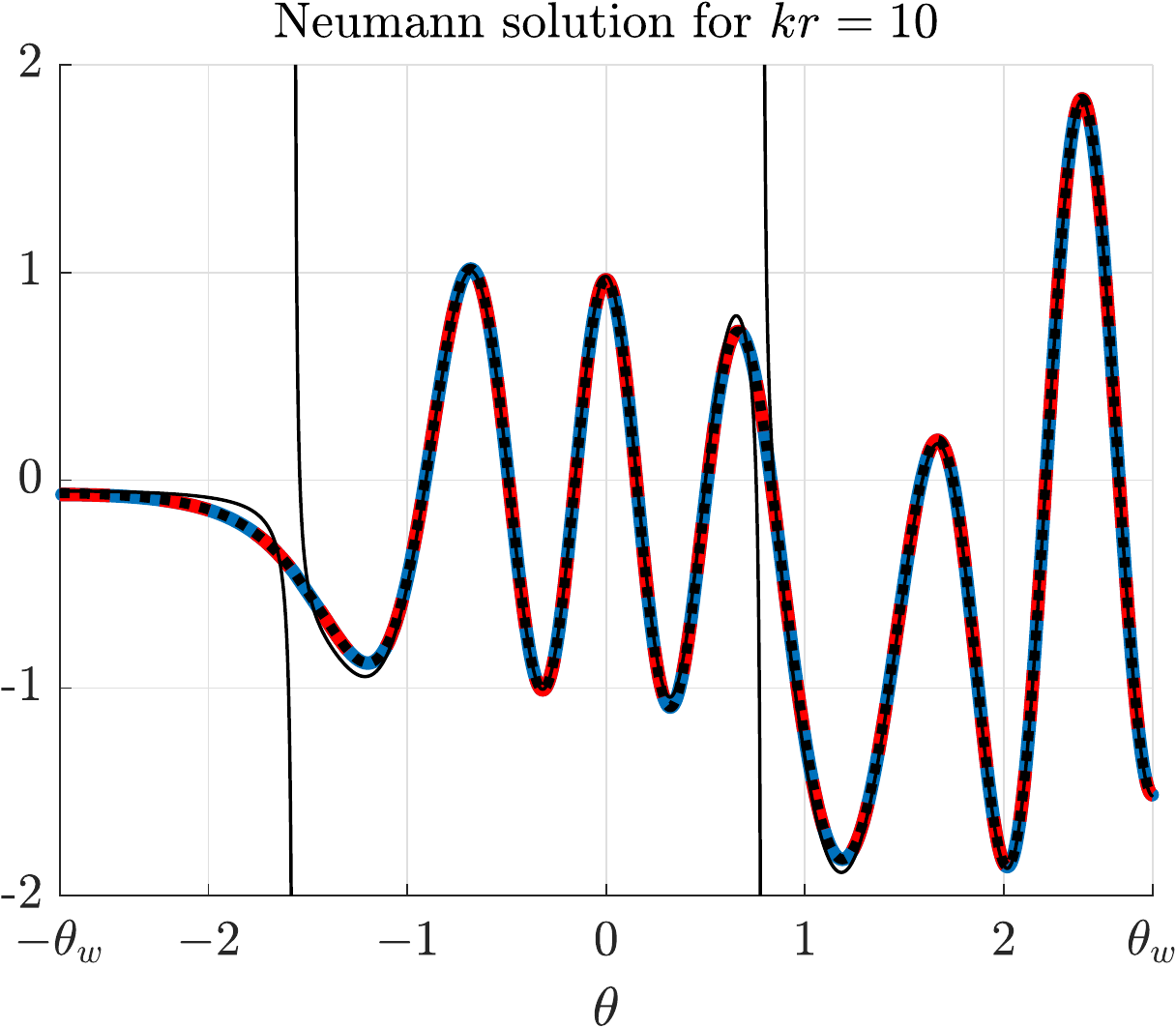}
\caption{Comparison between the real part of the exact solution (SI/SDC) and various approximations for Dirichlet and Neumann \RED{BCs}, for $kr=1,5,10$ and for a wedge characterised by $\theta_\text{w}=7\pi/8$ and an incident angle $\theta_\text{I}=\pi/2$.}
\label{fig:phicomparehalfpi}
\end{figure}

In both Figures \ref{fig:phicompare0} and \ref{fig:phicomparehalfpi}, we confirm that,
\begin{itemize}
\item{The series solution is very accurate despite the truncation. However if we want to consider larger values of $kr$, more terms will be required to remain accurate, which will slow down its computation.}
\item{The GTD approximation has the least overall accuracy and becomes invalid when $\theta$ is close to any GO discontinuities $\theta_\text{I}-\pi$, $2\theta_{\text{w}}-\theta_\text{I}-\pi$ and $-2\theta_{\text{w}}-\theta_\text{I}+\pi$. It does however satisfy the correct \RED{BCs}.}
\item{The \RED{UTD} approximation is a clear improvement to the standard GTD approximation away from the boundaries, in particular it does not have any singularities, but fails to satisfy the \RED{BCs}.}
\end{itemize}

\RED{Both the GTD and UTD approximations appear to improve their accuracy as $kr$ gets larger. To show this, we take the Dirichlet case with $\theta_\text{I}=0$ and look at the quantities $\text{GTD Error}=|\eqref{SMTSommIntDirichlet}-\eqref{SA-GTD}|$ and $\text{UTD Error}=|\eqref{SMTSommIntDirichlet}-\eqref{SA-UTD1st-sol}|$ against $\theta$ for $kr=1,5,10,25$. Figure \ref{fig:SA-error} (left) illustrates the GTD error and shows that it is a good approximation, provided that $kr$ is large enough and $\theta$ is not too close to one of the singular angles $\theta=2\theta_{\text{w}}-\theta_\text{I}-\pi$ and $-2\theta_{\text{w}}-\theta_\text{I}+\pi$ (which are indicated by a thin vertical dashed line). In Figure \ref{fig:SA-error} (right), it is clear that the UTD error at the boundary decreases significantly as $kr$ increases, rendering it a very good approximation everywhere if $kr$ is large enough. We also reconfirm that the UTD approximation is a large improvement in comparison to the GTD approximation. }

\begin{figure}[h!]\centering
\includegraphics[width=0.45\textwidth]{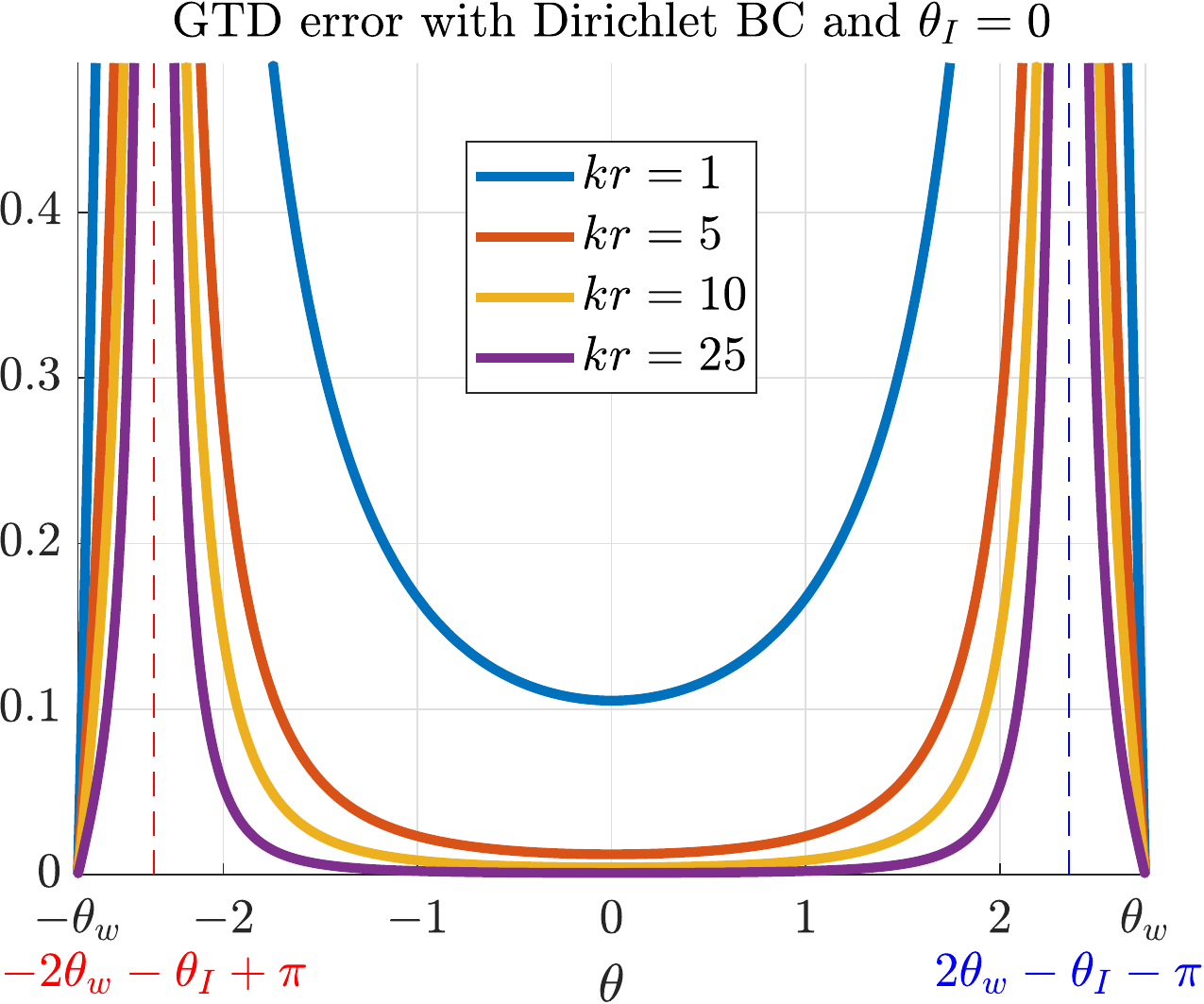}
\includegraphics[width=0.45\textwidth]{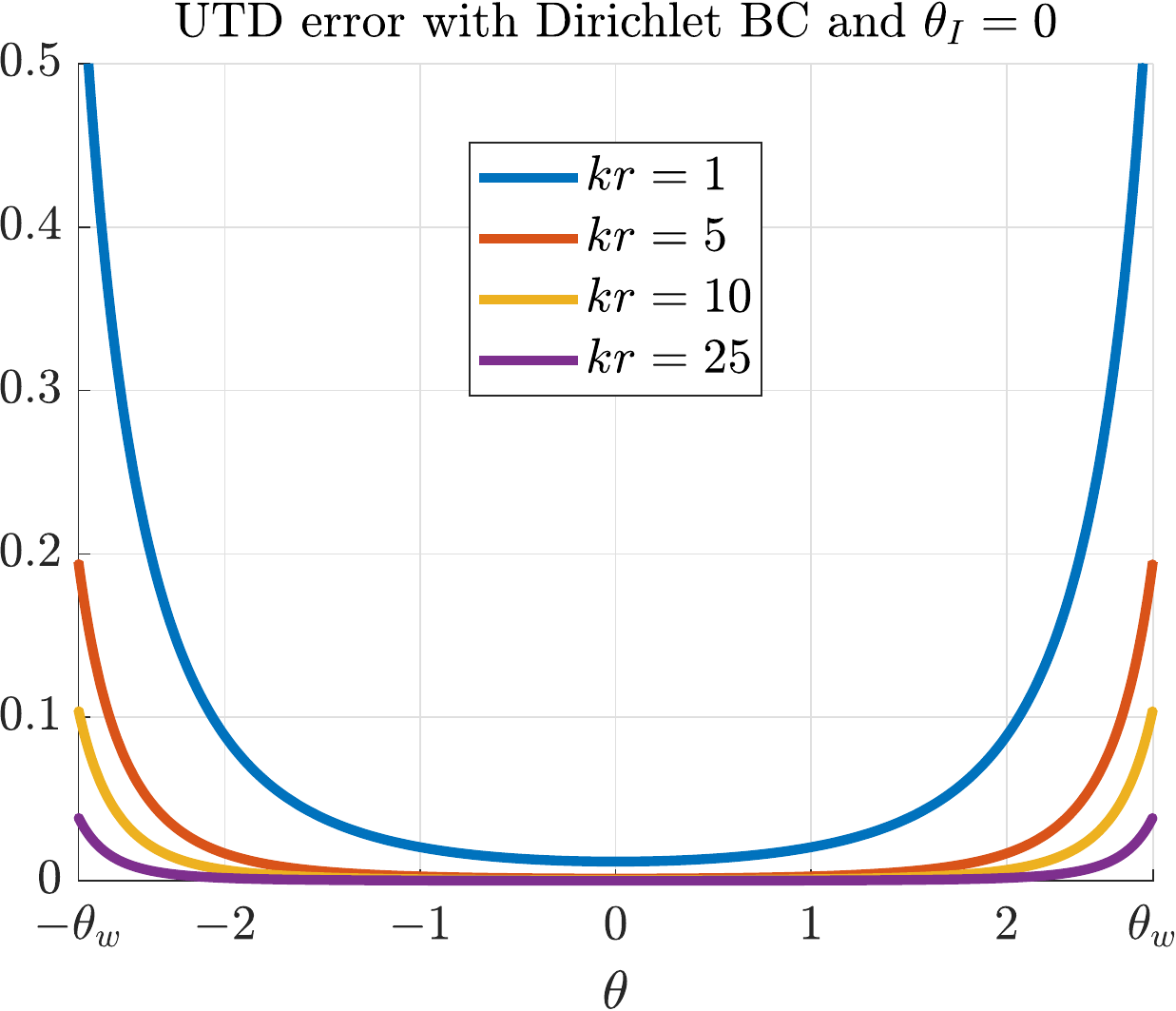}
\caption{\RED{Comparison of the GTD (left) and UTD (right) error for Dirichlet BCs, incident angle $\theta_\text{I}=0$ and increasing values of $kr$,} in the case of a wedge characterised by $\theta_\text{w}=7\pi/8$.}
\label{fig:SA-error}
\end{figure}

Finally, for completeness, we replicate some plots from existing literature using the \RED{UTD} approximation. \RED{Specifically, we replicate the first and last plots of} figure 5 in \citep{HSU2011} which is a comparison of \RED{an alternate definition for \eqref{SA-GammaInt2},} the series solution with 100 terms and a similar \RED{UTD approximation. In order to replicate these plots, we need to adapt to their geometric configuration by making} the substitutions $\theta=\theta_\text{w}-\hat{\theta}$ and $\theta_\text{I}=\theta_\text{w}-\hat{\theta}_\text{I}$. \RED{We use \eqref{SA-UTD2nd-sol} with $\bar{\theta}_{\text{w}}=\pi/36$ and $kr=10\pi$. Figure \ref{fig:HSUPapersCompare} (left) is the Dirichlet case with $\hat{\theta}_\text{I}=\pi/2$. Figure \ref{fig:HSUPapersCompare} (right) is the Neumann case with $\hat{\theta}_\text{I}=2\pi/3$.}

\begin{figure}[h!]\centering
	\includegraphics[width=0.95\textwidth]{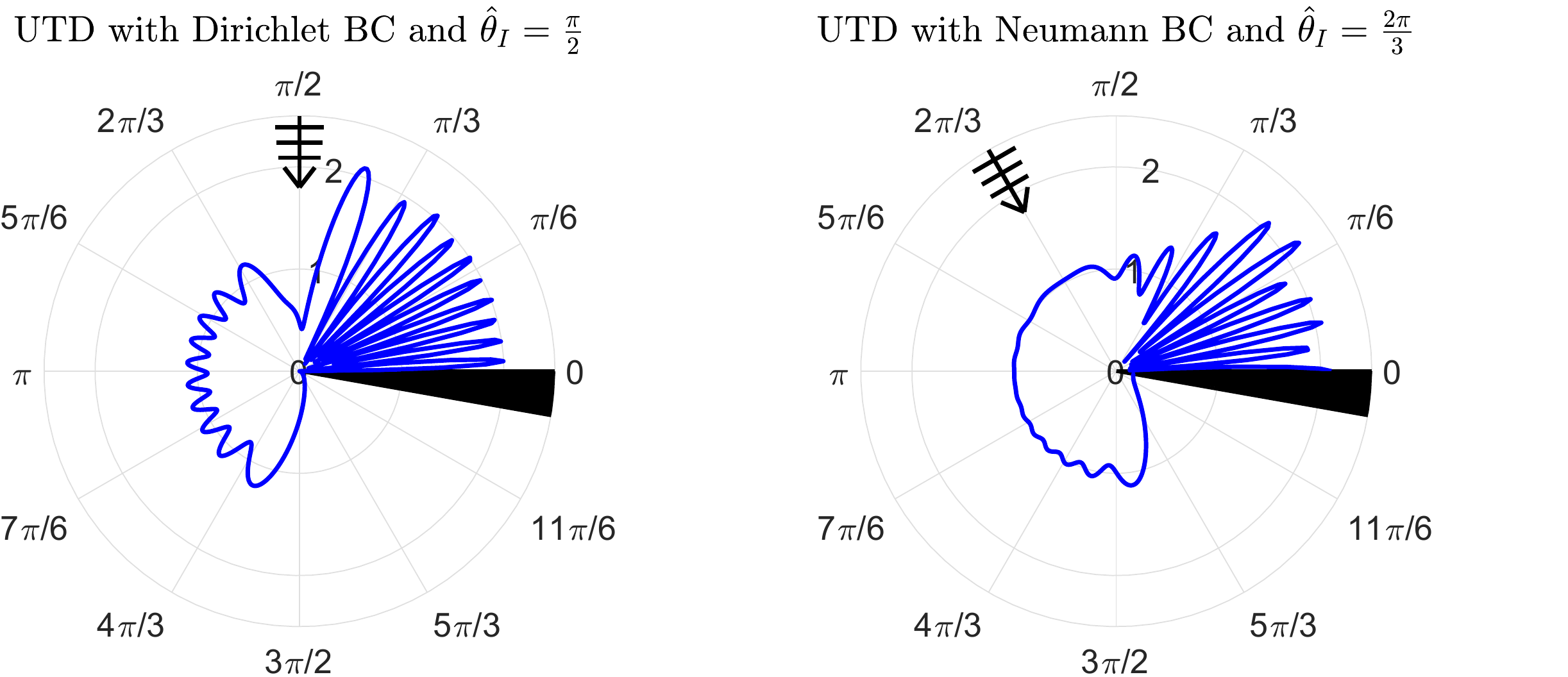}
	\caption{\RED{Replication of the top left and bottom right plots of figure 5 in \citep{HSU2011} using \eqref{SA-UTD2nd-sol}.}}
	\label{fig:HSUPapersCompare}
\end{figure}


\section{Alternative methods}\label{Alt}
Sections \ref{SMT}, \ref{WHT} and \ref{MKLT} cover methods that are most commonly used in diffraction theory. In this section, we will briefly present three alternative methods that have been tailored to tackle the perfect wedge problem.


\subsection{Embedding Formula technique}\label{Embedding} 

The \RED{first} method to be reviewed is based on the idea of embedding. This idea is relatively new in diffraction theory \citep{williams}, and has mainly been used for planar cracks and slits, as well as parallel combinations of these \citep{Gautesen1983,Martin1983,Biggs2001,Biggs2002}. Though, recently, in \citep{ShaninCraster2005} it was adapted to wedges with \RED{a} rational angle. We will here attempt to summarise this method and consider again our wedge region characterised by $\theta_{\text{w}}$. We seek the total field $\Phi$ satisfying the Helmholtz equation \eqref{Intro-Helmholtz}, subjected to Dirichlet \eqref{Intro-DBC} or Neumann \eqref{Intro-NBC} \RED{BCs, as well as radiation and edge conditions \eqref{Intro-2DSRcondition} and \eqref{Intro-Edge}} for a plane wave incidence $\Phi_\text{I} = e^{- ikr \cos (\theta - \theta_\text{I})}$, with incident angle $\theta_\text{I}$. The aim of the method is to recover the diffraction coefficient of the diffracted field $\Phi_{\text{Diff}}$.

\paragraph{The diffraction coefficient}Using classical separation of variables in the polar coordinates $(r,\theta)$ and the edge conditions, it can be shown that the total field $\Phi$ has an eigenfunction expansion of the form}
\begin{align}
\Phi(r,\theta)&=\sum_{m=0}^{\infty}(2/k)^{\nu_m}\Gamma(1+\nu_m)K_m(\theta_\text{I})u_m(r,\theta),\label{eq:EFeigenfunctionexpansion}
\end{align}
\mylinenum{where $\nu_m = m \delta = m \pi / \left( 2 \theta_{\text{w}} \right)$ and $u_m$ is a product of Bessel functions $J_{\nu_m} (kr)$ and some trigonometric functions of $\theta$ satisfying the \RED{BCs}\footnote{The multiplicative factor $(2 / k)^{\nu_m} \Gamma (1 + \nu_m)$ is just here to compensate the near-field behaviour of the Bessel functions, and, doing so, somehow normalise the expansion.}. In the Dirichlet case, the $m=0$ term in the sum is equal to zero. Note that using the series results (\ref{MKLT-SeriesSol-D}) and (\ref{MKLT-SeriesSol-N}) of Macdonald type, we can recover $K_m$ exactly, but we will not use this here.
The aim is to determine the diffraction coefficient $D (\theta,\theta_\text{I})$, already defined in (\ref{SA-diffraction-coefficient}), that is such that}
\begin{align}
\Phi_{\text{Diff}}(r,\theta)&\underset{r\rightarrow\infty}{\sim}D(\theta,\theta_\text{I})\frac{e^{ikr}}{\sqrt{kr}} 
\end{align}
\mylinenum{\paragraph{The edge Green's functions}In order to do this, as is customary
with embedding, we need to introduce an auxiliary problem\footnote{\RED{Here the auxiliary problems will be constructed from point sources. However, another type of embedding formulae can be obtained with plane wave auxiliary problems, see \citep{Biggs2006} for example.}}. In fact here, we will introduce infinitely many of them. Let $m \in \mathbb{N} \backslash \{0\}$, and consider the function $\widehat{u}_{m, \varepsilon}$ that is the tailored
Green's function (i.e. that satisfies the \RED{BCs}) for the Helmholtz equation resulting from $m$ point sources given by $z_j=\varepsilon e^{i \left( \varphi_j - \theta_{\text{w}} \right)}, j = 1, \ldots,m$, where \RED{$\varphi_j=(2j-1)\theta_{\text{w}}/m$ for Dirichlet BC and $\varphi_j=2j\theta_{\text{w}}/m$ for Neumann BC. The} strength $A_j$ of each source is given by $A_j=(-1)^j\pi\varepsilon^{-\nu_m}$, as illustrated in Figure \ref{fig:EF-sources}.

\begin{figure}[h!]\centering
\includegraphics[width=0.3\textwidth]{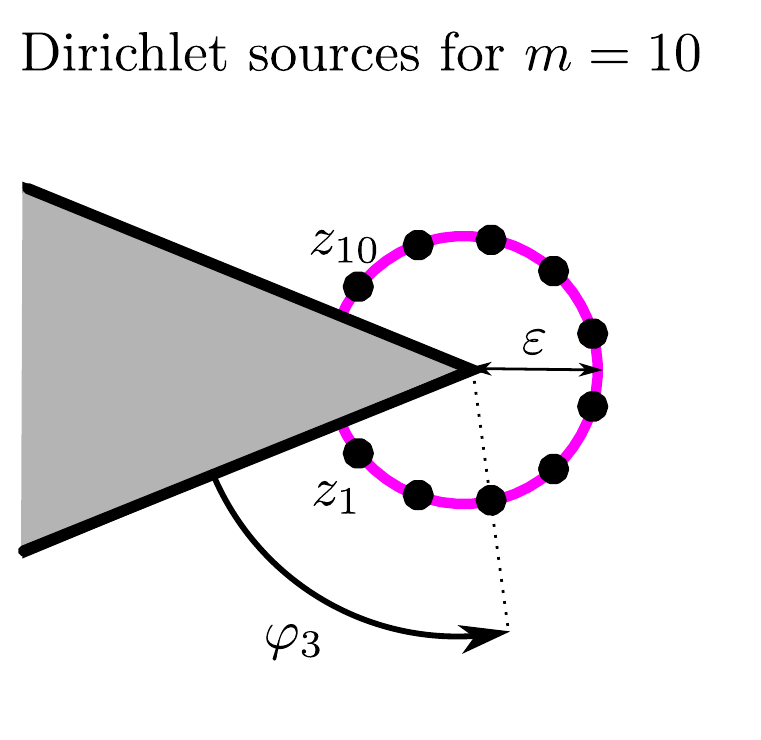}
\qquad \includegraphics[width=0.3\textwidth]{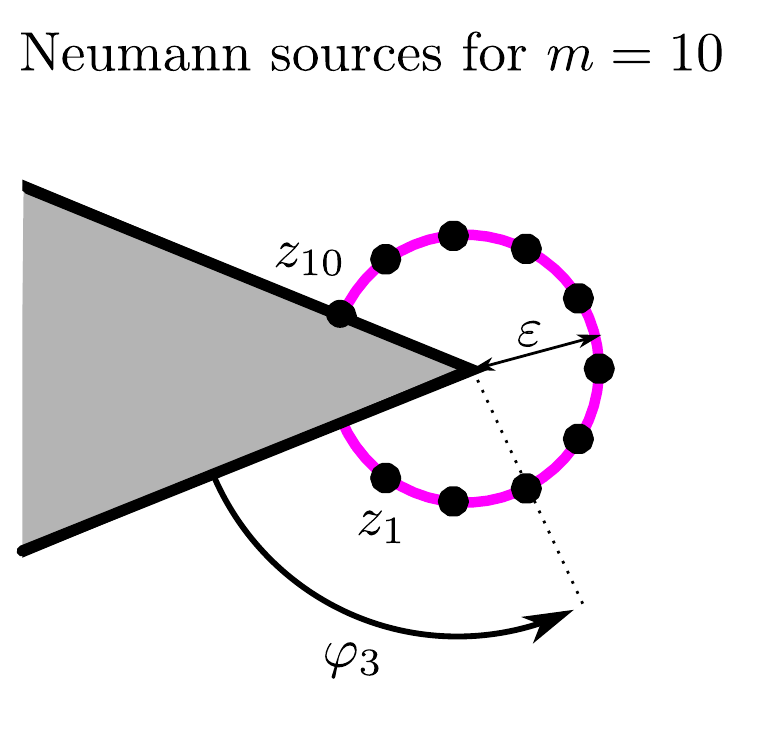}
\caption{Position of the sources in Dirichlet and Neumann cases for $m = 10$}
\label{fig:EF-sources}
\end{figure}

The $m$th edge Green's function is then defined by}
\begin{align}
\widehat{u}_m&=\lim_{\varepsilon\rightarrow0}\widehat{u}_{m,\varepsilon} 
\end{align}
\mylinenum{The near field behaviour of the edge Green's function can be studied by considering $\widehat{u}_{m, \varepsilon}$ for fixed $\varepsilon$, close to the wedge edge. In that vicinity, we can scale the space variables to show that $\widehat{u}_{m, \varepsilon}$ behaves locally like $\widehat{u}_{m,\varepsilon}^{\tmop{inner}}$, which is the exact same Green's function, but for Laplace's equation instead of Helmholtz. Using the method of images in a half-space, and the mapping $z = w^{1 / \delta}$, it is possible to find $\widehat{u}_{m, \varepsilon}^{\tmop{inner}}$ explicitly\footnote{Note that in {\citep{ShaninCraster2005}}, only the Dirichlet formula is given, and is slightly different from this one (the factors $i$ are missing), which we think is a typographical error. } as}
\RED{\begin{align*}
\widehat{u}_{m,\varepsilon}^{\tmop{inner}}=-\frac{\varepsilon^{-\nu_m}}{2}\text{Re}\left\{\ln\left(\frac{Z^{\nu_m}-a\varepsilon^{\nu_m}}{Z^{\nu_m}+a\varepsilon^{\nu_m}}\right)\right\} \text{ with } \left\{
\begin{array}{l}
a=i \text{ for Dirichlet BC} \\
a=1 \text{ for Neumann BC}
\end{array}
\right.
\end{align*}}
\mylinenum{where $Z=re^{i \varphi}$, $\varphi$ being the angle measured from the bottom face of the wedge, so that we have $\varphi=\theta+\theta_{\text{w}}$. Looking at the leading order of $\widehat{u}_{m,\varepsilon}^{\tmop{inner}}$ as $\varepsilon\rightarrow 0$, using the fact that $\RED{\ln}(z)\underset{z\rightarrow 1}{\sim}1-1/z$, we get}
\begin{eqnarray}
\widehat{u}_{m,\varepsilon}^{\tmop{inner}}&\underset{\varepsilon\rightarrow0}{\sim}&r^{-\nu_m}\left\{
\begin{array}{cl}
\sin(\nu_m\varphi)&\text{for Dirichlet BC}\\
\cos(\nu_m\varphi)&\text{for Neumann BC}
\end{array}\right.,
\label{EF-edge-green-edge-behaviour}\end{eqnarray}
\mylinenum{which by construction, is also the local behaviour of $\widehat{u}_m$ near the wedge edge. Note that the edge Green's function is singular on the wedge edge and does not satisfy the edge condition, we say that it is {\tmem{oversingular}}. It does however satisfy the Helmholtz equation everywhere outside the wedge. This leads to the exact representation of $\widehat{u}_m$ as:}
\begin{eqnarray}
\widehat{u}_m(r,\theta)&=&\frac{\pi i}{\Gamma(\nu_m)}(k/2)^{\nu_m}H_{\nu_m}^{(1)}(k_0r)\left\{
\begin{array}{cl}
\sin\left(\nu_m\left(\theta+\theta_{\text{w}}\right)\right)&\text{for Dirichlet BC}\\
\cos\left(\nu_m\left(\theta+\theta_{\text{w}}\right)\right)&\text{for Neumann BC}
\end{array}\right.,
\label{eq:exact-edge-green}\end{eqnarray}
\mylinenum{since it is clear that the above expression has the right type of singularity, and satisfies the boundary and radiation conditions, as well as the Helmholtz equation.

It is also natural to define the directivity $\widehat{D}_m (\theta)$ for each edge Green's function by}
\begin{eqnarray}
\widehat{u}_m(r,\theta)&\underset{r\rightarrow\infty}{\sim}\widehat{D}_m(\theta)\frac{e^{ikr}}{\sqrt{kr}},
\label{EF-def-directivity}\end{eqnarray}
\mylinenum{and using the asymptotic behaviour of the Hankel function for large argument, \eqref{eq:exact-edge-green} and (\ref{EF-def-directivity}) imply that}
\begin{eqnarray}
\widehat{D}_m(\theta)&=\frac{\sqrt{2\pi}(k/2)^{\nu_m}e^{-\frac{i\nu_m\pi}{2}}}{\Gamma(\nu_m)}\left\{
\begin{array}{cl}
\sin\left(\nu_m\left(\theta+\theta_{\text{w}}\right)\right)&\text{for Dirichlet BC}\\
\cos\left(\nu_m\left(\theta+\theta_{\text{w}}\right)\right)&\text{for Neumann BC}
\end{array}\right. .
\label{EF-exact-directivity}\end{eqnarray}
\mylinenum{It is important to note the main difference between the directivities $\widehat{D}_m (\theta)$ and the diffraction coefficient $D (\theta, \theta_\text{I})$: the former only depends on one angular variable, while the latter depends on two. Remarkably, using the reciprocity principle, it is possible to relate the far-field of the edge Green's functions to the near-field of each components of the eigenfunction expansion \eqref{eq:EFeigenfunctionexpansion} as follows:}
\begin{eqnarray}
\widehat{D}_m(\theta_\text{I})&=\frac{m\pi}{2}K_m(\theta_\text{I})
\end{eqnarray}
\mylinenum{\paragraph{The operator $\textbf{H}_p$} As mentioned above, this method can only be applied to rational angles\footnote{\RED{\citet{ShaninPseudoDiff} have extended this work by considering a pseudo-differential operator $\textbf{K}_{\mu}$ instead of the differential operator $\textbf{H}_p$. Note that for an integer $n$, $\textbf{K}_n$ reduces to $C_n T_n \left( \frac{i}{k} \frac{\partial}{\partial x} \right)$ for some constant $C_n$, which establishes the link with the theory developed here. This new operator can however be used when $\mu \nin \mathbb{N}$ to produce an embedding formula valid for wedges with non-rational angles, though it cannot be used for polygons.}}, so let us set $2\theta_{\text{w}}=\frac{q\pi}{p}$ for some positive integers $p$ and $q$. Now define, the operator $\textbf{H}_p$ as follows:}
\begin{align}\nonumber
\textbf{H}_p=(-ik)^p\left[T_p\left(\frac{i}{k}\frac{\partial}{\partial x}\right)-T_p(\cos(\theta_\text{I}))\right],
\end{align}
\mylinenum{where $T_p$ is the $p$th Chebyshev polynomial, and it is understood that for some integer $n$, $\left(a\frac{\partial}{\partial x}\right)^n=a^n\frac{\partial^n}{\partial x^n}$. 
From now on, for brevity, we will focus solely on the Dirichlet case. It is relatively easy to show that for every $m\in\mathbb{N}\backslash\{0\}$, $\textbf{H}_p[u_m]$ satisfies the Helmholtz equation, the correct boundary conditions and the radiation condition, and that $\textbf{H}_p [\Phi_\text{I}] = 0$. It is also possible to prove (though it is more difficult) that}
\begin{align}
\textbf{H}_p[\Phi]\underset{r\rightarrow 0}{\sim}&\ 2^{p-1}(-1)^{q-p+1}\sum_{m=1}^{q-1}K_m(\theta_\text{I})\nu_m(\nu_m-1)\ldots(\nu_m-p+1)r^{-\nu_{q-m}}\sin(\nu_{q-m}\varphi)\nonumber\\
&+\text{terms that satisfy the edge conditions}
\label{EF-exp-H-Phi}\end{align}
\mylinenum{We refer to {\citep{ShaninCraster2005}} for the details of the proof, but it relies on a careful analysis of the near-field and far-field behaviour of $\textbf{H}_p[u_m]$. It also uses the identity $\nu_m \pm p = \nu_{m \pm q}$, which explains how $q$ enters the scene.

\paragraph{Embedding formula} Note now that the behaviour of each term in (\ref{EF-exp-H-Phi}) reminds of that of the $(q - m)$th edge Green's function (see (\ref{EF-edge-green-edge-behaviour})). This motivates the introduction of the auxiliary function}
\begin{align}\nonumber
W=\textbf{H}_p[\Phi]-2^{p-1}(-1)^{q-p+1}\sum_{m=1}^{q-1}K_m(\theta_\text{I})\nu_m(\nu_m-1)\ldots(\nu_m-p+1)\widehat{u}_{q-m}. 
\end{align}
\mylinenum{By construction, $W$ satisfies the edge condition, and it is also clear from what has been done above, that it satisfies the Helmholtz equation, the boundary and the radiation conditions. Hence, by uniqueness, we conclude that $W \equiv 0$, and we obtain the \RED{\textit{weak} form of the embedding} formula}
\begin{eqnarray}
\textbf{H}_p[\Phi]=2^{p-1}(-1)^{q-p+1}\sum_{m=1}^{q-1}K_m(\theta_\text{I})\nu_m(\nu_m-1)\ldots(\nu_m-p+1)\widehat{u}_{q-m},
\label{EF-weak-EF}\end{eqnarray}
\mylinenum{valid everywhere, that relates the total field $\Phi$ to the edge Green's functions. Focusing now on the far-field, (\ref{EF-weak-EF}) makes it possible to express the diffraction coefficient $D(\theta,\theta_\text{I})$ in terms of the directivities of some of the edge Green's functions, as summarised in the equation below:}
\begin{eqnarray}
D(\theta,\theta_\text{I})=\sum_{m=1}^{q-1}\frac{(-1)^{q-p+1}\nu_m(\nu_m-1)\ldots(\nu_m-p+1)}{m\pi(ik/2)^p(\cos(p\theta)-(-1)^p\cos(p\theta_\text{I}))}\widehat{D}_m(\theta_\text{I})\widehat{D}_{q-m}(\theta)
\label{EF-strong-EF}\end{eqnarray}
\mylinenum{The formula (\ref{EF-strong-EF}) is the main result of {\citep{ShaninCraster2005}} and is referred to as the {\tmem{Embedding formula}}. It is remarkable in the sense that it allows to express the diffraction coefficient, depending on two angular variables, in terms of a sum of products of simpler directivities depending on one angular variable only. Moreover, in that case, thanks to (\ref{EF-exact-directivity}), we know the directivities exactly and we can then recover a new analytical expression for
the diffraction coefficient. For a given rational angle, it is possible to show that it is indeed equal to that given in (\ref{SA-diffraction-coefficient}).


\RED{
\paragraph{Critical analysis} The concept of embedding is very general in diffraction theory, which makes this method very adaptable to all kinds of geometries such as slits, wedges, plane sectors, cubes \citep{wedgecraster} and curved geometries \citep{Moran2016}. In that respect, instead of being seen as a method, one can consider the embedding structure as an inherent property of diffraction problems.

Its main advantage is that once derived explicitly, the embedding formula of a given diffraction problem allows one to obtain a very efficient way of computing the diffraction coefficient resulting from a incident plane wave for all observation and incident angles.

Even though the weak embedding formulae of the type (\ref{EF-weak-EF}) are valid everywhere, the power of embedding formulae only becomes apparent in the far-field, where it can be written in its strong form (see (\ref{EF-strong-EF}) for the present case). In that sense, such formula is not particularly helpful to shed some light on the near field behaviour of diffraction problems. Though because of this emphasis on the far-field, one can consider structures with multiple diffracting parts such as polygons for example.
}


\subsection{Random Walk method}\label{RWM} 

This method, developed in a series of papers {\citep{BudaevBogy2001,BudaevBogy2002a,BudaevBogy2002b}}, and applied to the wedge problem in {\citep{BudaevBogy2003}}, is based on the Feynman-Kac formula (see e.g. {\citep{Feynman1948}},
{\citep{Kac1949}} and {\citep{Freidlin1985}}). This formula implies, in particular, that the solution $u$ of a deterministic PDE on a domain $\Omega$ with Dirichlet condition $ u|_{_{\partial \Omega}} = f (r, \theta)$ on the boundary $\partial \Omega$}
\begin{eqnarray}
\frac{\sigma_1^2(r,\theta)}{2}\frac{\partial^2 u}{\partial r^2}+\frac{\sigma_2^2(r,\theta)}{2}\frac{\partial^2 u}{\partial \theta^2}+A_1(r,\theta)\frac{\partial u}{\partial r}+A_2(r,\theta)\frac{\partial u}{\partial\theta}+B(r,\theta)u=0
\label{eq:rw0}\end{eqnarray}
\mylinenum{with real-valued coefficients $\sigma_{1, 2}$, $A_{1, 2}$ and $B$, can be written as}
\begin{eqnarray}
u(r,\theta)=\textbf{E}\left\{f(\xi_{\tau}^1,\xi_{\tau}^2)e^{\int_0^{\tau}B(\xi_s^1,\xi_s^2)\mathd s}\right\},
\end{eqnarray}
\mylinenum{where $\textbf{E}$ represents the mean operator, and $\xi_t^{1, 2}$ are random motions governed by the two coupled stochastic differential equations (SDE) with drift coefficient $A_{1, 2}$ and diffusion coefficient $\sigma_{1,2}$}
\begin{align}
\mathd\xi_t^1=A_1(\xi_t^1,\xi_t^2)\mathd t+\sigma_1 (\xi_t^1,\xi_t^2)\mathd W_t^1 \text{ and } \mathd\xi_t^2=A_2(\xi_t^1,\xi_t^2)\mathd t+\sigma_2(\xi_t^1,\xi_t^2)\mathd W_t^2,
\label{eq:firstSDEs}\end{align}
\mylinenum{with initial conditions \RED{(ICs)} $\xi_0^1 = r$ and $\xi_0^2 = \theta$, where $W_t^{1,2}$ are Brownian motions (also known as Wiener processes)\footnote{\RED{We do not intend to insist on the rigorous mathematical definitions of these objects here, however we refer the interested reader to general textbooks on the topic, such as \citep[Chapter 6]{Voss_2013} for example, where Brownian motions, SDEs (and their resolution via the Euler-Maruyama scheme) and It{\^o} calculus are introduced.}}. The exit time $\tau$ is the time when each computation should be stopped and it corresponds to the first time $t$ such that $(\xi_t^1, \xi_t^2) \in \partial \Omega$.

If the coefficients in \eqref{eq:rw0} and \eqref{eq:firstSDEs} are complex-valued (which as we will see will be the case for the problem at hand), then the \RED{Feynman}-Kac representation is still valid, but it becomes difficult to determine and define the exit time $\tau$. In fact if the coefficients are complex, then so will be the random motions $\xi_t^{1, 2}$, and since the points of $\partial \Omega$ belong to $\mathbb{R}^2$, we cannot easily characterise the fact that $(\xi_t^1, \xi_t^2)$ hits this boundary. This can be addressed by considering the ``continuation'' of the boundary $\partial \Omega$ into a manifold of real dimension 2 within the space $\mathbb{C}^2$ and by multiplying \eqref{eq:rw0} by $q^2 (r, \theta)$, where $q$ is a complex-valued function. For a suitable function $q$, it becomes possible to define an exit time $\tau$, and the solution to the PDE is given by}
\begin{eqnarray}
  u (r, \theta) & = & \textbf{E} \left\{ f (\xi_{\tau}^1, \xi_{\tau}^2)
  e^{\int_0^{\tau} q^2 (\xi_s^1, \xi_s^2) B (\xi_s^1, \xi_s^2) \mathd s}
  \right\},  \label{eq:solpolar}
\end{eqnarray}
\mylinenum{where $\xi_t^{1, 2}$ are random motions governed by the two coupled stochastic differential equations (SDE)}
\begin{align}\nonumber
\mathd\xi_t^{1,2}=q^2(\xi_t^1,\xi_t^2)A_{1,2}(\xi_t^1,\xi_t^2)\mathd t+q(\xi_t^1,\xi_t^2)\sigma_{1,2}(\xi_t^1,\xi_t^2)\mathd W_t^{1,2}
\end{align}
\mylinenum{with \RED{ICs} $\xi_0^1 = r$ and $\xi_0^2 = \theta$.

In order to fit within this framework, for the wedge problem in {\citep{BudaevBogy2003}}, the authors aim to solve the Helmholtz equation (\ref{Intro-Helmholtz}), subject to the radiation condition and to Dirichlet \RED{BCs} of the type $\Phi\left(r,\pm\theta_{\text{w}}\right)=F\left(r,\pm\theta_{\text{w}}\right)$. They seek a solution of the form $\Phi = u e^{i S}$ for some unknown functions $u$ and $S$. The Helmholtz equation becomes}
\begin{eqnarray}
\Delta u+2i\nabla u\cdot\nabla S+iu\Delta S+u(k^2-\nabla S\cdot\nabla S)=0
\label{eq:rw1}\end{eqnarray}
\mylinenum{If the solution is in the Liouville form, we choose $S (r, \theta) = k r$. In this case, $S$ automatically satisfies the eikonal equation $\nabla S \cdot \nabla S = k^2$ and, after multiplication by $\frac{i}{2 k}$, \eqref{eq:rw1} becomes}
\begin{align}\nonumber
\frac{i}{2k}\left(\frac{\partial^2u}{\partial r^2}+\frac{1}{r^2}\frac{\partial^2u}{\partial\theta^2}\right)+\left(\frac{i}{2kr}-1\right)\frac{\partial u}{\partial r}-\frac{1}{2r}u=0
\end{align}
\mylinenum{and we can write the \RED{BCs} $u\left(r,\pm\theta_{\text{w}}\right)=e^{-ikr}F(r,\pm\theta_{\text{w}})=f(r,\pm\theta_{\text{w}})$. This fits exactly within the realm of \eqref{eq:rw0}, but with complex coefficients. It is shown in {\citep{BudaevBogy2003}} that a suitable choice of the function $q$ is $q^2 (r, \theta) = - ikr^2$. The manifold extending the boundary is chosen as $\partial \mathfrak{G}$, where}
$$\mathfrak{G}=\left\{(r,\theta)\in\mathbb{C}^2,r\in\mathbb{C},-\theta_{\text{w}}<\tmop{Re}\{\theta\}<\theta_{\text{w}}\right\}.$$ 
\mylinenum{In this particular case, the two SDEs to consider become}
\begin{eqnarray}
\mathd\xi_t^1=\xi_t^1\left(ik\xi_t^1+\frac{1}{2}\right)\mathd t+\xi_t^1\mathd W_t^1\ \ \text{and}\ \ \mathd\xi_t^2=\mathd W_t^2
\label{eq:finalSDEs}\end{eqnarray}
\mylinenum{with \RED{ICs} $\xi_0^1 = r$ and $\xi_0^2 = \theta$ and exit time $\tau$ defined such that $\xi_{\tau}^2 = \pm \theta_{\text{w}}$, which now makes sense since the random process $\xi_t^2$ is now real for all times. Using \eqref{eq:solpolar}, we can hence write the solution\footnote{Note that in \citep{BudaevBogy2003}, in their equivalent of the second part of \eqref{eq:cont-sol} (their equation (26)), the argument of the exponential is $-\tfrac{1}{2}\mathd W_\tau^1$. We believe it to be a typographical error.} as}
\begin{eqnarray}
u(r,\theta)=\textbf{E}\left\{f(\xi_{\tau}^1,\xi_{\tau}^2)e^{\frac{ik}{2}\int_0^{\tau}\xi_s^1\mathd s}\right\}\ \ \text{or}\ \ u(r,\theta)=\frac{1}{\sqrt{r}}\textbf{E}\left\{f(\xi_{\tau}^1,\xi_{\tau}^2)\sqrt{\xi_{\tau}^1}e^{-\frac{1}{2}W_{\tau}^1}\right\}, \label{eq:cont-sol}\end{eqnarray}
\mylinenum{where the second part of \eqref{eq:cont-sol} is derived from the first using It{\^o} calculus. Note that for this to be valid, $f$ should be chosen such that it can be analytically continued for $r \in \mathbb{C}$. The second arguments in \eqref{eq:cont-sol} do not pose any problem since by definition $\xi_{\tau}^2 = \pm \theta_{\text{w}}$. The two SDEs \eqref{eq:finalSDEs} are reasonably straightforward to solve numerically (see Figure \ref{fig:random-hankel}, left) using Euler-Maruyama with time step $\Delta t = 0.01$. If the solution we are trying to find is continuous everywhere, the solution \eqref{eq:cont-sol} can be implemented and works well. To illustrate this point we use the same example as in {\citep{BudaevBogy2003}} and apply this method to reproduce the function $H_0^{(1)} (kr)$, which is well known to satisfy the Helmholtz equation and the radiation condition. In order to do so we tailored the \RED{BCs} to be $f (r, \pm \theta_{\text{w}}) = e^{- ikr} H_0^{(1)} (kr)$ and plotted an illustration of the result in Figure \ref{fig:random-hankel}.

\begin{figure}[h!]\centering
	\includegraphics[width=0.45\textwidth]{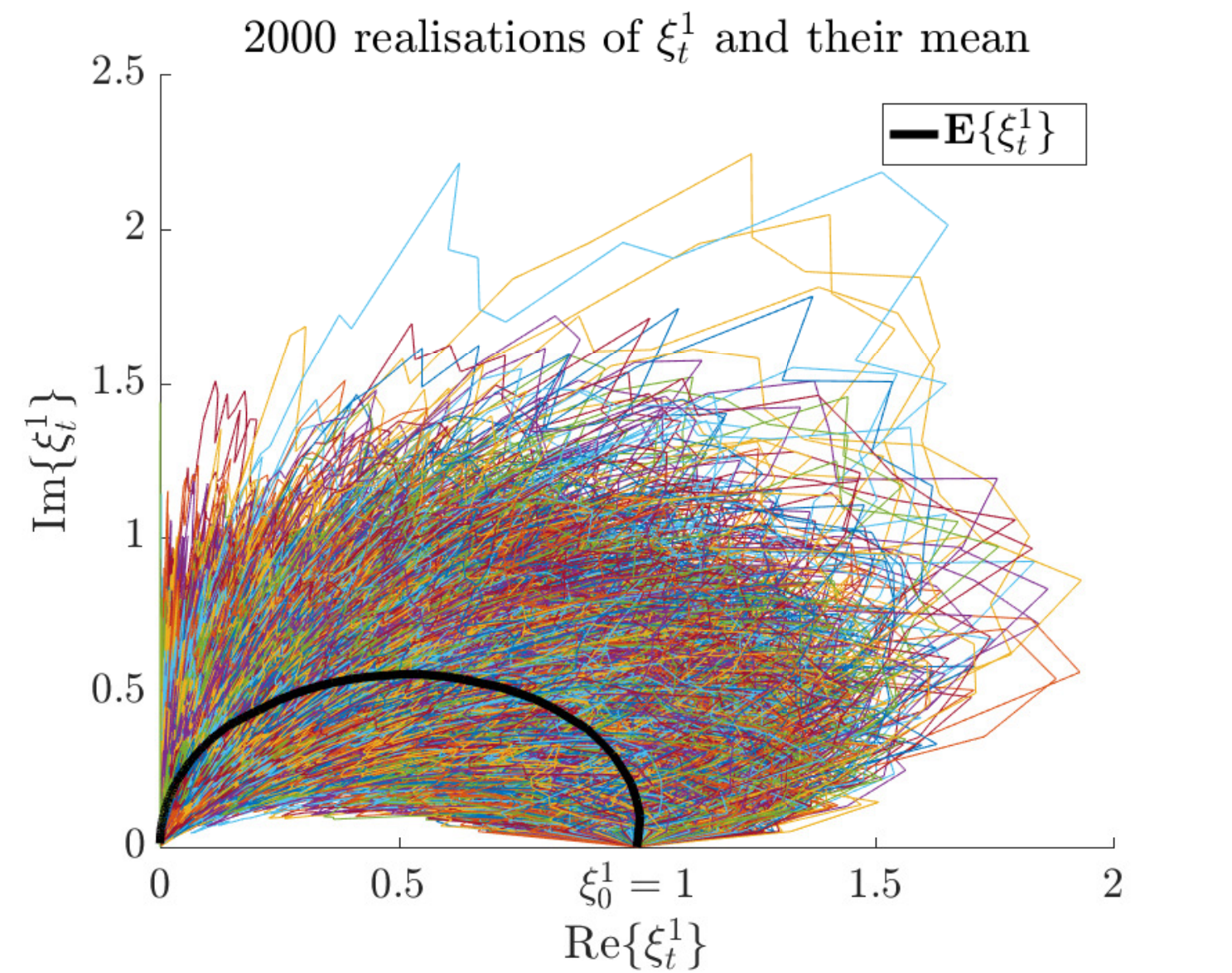}
	\includegraphics[width=0.45\textwidth]{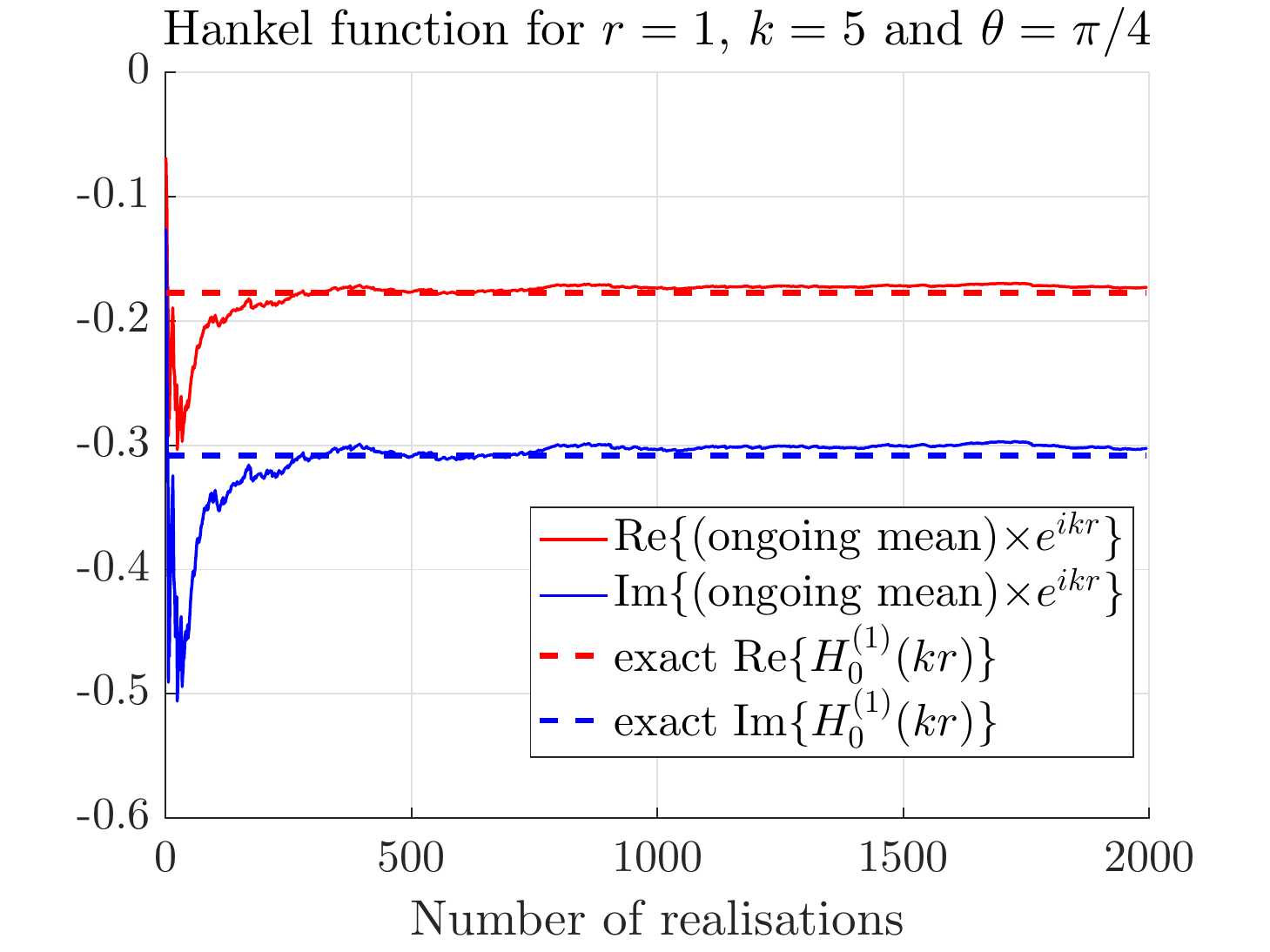}
	\caption{(Left) 2000 realisations of the SDE governing $\xi_t^1$ and their mean, plotted in the complex plane with initial condition $\xi_0^1 = r = 1$. (right) The ongoing mean computed by \eqref{eq:cont-sol} for the Hankel function up to 2000 realisations for $\theta=\pi/4$, $r=1$ and $k=5$.}
	\label{fig:random-hankel}
\end{figure}

If the solution we are seeking has some discontinuities, then the method should be adapted slightly. We are interested in this since what we want to compute is the diffracted field $\Phi_{\text{Diff}} (r, \theta)$ resulting from an incident plane wave with incident angle $\theta_\text{I}$, which satisfies homogeneous Dirichlet \RED{BCs} and the radiation condition. \RED{In what follows, we choose $\theta_\text{I}$ such that both wedge faces are illuminated.}
As shown in Section \ref{SolAnalysis}, the field $\Phi_{\text{Diff}}$ has GO discontinuities\footnote{Note that in {\citep{BudaevBogy2003}}, the convention to choose the index of $\theta_1$ or $\theta_2$ is different, but we have made that choice in order to be consistent with the rest of the review.} at $\theta=\theta_1=2\theta_{\text{w}}-\theta_\text{I}-\pi$ and $\theta=\theta_2=-2\theta_{\text{w}}-\theta_\text{I}+\pi$, and the knowledge of the $\tmop{GO}$ field allows us to derive the following jump conditions across $\theta_{1,2}$:}
\begin{align}\nonumber
[u]_{\theta_1}=1,\ \ [u]_{\theta_2}=-1,\ \ [\partial u/\partial\theta]_{\theta_{1,2}}=0,
\end{align}
\mylinenum{where $u$ is defined such that $\Phi_{\text{Diff}}(r,\theta)=u(r,\theta)e^{ikr}$, and the bracket $[u]_{\theta_{1,2}}=u(r,\theta_{1,2}+0)-u(r,\theta_{1,2}-0)$. Using these jump conditions, it can be shown that \eqref{eq:cont-sol} can be rewritten as}
\begin{align}
u(r,\theta)=\textbf{E}\left\{\sum_{\tau_{\nu}<\tau}(-1)^{m_{\nu}}\delta_{\nu}e^{\frac{ik}{2}\int_0^{\tau_{\nu}}\xi_s^1\mathd s}\right\} \text{ or } 
u(r,\theta)=\frac{1}{\sqrt{r}}\textbf{E}\left\{\sum_{\tau_{\nu}<\tau}(-1)^{m_{\nu}}\delta_{\nu}\sqrt{\xi_{\tau_{\nu}}^1}e^{-\frac{1}{2}W_{\tau_{\nu}}^1}\right\},
\label{eq:disc-sol}\end{align}
\mylinenum{where for each simulation, the $\tau_{\nu}$ represent the times of crossings between $\xi_t^2$ and the discontinuous lines $\theta_{1, 2}$. If $\theta_1$ (resp. $\theta_2$) is crossed, then $m_{\nu} = 2$ (resp. 1). If the crossing is from above (resp. below), then $\delta_{\nu} = 1$ (resp. $- 1$). As illustrated in Figure \ref{fig:plane-wave-random} (left), many such crossings can occur before the exit time $\tau$ is reached. The method has been implemented for a wedge characterised by $\theta_{\text{w}} = 7 \pi / 8$, and the results, obtained for 2000 realisations (simulated by Euler-Maruyama with time step $\Delta t = 0.01$), are shown at an observation point $r = 1$, $\theta = \pi / 4$ for $k = 5$. Note that if the method was described in {\citep{BudaevBogy2003}}, it was only implemented for a half-plane, and not for a wedge. Though, as predicted in {\citep{BudaevBogy2003}} the error is of the order of 0.01 and the method works well\footnote{\RED{ Note that in {\citep{BudaevBogy2003}}, there is a factor $\frac{1}{2}$ in front of \textbf{E} in the formulae \eqref{eq:disc-sol}. This was a typographical error.} 
}.

\begin{figure}[h!]
\centering
  \includegraphics[width=0.45\textwidth]{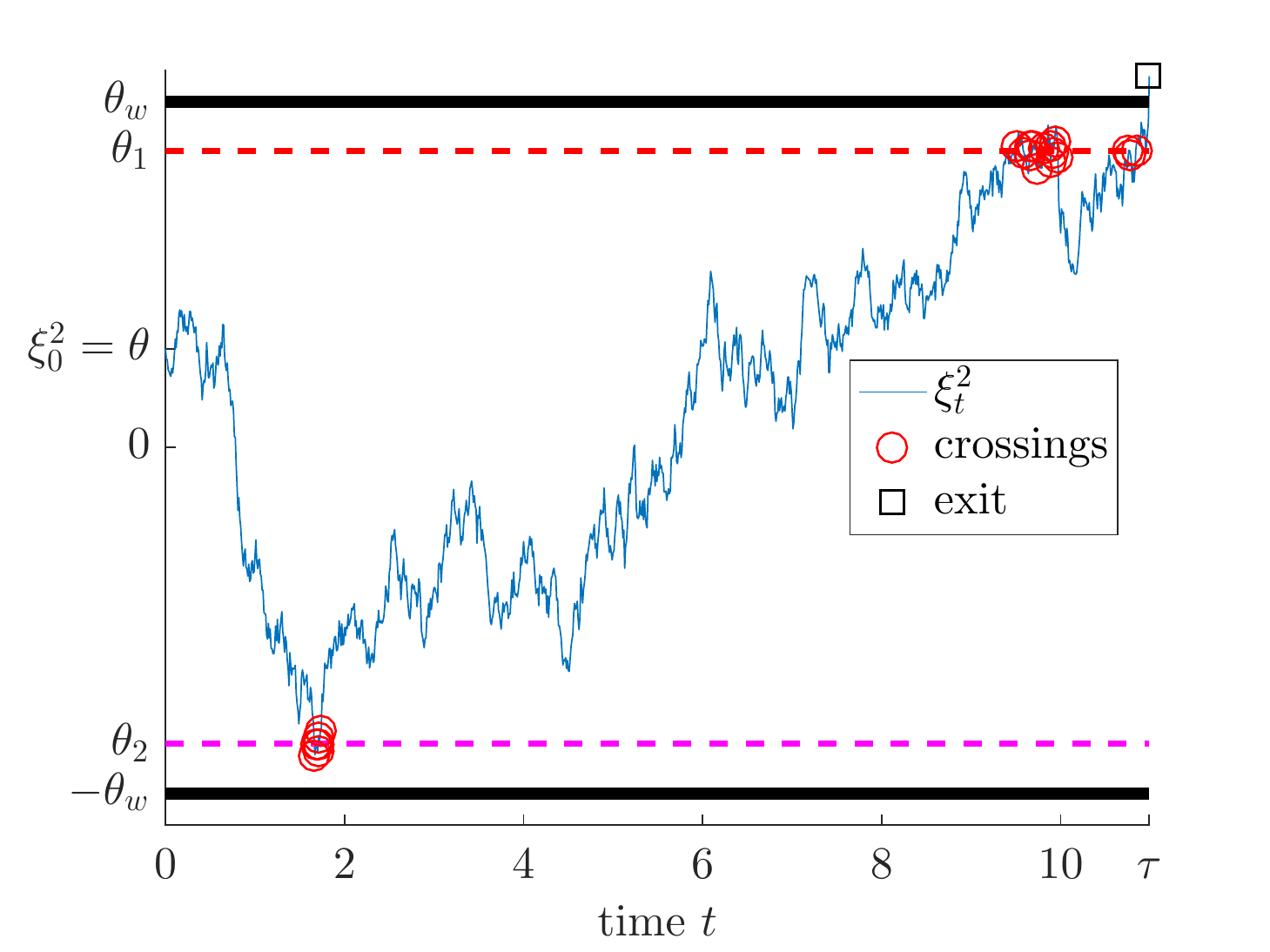}
  \includegraphics[width=0.45\textwidth]{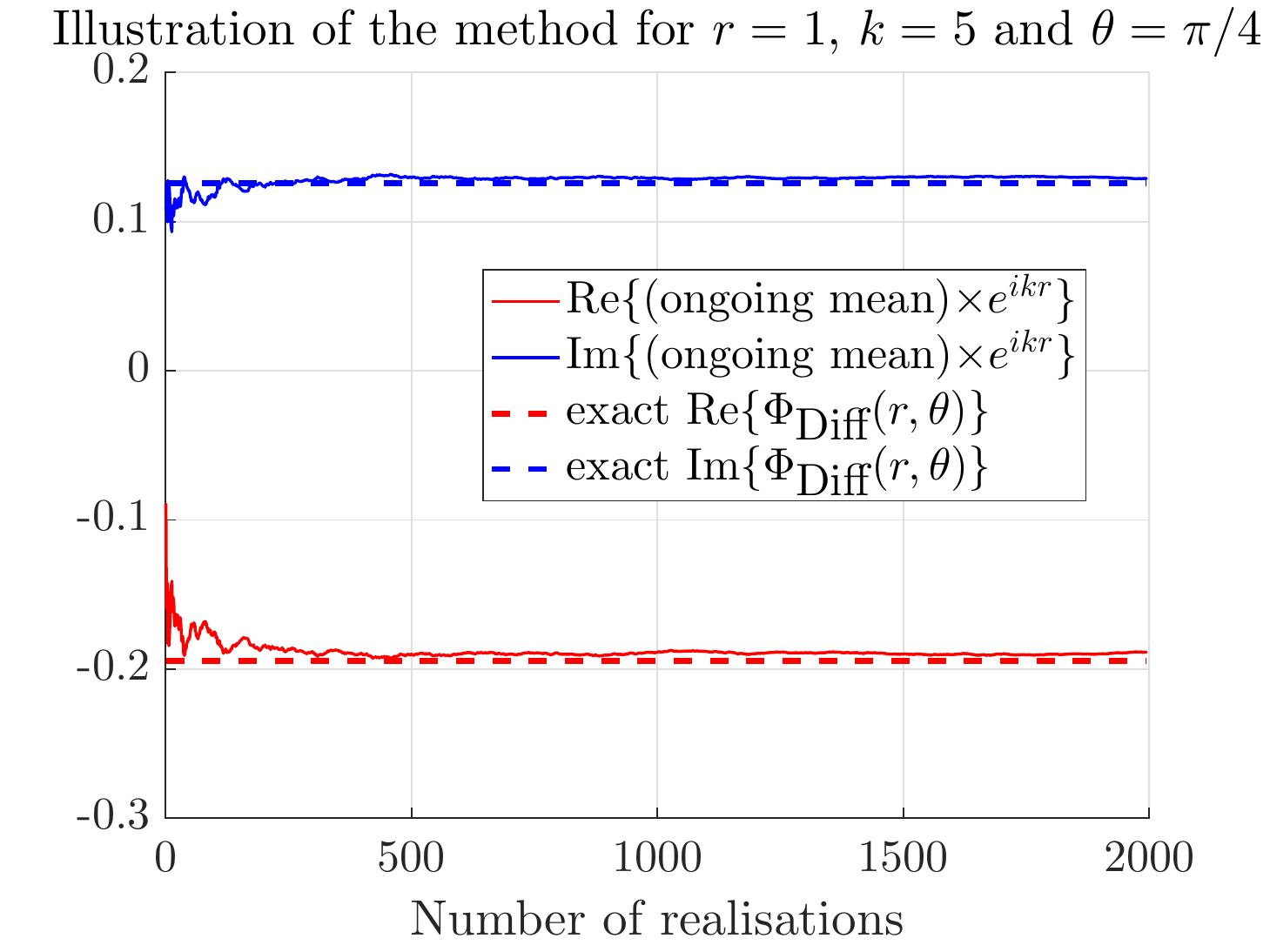}
  \caption{(Left) Illustration of one realisation of $\xi_t^2$, its crossings
  with $\theta_{1, 2}$ and its exit time $\tau$. (Right) The ongoing mean up
  to 2000 realisations computed by \eqref{eq:disc-sol} for a plane wave
  incident at an angle $\theta_\text{I} = 0$ on a wedge characterised by
  $\theta_{\text{w}}=7\pi/8$ for $\theta = \pi / 4$, $r = 1$ and $k
  = 5$.}
\label{fig:plane-wave-random}
\end{figure}


\RED{
\paragraph{Critical analysis} This method has the advantage of being very adaptable to all sorts of geometries since it is based on the Feynman-Kac theorem that is a very general result (both in terms of geometry and in terms of equation). This adaptability
is confirmed by the fact that it has been used in the context of cones \citep{BudaevBogy2003}, quarter-plane \citep{BudaevBogy2005} and other geometries.

It has to be said however, that for the Helmholtz equation, the PDE and SDE coefficients become complex. This renders the determination of the end time rather more complicated than the real coefficient case. It necessitates to find a convenient complex coefficient to multiply our equations by, and also to find a way of somehow extending the real geometries in a higher dimension complex space.

This method can also become very computational very quickly. Indeed, if one would like for example to recreate a heat map similar to those presented in Figure \ref{fig:SA-UTD}, one would need about 2000 simulations of the SDE system per point, which for a good resolution may lead to a very long computational time.

Another comment that can be made about this method, is that it stands out from all the other methods presented here in terms of the type of mathematics used. This can be considered as an advantage for researchers open to exploring many areas of mathematics, though, this also means that for the usual specialists in diffraction theory, this may result in a steep learning curve.
}



\subsection{The method of functionally-invariant solutions}\label{SSM} 
The \RED{third and final} alternative method to be reviewed is also known as the Sobolev-Smirnov method. Some recent publications using this method include \citep{KMM2015,Babich2015}. The former studies wedge diffraction with a number of different combinations of Dirichlet and Neumann \RED{BCs}, while the latter studies the impedance wedge problem.

The idea behind this method is to identify the time-harmonic problem with an elementary time-dependent problem where the incident plane wave is replaced with a Heaviside step function such that no diffraction occurs before the time $t=0$. This means that the solution to this elementary problem (call it $u(r,\theta,t)$) satisfies the following conditions,
\begin{itemize}
\item{The governing equation is the linear wave equation $\nabla^2u-\frac{1}{c^2}\frac{\partial^2u}{\partial t^2}=0$.}
\item{Dirichlet or Neumann \RED{BCs} at $\theta=\pm\theta_\text{w}$.}
\item{$u$ can be linearly decomposed into incident and scattered parts, $u=u_\text{I}+u_\text{S}$, where $u_\text{I}(r,\theta,t)=\mathcal{H}\left(t+\frac{r}{c}\cos(\theta-\theta_\text{I})\right)$.}
\end{itemize}

For simplicity, we shall restrict values of the incident angle and the wedge angle such that, $\pi-\theta_\text{w}<\theta_\text{I}<\theta_\text{w}-\frac{\pi}{2}$. This restriction means that for $t<0$, the incident wave does not reach the wedge until $t=0$ when it first touches the wedge at its corner. For $t>0$, the incident wave has passed the wedge corner and reflected and diffracted waves have appeared. Figure \ref{SSM-geometry} describes this configuration and gives known values of $u$ outside the diffraction circle which are found by Geometrical Optics.

\begin{figure}[ht]\centering
\includegraphics[width=0.7\textwidth]{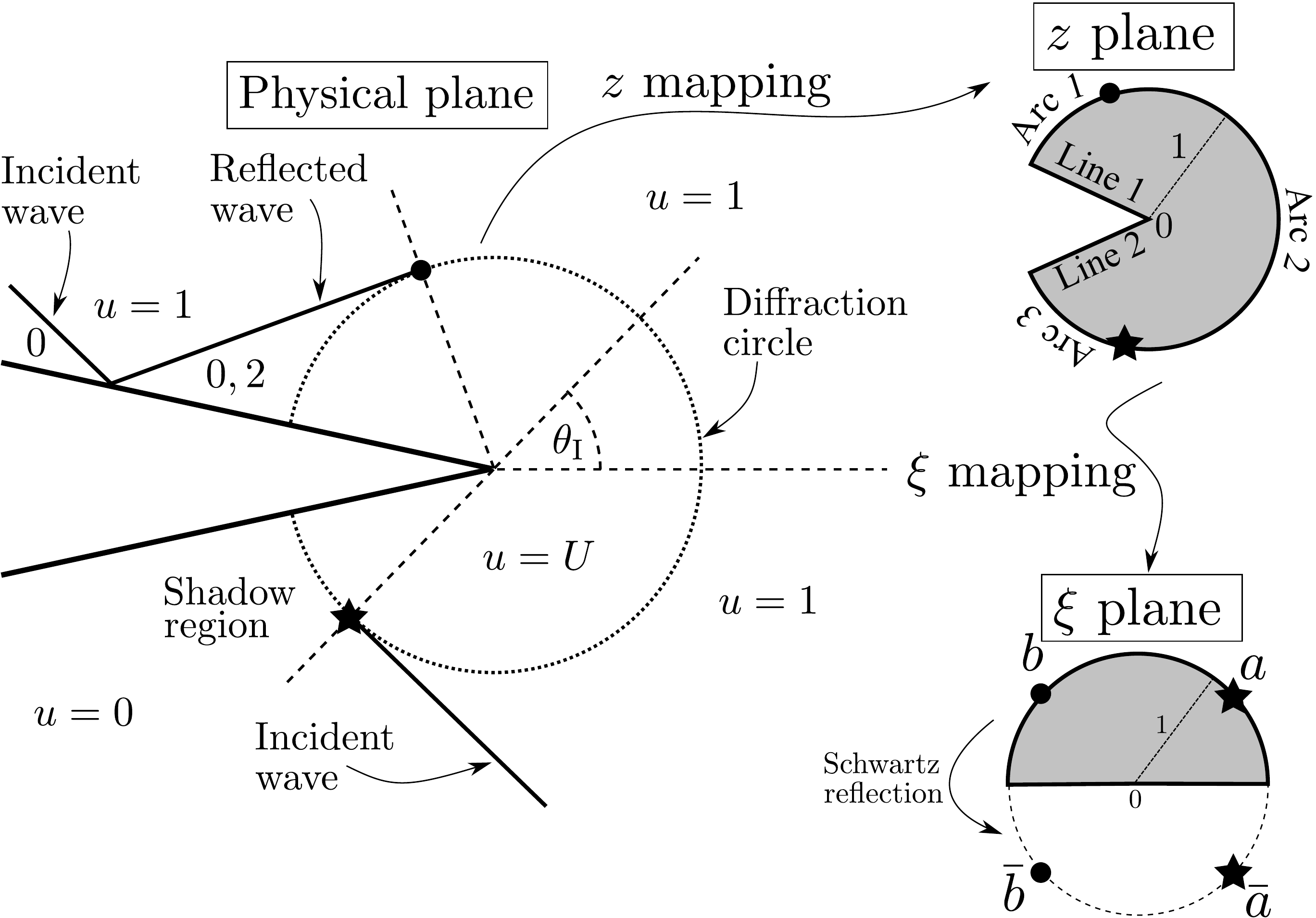}
\caption{Physical diagram and images of the diffraction disc in the $z$ and $\xi$ planes.}
\label{SSM-geometry}
\end{figure}

The radius of the diffraction circle is $ct$ and we call the unknown solution inside, i.e. within the diffraction disc, $U(r,\theta,t)$. We need to look at a particular class of solutions to the wave equation and express $U$ in terms of a complex variable. Noting that a real solution is required, we write}
\begin{align}\label{SSM-invariance}
U(r,\theta,t)=\RE{V(z)},\quad\text{where,}\quad z=\frac{ct}{r}\left(1-\sqrt{1-\left(\frac{r}{ct}\right)^2}\right)e^{i\theta},
\end{align}
\mylinenum{ensuring that the wave equation is automatically satisfied. Sections 52-53 in \citep{Smirnov1964} gives a detailed explanation as to why this is the case. Note also that within the diffraction disc (i.e. $0\leq r \leq ct$), the pre-exponential factor in \eqref{SSM-invariance} is real and positive, varying from $0$ to $1$. As a result, $z(r,\theta,t)$ maps the diffraction disc onto the unit disc $|z|\leq1$. Therefore, we need to find a function $V(z)$ that is analytic inside the unit disc and has the following boundary values,}
\begin{align}
\nonumber\text{Both: } & \RE{V(z)}=0 \text{ on Arc 3,}\ \ \RE{V(z)}=1 \text{ on Arc 2,} \\
\nonumber\text{Dirichlet: } & \RE{V(z)}=0 \text{ on Arc 1,}\ \ \RE{V(z)}=0 \text{ on Lines 1 \& 2,} \\
\nonumber\text{Neumann: } & \RE{V(z)}=2 \text{ on Arc 1,}\ \ \RE{izV'(z)}=0 \text{ on Lines 1 \& 2,}
\end{align}
\mylinenum{where the arc and line numbers are given in Figure \ref{SSM-geometry}.

We will now use a conformal mapping to transform this problem into a Riemann-Hilbert problem. In order to do so, we define $\xi=e^{\frac{i\pi}{2}}z^\delta$, where as before $\delta=\tfrac{\pi}{2\theta_{\text{w}}}$, which transforms the indented unit disc of the $z$ plane onto the unit upper-half semi-disc of the $\xi$ plane, as illustrated in Figure \ref{SSM-geometry}. The branch of the root is defined such that the cut is on the negative real axis and $\xi(z=1)=e^{\frac{i\pi}{2}}$. Let $\widetilde{V}(\xi)=V(z(\xi))$, then we analytically continue $\widetilde{V}$ into the unit lower-half semi-disc by \RED{Schwarz} reflection principle (see Figure \ref{SSM-geometry}) using anti-symmetry (Dirichlet case) or symmetry (Neumann case). Let $a=e^{i(\frac{\pi}{2}-\delta(\pi-\theta_\text{I}))}$ and $b=e^{i(\frac{3\pi}{2}-\delta(\pi+\theta_\text{I}))}$ then $\widetilde{V}(\xi)$ has the following boundary values for the Dirichlet and Neumann cases,}
\begin{eqnarray}
   & \tmop{Dir:} \left\{ \begin{array}{ll}
    \tmop{Re} \{ \widetilde{V} (\xi) \} = - 1 \text{ on arc } (\bar{b},
    \bar{a}), & \tmop{Re} \{ \widetilde{V} (\xi) \} = 1 \text{ on arc } (a,
    b),\\
    \tmop{Re} \{ \widetilde{V} (\xi) \} = 0 \text{ on real line } (- 1, 1), &
    \tmop{Re} \{ \widetilde{V} (\xi) \} = 0 \text{ on arcs } (\bar{a}, a)
    \text{ and } (b, \bar{b}),
  \end{array} \right. & \label{SSM-Dproblem}\\
   & \tmop{Neu:} \left\{ \begin{array}{ll}
    \tmop{Re} \{ \widetilde{V} (\xi) \} = 0 \text{ on arc } (\bar{a}, a), &
    \tmop{Re} \{ \widetilde{V} (\xi) \} = 2 \text{ on arc } (b, \bar{b}),\\
    \tmop{Re} \{ i \xi \widetilde{V}' (\xi) \} = 0 \text{ on real line } (- 1,
    1), & \tmop{Re} \{ \widetilde{V} (\xi) \} = 1 \text{ on arcs } (a, b)
    \text{ and } (\bar{b}, \bar{a}).
  \end{array} \right. \label{SSM-Nproblem}& 
\end{eqnarray}

\mylinenum{The method to solve these two Riemann-Hilbert problems is detailed in section 54 in \citep{Smirnov1964}. The respective solutions to \eqref{SSM-Dproblem} and \eqref{SSM-Nproblem} are,}
\begin{align}
\label{SSM-hilbert-Dsol}\widetilde{V}(\xi)&=\frac{1}{\pi i}\ln\left(\frac{\bar{b}-\xi}{a -\xi}\right)-\frac{1}{\pi i}\ln\left( \frac{\bar{a}-\xi}{b -\xi}\right),\\
\label{SSM-hilbert-Nsol}\widetilde{V}(\xi)&=\frac{1}{\pi i}\ln\left(\frac{\bar{b}-\xi}{a -\xi}\right)+\frac{1}{\pi i}\ln\left(\frac{\bar{a}-\xi}{b-\xi}\right)-2\delta,\end{align}
\mylinenum{where the \RED{used logarithm $\ln(Z)$} has a branch point at $Z=0$ with a branch cut along the positive real axis. Using this solution, it is easy to recover the physical solution $U(r,\theta,t)$ inside the diffraction disc, and hence the whole solution $u(r,\theta,t)$. We will now see that using a simple Fourier transform, we can recover the sought-after time-harmonic problem from this solution $u(r,\theta,t)$. Consider the evaluation to the following integral, assuming that \RED{$\omega$ has a small positive imaginary increment so that} $e^{i\omega t}|_{t=\infty}$ is zero,}
\begin{align}\label{SSM-uI->PhiI}
-\int_{-\infty}^\infty u_\text{I}(r,\theta,t)\Der{}{t}\left(e^{i\omega t}\right)\text{d}t=\left[e^{i\omega t}\right]_\infty^{-\frac{r}{c}\cos(\theta-\theta_\text{I})}=\Phi_\text{I}.
\end{align}
\mylinenum{With this in mind, we can determine the total field $\Phi$ from $u$ by using a similar integral,}
\begin{align}\label{SSM-u->Phi}
\Phi(r,\theta)=-\int_{-\infty}^\infty u(r,\theta,t)\Der{}{t}\left(e^{i\omega t}\right)\text{d}t,
\end{align}
\mylinenum{and thus, we have found the solution to the time-harmonic problem.

\RED{
\paragraph{Critical analysis} All of the methods presented so far were tailored to the time-harmonic problem, this means that if one is interested in a time-dependent problem, using these methods would involve taking the inverse Fourier transform of
our time-harmonic solutions, which can be expensive computationally. This present method however is tailored to the time dependent problem, which is great if one is interested in the tracking of wave fronts in time for example. It means however that if one is interested in the time-harmonic problem, one would have to take the Fourier transform in time of the solution, as per (\ref{SSM-u->Phi}), which can prove quite expensive numerically. This method, though in essence designed for the wedge geometry, has been shown to be adaptable to various BCs. One can refer to \citep{Babich2015} for example for the case of impedance BCs.
}


\section{Final plots and conclusions} 
In this review article, we have discussed six different methods that have been applied to the problem of diffraction by wedges with perfect Dirichlet or Neumann boundary conditions. The three main methods discussed were the Sommerfeld-Malyuzhinets technique, the Wiener-Hopf technique and the Kontorovich-Lebedev transform technique. The three alternative methods reviewed were the \RED{embedding formula, the random walk method and the method of functionally-invariant solutions (Sobolev-Smirnov).} We also looked at two approximation methods, the Geometrical Theory of Diffraction and the Uniform Geometrical Theory of Diffraction.

This list is by no means exhaustive and we should also mention Budaev's method for elastic wedge scattering \citep{BudaevBogy1998}, the Physical Theory of Diffraction \citep{Ufimtsev2014} and an interesting method called the Wiener--Hopf--Hankel formulation \citep{Teixeira1991,Castro2010}. We note that \citep{Israilov2013} could also be applied to the wedge geometry.

We evaluated numerically the exact solution and the associated approximations (series, GTD, UTD) for several configurations and studied their relative performances. We found that the best way to evaluate the exact solution was to consider the integral defined on the steepest descent contour. As regard to the approximations, the truncated series solutions \RED{performs very well with low wavenumbers}, and we found that while the \RED{UTD approximation} takes longer to compute, it is a better approximation compared with the GTD because it is uniformly valid and more accurate at lower values of $kr$. It has however two main disadvantages, the inaccuracy at the wedge boundary, and also the fact that \eqref{SA-UTD1st-sol} and \eqref{SA-UTD2nd-sol} are not continuous across the Geometrical Optics limit $\theta_\text{I}=\pi-\theta_{\text{w}}$.

\RED{As emphasised in the critical analysis of each section, the} use of the six techniques in this review is not limited to the perfect wedge problem. Examples of extensions include for example impedance wedges \citep{Malyuzhinets1958-3,Babich2015}, penetrable wedges \citep{Rawlins1999,Lyalinov1999,DanieleLombardi2011} and quarter-plane diffraction \citep{shanin1,AssierPeake2012,BudaevBogy2005,Lyalinov2013}.

To conclude this review, using the UTD approximation, we produce some density plots of the real part of the diffracted field $\Phi_{\text{Diff}}$ and the total field $\Phi$ in Figure \ref{fig:SA-UTD} for a wedge defined by $2\bar{\theta}_\text{w}=\pi/4$ and  two incident waves, $\theta_\text{I}=0$ and $\theta_\text{I}=\pi/2$ and a wavenumber $k=2$. As expected, we see clear discontinuities in the diffracted wave, which are due to GO discontinuities. For the total field, as expected, we see the GO behaviour in the relevant regions, the boundary conditions, and a decaying diffracted field.

\begin{figure}[h!]\centering
\includegraphics[width=0.9\textwidth]{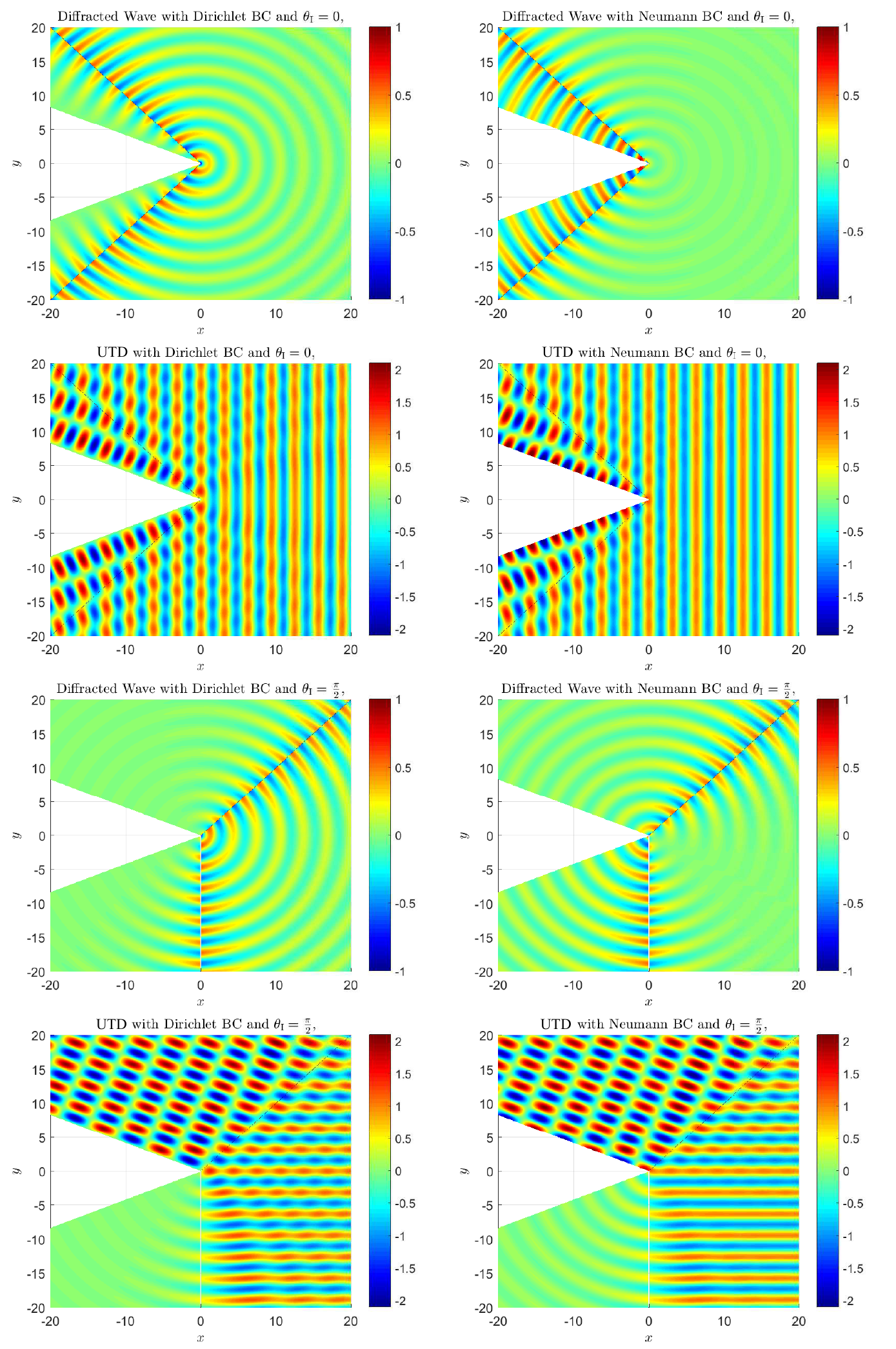}
\caption{Density plots of the UTD approximation of $\RE{\Phi_{\text{Diff}}}$ and $\RE{\Phi}$ for Dirichlet and Neumann \RED{BCs}, for $k=2$ and for a wedge characterised by $\theta_\text{w}=7\pi/8$ and two incident angles $\theta_\text{I}=0$ and $\theta_\text{I}=\pi/2$.}
\label{fig:SA-UTD}
\end{figure}

\section*{Acknowledgements}
Nethercote acknowledges the financial support from EPSRC (DTA studentship). Assier acknowledges that part of this work was supported by EPSRC (EP/N013719/1). Abrahams acknowledges the support by UKRI under grant no EP/K032208/1. \RED{The authors thank the anonymous reviewers for their help in improving this work. Assier would also like to thank A.V. Shanin for fruitful discussions, especially on the topic of \ref{app:app1}, and B. Budaev for being so responsive when queried about detailed aspects of the random walk method.}}


\clearpage
{
 \def\section*#1{}
	\bibliographystyle{elsarticle-harv_raph} 
	\bibliography{IWRlibrary}

\begin{thebibliography}{103}
\expandafter\ifx\csname natexlab\endcsname\relax\def\natexlab#1{#1}\fi
\expandafter\ifx\csname url\endcsname\relax
  \def\url#1{\texttt{#1}}\fi
\expandafter\ifx\csname urlprefix\endcsname\relax\def\urlprefix{URL }\fi

\bibitem[{Abrahams(1986)}]{Abrahams1986}
Abrahams, I.~D., 1986. {Diffraction by a semi-infinite membrane in the presence
  of a vertical barrier}. J. Sound Vib. 111~(2), 191--207.

\bibitem[{Abrahams(1987)}]{Abrahams1987}
Abrahams, I.~D., 1987. {On the Sound Field Generated by Membrane Surface Waves
  on a Wedge-Shaped Boundary}. Proc. R. Soc. A Math. Phys. Eng. Sci. 411,
  239--250.

\bibitem[{Abramowitz and Stegun(1964)}]{Handbook}
Abramowitz, M., Stegun, I.~A., 1964. {Handbook of Mathematical Functions}.
  Courier Corporation.

\bibitem[{Assier and Peake(2012{\natexlab{a}})}]{AssierPeake2012}
Assier, R.~C., Peake, N., 2012{\natexlab{a}}. {On the diffraction of acoustic
  waves by a quarter-plane}. Wave Motion 49~(1), 64--82.

\bibitem[{Assier and Peake(2012{\natexlab{b}})}]{Assier2012b}
Assier, R.~C., Peake, N., 2012{\natexlab{b}}. {Precise description of the
  different far fields encountered in the problem of diffraction of acoustic
  waves by a quarter-plane}. IMA J. Appl. Math. 77~(5), 605--625.

\bibitem[{Babich(2015)}]{Babich2015}
Babich, V.~M., 2015. {Solution of the Diffraction Problem of a Plane Wave by an
  Impedance Wedge (Non-Stationary Case, Smirnov-Sobolev Method).} Russ. J.
  Math. Phys. 22~(2), 145--152.

\bibitem[{Babich et~al.(2007)Babich, Lyalinov, and Grikurov}]{SMtechnique2007}
Babich, V.~M., Lyalinov, M.~A., Grikurov, V.~E., 2007. {Diffraction Theory: The
  Sommerfeld-Malyuzhinets Technique}. Alpha Science.

\bibitem[{Bender and Orszag(1999)}]{BenderOrszag1999}
Bender, C.~M., Orszag, S.~A., 1999. {Advanced Mathematical Methods for
  Scientists and Engineers I: Asymptotic Methods and Perturbation Theory}.
  Springer-Verlag, New York.

\bibitem[{Biggs(2001)}]{Biggs2001}
Biggs, N. R.~T., 2001. {Wave Diffraction Through a Perforated Barrier of
  Non-Zero Thickness}. Q. J. Mech. Appl. Math. 54~(4), 523--547.

\bibitem[{Biggs(2002)}]{Biggs2002}
Biggs, N. R.~T., 2002. {Wave Scattering by a Perforated Duct}. Q. J. Mech.
  Appl. Math. 55~(2), 249--272.

\bibitem[{Biggs(2006)}]{Biggs2006}
Biggs, N. R.~T., 2006. {A new family of embedding formulae for diffraction by
  wedges and polygons}. Wave Motion 43~(7), 517--528.

\bibitem[{Bleistein and Handelsman(2010)}]{Bleistein}
Bleistein, N., Handelsman, R.~A., 2010. Asymptotic Expansions of Integrals
  (Dover Books on Mathematics). Dover Publications.

\bibitem[{Borovikov and Kinber(1994)}]{BorovikovKinber1994}
Borovikov, V.~A., Kinber, B.~Y., 1994. {Geometrical Theory of Diffraction}.
  Institution of Electrical Engineers, London.

\bibitem[{Bowman et~al.(1987)Bowman, Senior, and Uslenghi}]{Bowman1987}
Bowman, J.~J., Senior, T. B.~A., Uslenghi, P. L.~E., 1987. {Electromagnetic and
  Acoustic Scattering by Simple Shapes}. Hemisphere Publishing Corporation.

\bibitem[{Bromwich(1915)}]{Bromwich1915}
Bromwich, T.~J., 1915. {Diffraction of Waves by a Wedge}. Proc. London Math.
  Soc. 14~(1), 450--463.

\bibitem[{Budaev(1995)}]{Budaev1995}
Budaev, B.~V., 1995. {Diffraction by Wedges}, 1st Edition. CRC Press.

\bibitem[{Budaev and Bogy(1995)}]{BudaevBogy1995}
Budaev, B.~V., Bogy, D.~B., 1995. {Rayleigh wave scattering by a wedge}. Wave
  Motion 22~(1), 239--257.

\bibitem[{Budaev and Bogy(1996)}]{BudaevBogy1996}
Budaev, B.~V., Bogy, D.~B., 1996. {Rayleigh wave scattering by a wedge II}.
  Wave Motion 24~(3), 307--314.

\bibitem[{Budaev and Bogy(1998)}]{BudaevBogy1998}
Budaev, B.~V., Bogy, D.~B., 1998. {Rayleigh wave scattering by two adhering
  elastic wedges}. Proc. R. Soc. A Math. Phys. Eng. Sci. 454, 2949--2996.

\bibitem[{Budaev and Bogy(2001)}]{BudaevBogy2001}
Budaev, B.~V., Bogy, D.~B., 2001. {Probabilistic solutions of the Helmholtz
  equation}. J. Acoust. Soc. Am. 109~(5), 2260--2262.

\bibitem[{Budaev and Bogy(2002{\natexlab{a}})}]{BudaevBogy2002a}
Budaev, B.~V., Bogy, D.~B., 2002{\natexlab{a}}. {Application of random walk
  methods to wave propagation}. Q. J. Mech. Appl. Math. 55~(2), 209--226.

\bibitem[{Budaev and Bogy(2002{\natexlab{b}})}]{BudaevBogy2002b}
Budaev, B.~V., Bogy, D.~B., 2002{\natexlab{b}}. {Random walk methods and wave
  diffraction}. Int. J. Solids Struct. 39~(21-22), 5547--5570.

\bibitem[{Budaev and Bogy(2003)}]{BudaevBogy2003}
Budaev, B.~V., Bogy, D.~B., 2003. {Random walk approach to wave propagation in
  wedges and cones}. J. Acoust. Soc. Am. 114~(4), 1733--1741.

\bibitem[{Budaev and Bogy(2005)}]{BudaevBogy2005}
Budaev, B.~V., Bogy, D.~B., 2005. {Diffraction of a plane wave by a sector with
  Dirichlet or Neumann boundary conditions}. IEEE Trans. Antennas Propag.
  53~(2), 711--718.

\bibitem[{Busemann(1947)}]{Busemann1947}
Busemann, A., 1947. {Infinitesimal conical supersonic flow}. Tech. rep.,
  National Advisory Committee for Aeronautics, Washington D.C.

\bibitem[{Castro and Kapanadze(2010)}]{Castro2010}
Castro, L.~P., Kapanadze, D., 2010. {Exterior Wedge Diffraction Problems with
  Dirichlet, Neumann and Impedance Boundary Conditions}. Acta Appl. Math. 110,
  289--311.

\bibitem[{Copson(1946)}]{Copson1946}
Copson, E.~T., 1946. {On an integral equation arising in the theory of
  diffraction}. Q. J. Math. 17~(1), 19--34.

\bibitem[{Craster and Shanin(2005)}]{ShaninCraster2005}
Craster, R., Shanin, A.~V., 2005. {Embedding formulae for diffraction by
  rational wedge and angular geometries}. Proc. R. Soc. A Math. Phys. Eng. Sci.
  461, 2227--2242.

\bibitem[{Croisille and Lebeau(1999)}]{CroisilleLebeau1999}
Croisille, J.-P., Lebeau, G., 1999. {Diffraction by an Immersed Elastic Wedge}.
  Springer-Verlag.

\bibitem[{Daniele(2000)}]{Daniele2000}
Daniele, V.~G., 2000. {Generalized Wiener-Hopf Technique for Wedge Shaped
  Regions of Arbitrary Angles}. In: VIII-th Int. Conf. Math. Methods
  Electromagn. Theory. Kharkov Ukraine, pp. 432--434.

\bibitem[{Daniele(2001)}]{Daniele2001}
Daniele, V.~G., 2001. {New analytical Methods for wedge problems}. In: Proc.
  2001 Int. Conf. Electromagn. Adv. Appl. pp. 385--393.

\bibitem[{Daniele(2003{\natexlab{a}})}]{Daniele2003-2}
Daniele, V.~G., 2003{\natexlab{a}}. {Rotating waves in the Laplace domain for
  angular regions}. Electromagnetics 23~(3), 223--236.

\bibitem[{Daniele(2003{\natexlab{b}})}]{Daniele2003}
Daniele, V.~G., 2003{\natexlab{b}}. {The Wiener-Hopf Technique for Impenetrable
  Wedges having Arbitrary Aperture Angle}. SIAM J. Appl. Math. 63~(4),
  1442--1460.

\bibitem[{Daniele and Lombardi(2006)}]{DanieleLombardi2006}
Daniele, V.~G., Lombardi, G., 2006. {Wiener-Hopf solution for impenetrable
  wedges at skew incidence}. IEEE Trans. Antennas Propag. 54~(9), 2472--2485.

\bibitem[{Daniele and Lombardi(2011)}]{DanieleLombardi2011}
Daniele, V.~G., Lombardi, G., 2011. {The Wiener-Hopf solution of the isotropic
  penetrable wedge problem: Diffraction and total field}. IEEE Trans. Antennas
  Propag. 59~(10), 3797--3818.

\bibitem[{Daniele and Zich(2014)}]{DanieleZich2014}
Daniele, V.~G., Zich, R.~S., 2014. {The Wiener-Hopf Method in
  Electromagnetics}. Scitech, Edison, New Jersey.

\bibitem[{Felsen and Marcuvitz(1994)}]{Felsen1994}
Felsen, L.~B., Marcuvitz, N., 1994. {Radiation and scattering of waves}, 2nd
  Edition. John Wiley {\&} Sons, Hoboken, N.J.

\bibitem[{Feynman(1948)}]{Feynman1948}
Feynman, R.~P., 1948. {Space-time approach to non-relativistic quantum
  mechanics}. Rev. Mod. Phys. 20, 367--387.

\bibitem[{Filippov(1964)}]{Filippov1964}
Filippov, A.~F., 1964. {Diffraction of an arbitrary acoustic wave by a wedge}.
  J. Appl. Math. Mech. 28~(2), 372--388.

\bibitem[{Fock(1965)}]{Fock1965}
Fock, V.~A., 1965. {Electromagnetic Diffraction and Propagation Problems}.
  Pergamon Press.

\bibitem[{Freidlin(1985)}]{Freidlin1985}
Freidlin, M., 1985. Functional Integration and Partial Differential Equations.
  (AM-109), Volume 109 (Annals of Mathematics Studies). Princeton University
  Press.

\bibitem[{Gautesen(1983)}]{Gautesen1983}
Gautesen, A.~K., 1983. {On the Green's function for acoustical diffraction by a
  strip}. The Journal of the Acoustical Society of America 74~(2), 600.

\bibitem[{Gradshteyn and Ryzhik(2014)}]{TablesISP8th}
Gradshteyn, I.~S., Ryzhik, I.~M., 2014. {Table of Integrals, Series, and
  Products}, 8th Edition. Academic Press.

\bibitem[{Hacivelioglu et~al.(2011)Hacivelioglu, Sevgi, and Ufimtsev}]{HSU2011}
Hacivelioglu, F., Sevgi, L., Ufimtsev, P.~Y., 2011. {Electromagnetic Wave
  Scattering from a Wedge with Perfectly Reflecting Boundaries : Analysis of
  Asymptotic Techniques}. IEEE Antennas Propag. Mag. 53~(3), 232--253.

\bibitem[{Israilov(2013)}]{Israilov2013}
Israilov, M.~S., 2013. {Diffraction of acoustic and elastic waves on a
  half-plane for boundary conditions of various types}. Mechanics of Solids
  48~(3), 337--347.

\bibitem[{James(1986)}]{GLJames1986}
James, G.~L., 1986. {Geometrical Theory of Diffraction for Electromagnetic
  Waves}, 3rd Edition. Institution of Engineering and Technology, London.

\bibitem[{Jones(1952)}]{DSJones1952}
Jones, D.~S., 1952. {A simplifying technique in the solution of a class of
  diffraction problems}. Q. J. Math. 3~(1), 189--196.

\bibitem[{Jones(1964)}]{DSJones1964}
Jones, D.~S., 1964. {The Theory of Electromagnetism}. Pergamon Press, London.

\bibitem[{Jones(1980)}]{DSJones1980}
Jones, D.~S., 1980. {The Kontorovich-Lebedev Transform}. IMA J. Appl. Math.
  26~(2), 133--141.

\bibitem[{Jones(1986)}]{DSJones1986}
Jones, D.~S., 1986. {Acoustic and Electromagnetic Waves}. Oxford University
  Press, Oxford.

\bibitem[{Kac(1949)}]{Kac1949}
Kac, M., 1949. {On the distribution of certain Wiener functionals}. Trans. Am.
  Math. Soc. 65~(1), 1--13.

\bibitem[{Keller(1962)}]{Keller1962}
Keller, J.~B., 1962. {Geometrical Theory of Diffraction}. J. Opt. Soc. Am.
  52~(2), 116--130.

\bibitem[{Keller and Blank(1951)}]{KellerBlank1951}
Keller, J.~B., Blank, A., 1951. {Diffraction and reflection of pulses by wedges
  and corners}. Commun. Pure Appl. Math. 4~(1), 75--94.

\bibitem[{Kisil(2015)}]{Kisil2015}
Kisil, A.~V., 2015. {The relationship between a strip Wiener-Hopf problem and a
  line Riemann-Hilbert problem}. IMA J. Appl. Math. 80, 1569--1581.

\bibitem[{Knopoff(1969)}]{Knopoff1969}
Knopoff, L., 1969. {Elastic Wave Propagation in a Wedge}. In: Miklowitz, J.
  (Ed.), Wave Propag. Solids. ASME, Los Angeles, pp. 3--43.

\bibitem[{Komech et~al.(2015)Komech, Merzon, and Mendez}]{KMM2015}
Komech, A.~I., Merzon, A.~E., Mendez, J. E. D. l.~P., 2015. {Time-dependent
  scattering of generalized plane waves by a wedge}. Math. Methods Appl. Sci.
  38~(18), 4774--4785.

\bibitem[{Kontorovich and Lebedev(1939)}]{KL1939}
Kontorovich, M.~J., Lebedev, N.~N., 1939. {On a method of solution of some
  problems of the diffraction theory}. J. Phys. (Academy Sci. U.S.S.R.) 1~(3),
  229--241.

\bibitem[{Kouyoumjian and Pathak(1974)}]{KP1974}
Kouyoumjian, R.~G., Pathak, P.~H., 1974. {A uniform GTD for an edge in a
  perfectly conducting surface}. Proc. IEEE 62~(11), 1448--1461.

\bibitem[{Kythe(2011)}]{greensfunctionbook}
Kythe, P.~K., 2011. {Green's Functions and Linear Differential Equations}.
  Chapman and Hall/CRC.

\bibitem[{Lawrie and Abrahams(2007)}]{LawrieAbrahams2007}
Lawrie, J.~B., Abrahams, I.~D., 2007. {A brief historical perspective of the
  Wiener-Hopf technique}. J. Eng. Math. 59~(4), 351--358.

\bibitem[{Lebedev(1965)}]{Lebedev1965}
Lebedev, N.~N., 1965. {Special Functions and Their Application}. Eaglewood
  Cliffs.

\bibitem[{Lyalinov(1999)}]{Lyalinov1999}
Lyalinov, M.~A., 1999. {Diffraction by a highly contrast transparent wedge}. J.
  Phys. A. Math. Gen. 32~(11), 2183--2206.

\bibitem[{Lyalinov(2013)}]{Lyalinov2013}
Lyalinov, M.~A., 2013. {Scattering of acoustic waves by a sector}. Wave Motion
  50~(4), 739--762.

\bibitem[{Lyalinov and Zhu(2013)}]{LyalinovZhu2013}
Lyalinov, M.~A., Zhu, N.~Y., 2013. {Scattering of Waves by Wedges and Cones
  with Impedance Boundary Conditions}. Scitech, Edison, New Jersey.

\bibitem[{Macdonald(1902)}]{Macdonald1902}
Macdonald, H.~M., 1902. {Electric Waves}. University Press.

\bibitem[{Malyuzhinets(1955{\natexlab{a}})}]{Malyuzhinets1955-1-russian}
Malyuzhinets, G.~D., 1955{\natexlab{a}}. {Sound radiation by vibrating faces of
  arbitrary wedge. Part I}. Akust. Zhurnal (in Russ.) 1~(2), 144--164.

\bibitem[{Malyuzhinets(1955{\natexlab{b}})}]{Malyuzhinets1955-2-russian}
Malyuzhinets, G.~D., 1955{\natexlab{b}}. {Sound radiation by vibrating faces of
  arbitrary wedge. Part II}. Akust. Zhurnal (in Russ.) 1~(3), 226--234.

\bibitem[{Malyuzhinets(1958{\natexlab{a}})}]{Malyuzhinets1958-3}
Malyuzhinets, G.~D., 1958{\natexlab{a}}. {Excitation, reflection and emission
  of surface waves from a wedge with given face impedances}. Sov. Phys. Dokl.
  3, 752--755.

\bibitem[{Malyuzhinets(1958{\natexlab{b}})}]{Malyuzhinets1958-1}
Malyuzhinets, G.~D., 1958{\natexlab{b}}. {Inversion formula for Sommerfeld
  integral}. Sov. Phys. Dokl. 3, 52--56.

\bibitem[{Malyuzhinets(1958{\natexlab{c}})}]{Malyuzhinets1958-2}
Malyuzhinets, G.~D., 1958{\natexlab{c}}. {Relation between the inversion
  formulas for the Sommerfeld Integral and the formulas of
  Kontorovich-Lebedev}. Sov. Phys. Dokl. 3, 266--268.

\bibitem[{Martin and Wickham(1983)}]{Martin1983}
Martin, P.~A., Wickham, G.~R., 1983. {Diffraction of Elastic Waves by a
  Penny-Shaped Crack: Analytical and Numerical Results}. Proc. R. Soc. A Math.
  Phys. Eng. Sci. 390~(1798), 91--129.

\bibitem[{Mcnamara et~al.(1990)Mcnamara, Pistorius, and
  Malherbe}]{Pistorius1990}
Mcnamara, D.~A., Pistorius, C. W.~I., Malherbe, J. A.~G., 1990. {Introduction
  to the Uniform Geometrical Theory of Diffraction}. Artech House.

\bibitem[{Miles(1952)}]{Miles1952-1}
Miles, J.~W., 1952. {On the diffraction of an acoustic pulse by a wedge}. Proc.
  R. Soc. A Math. Phys. Eng. Sci. 212~(1111), 543--547.

\bibitem[{Moran et~al.(2016)Moran, Biggs, and Chamberlain}]{Moran2016}
Moran, C. A.~J., Biggs, N. R.~T., Chamberlain, P.~G., 2016. {Embedding formulae
  for wave diffraction by a circular arc}. Wave Motion 67, 32--46.

\bibitem[{Noble(1958)}]{Noble1958}
Noble, B., 1958. {Methods based on the Wiener-Hopf Technique for the solution
  of partial differential equations}, 2nd Edition. Chelsea Publishing Company,
  New York.

\bibitem[{Oberhettinger(1954)}]{Oberhettinger1954}
Oberhettinger, F., 1954. {Diffraction of Waves by a Wedge}. Commun. Pure Appl.
  Math. 7~(3), 551--563.

\bibitem[{Oberhettinger(1958)}]{Oberhettinger1958}
Oberhettinger, F., 1958. {On the diffraction and reflection of waves and pulses
  by wedges and corners}. J. Res. Natl. Bur. Stand. 61~(5), 343--365.

\bibitem[{Osipov and Norris(1999)}]{NorrisOsipov1999}
Osipov, A.~V., Norris, A.~N., 1999. {The Malyuzhinets theory for scattering
  from wedge boundaries: a review}. Wave Motion 29, 313--340.

\bibitem[{Poincar{\'{e}}(1892)}]{Poincare1892part1}
Poincar{\'{e}}, H., 1892. {Sur la Polarisation par Diffraction}. Acta Math. (in
  French) 16, 297--339.

\bibitem[{Poincar{\'{e}}(1897)}]{Poincare1892part2}
Poincar{\'{e}}, H., 1897. {Sur la Polarisation par Diffraction: Seconde
  partie}. Acta Math. (in French) 20, 313--355.

\bibitem[{Rawlins(1987)}]{Rawlins1987}
Rawlins, A.~D., 1987. {Plane-wave diffraction by a rational wedge}. Proc. R.
  Soc. A Math. Phys. Eng. Sci. 411, 265--283.

\bibitem[{Rawlins(1989)}]{Rawlins1989}
Rawlins, A.~D., 1989. {A Green's function for diffraction by a rational wedge}.
  Math. Proc. Cambridge Philos. Soc. 105~(1), 185--192.

\bibitem[{Rawlins(1999)}]{Rawlins1999}
Rawlins, A.~D., 1999. {Diffraction by, or diffusion into, a penetrable wedge}.
  Proc. R. Soc. A Math. Phys. Eng. Sci. 455, 2655--2686.

\bibitem[{Schot(1992)}]{SCHOT1992}
Schot, S., 1992. {Eighty years of Sommerfeld's radiation condition}. Historia
  Mathematica 19~(4), 385--401.

\bibitem[{Senior(1959)}]{Senior1959}
Senior, T. B.~A., 1959. {Diffraction by an Imperfectly Conducting Wedge}.
  Commun. Pure Appl. Math. 12~(2), 337--372.

\bibitem[{Shanin(1996)}]{AVShanin1996}
Shanin, A.~V., 1996. {On Wave Excitation in a Wedge-Shaped Region}. Acoust.
  Phys. 42~(5), 612--617.

\bibitem[{Shanin(1998)}]{AVShanin1998}
Shanin, A.~V., 1998. {Excitation of Waves in a Wedge-Shaped Region}. Acoust.
  Phys. 44~(5), 592--597.

\bibitem[{Shanin(2005)}]{shanin1}
Shanin, A.~V., 2005. {Modified Smyshlyaev's formulae for the problem of
  diffraction of a plane wave by an ideal quarter-plane}. Wave Motion 41~(1),
  79--93.

\bibitem[{Shanin and Craster(2010)}]{ShaninPseudoDiff}
Shanin, A.~V., Craster, R.~V., 2010. {Pseudo-differential operators for
  embedding formulae}. J. Comput. Appl. Math. 234, 1637--1646.

\bibitem[{Skelton et~al.(2010)Skelton, Craster, Shanin, and
  Valyaev}]{wedgecraster}
Skelton, E.~A., Craster, R.~V., Shanin, A.~V., Valyaev, V.~Y., 2010. {Embedding
  formulae for scattering by three-dimensional structures}. Wave Motion 47~(5),
  299--317.

\bibitem[{Smirnov(1964)}]{Smirnov1964}
Smirnov, V.~I., 1964. {A Course of Higher Mathematics. Volume 3, Part 2.}
  Pergamon Press.

\bibitem[{Sobolev(1935)}]{Sobolev1935}
Sobolev, S.~L., 1935. {General Theory of Diffraction of Waves on Riemann
  Surfaces}. In: Sel. Work. S. L. Sobolev Vol I. New York, pp. 201--262.

\bibitem[{Sommerfeld(1896)}]{Sommerfeld1896}
Sommerfeld, A., 1896. {Mathematische Theorie der Diffraction}. Math. Ann. (in
  Ger.) 47, 317--374.

\bibitem[{Sommerfeld(1901)}]{Sommerfeld1901}
Sommerfeld, A., 1901. {Theoretical about the diffraction of X-rays}.
  Zeitschrift f{\"{u}}r Math. und Phys. (in Ger.) 46, 11--97.

\bibitem[{Sommerfeld(2003)}]{Sommerfeld2003}
Sommerfeld, A., 2003. {Mathematical Theory of Diffraction}. Birkhauser, Boston.

\bibitem[{Teixeira(1991)}]{Teixeira1991}
Teixeira, F., 1991. {Diffraction by a rectangular wedge: Wiener-Hopf-Hankel
  formulation}. Integr. Equations Oper. Theory 14~(3), 436--454.

\bibitem[{Ufimtsev(1971)}]{Ufimtsev1971}
Ufimtsev, P.~Y., 1971. {Method of edge waves in the physical theory of
  diffraction}. (from the Russian ``Metod krayevykh voln v fizicheskoy teorii
  difraktsii'' Izd-Vo Sov. Radio, pp-243, 1962 ), translation prepared by the
  U.S. Air Force Foreign Technology Division.

\bibitem[{Ufimtsev(2014)}]{Ufimtsev2014}
Ufimtsev, P.~Y., 2014. {Fundamentals of the Physical Theory of Diffraction},
  2nd Edition. John Wiley {\&} Sons.

\bibitem[{Voss(2013)}]{Voss_2013}
Voss, J., 2013. An Introduction to Statistical Computing: A Simulation-based
  Approach. Wiley.

\bibitem[{Wegert(2012)}]{Wegert2012}
Wegert, E., 2012. {Visual Complex Functions}. Birkhauser Basel.

\bibitem[{Wiener and Hopf(1931)}]{WienerHopf1931}
Wiener, N., Hopf, E., 1931. {{\"{U}}ber eine klasse singul{\"{a}}rer
  integralgleichungen}. Sitzungsberichte der Preuss. Akad. der Wissenschaften,
  Phys. Klasse (in Ger.) 31, 696--706.

\bibitem[{Williams(1982)}]{williams}
Williams, M.~H., 1982. {Diffraction by a finite strip}. Q. J. Mech. Appl. Math.
  35, 103--124.

\bibitem[{Williams(1959)}]{Williams1959}
Williams, W.~E., 1959. {Diffraction of an E-polarized plane wave by an
  imperfectly conducting wedge}. Proc. R. Soc. A Math. Phys. Eng. Sci.
  252~(1270), 376--393.

\end{thebibliography}
}


\appendix
\mylinenum{
\section{Macdonald's series solution}\label{Mac} 
In this section we shall briefly discuss the separation of variables method applied by \citet{Macdonald1902} to the wedge problem with line source incidence. After this, departing slightly from Macdonald's approach, we shall use a limiting procedure in order to recover the series solutions \eqref{MKLT-SeriesSol-D} and \eqref{MKLT-SeriesSol-N} to the plane wave incidence problem. 

The wedge problem forced by a line source of strength $\mathcal{A}$ with polar coordinates $(r_\textrm{I},\theta_\textrm{I})$, has the following governing equation,}
\begin{align}
\label{Mac-goveqn}\nabla^2\Phi+k^2\Phi=\frac{\mathcal{A}}{r}\hat{\delta}(r-r_\textrm{I})\hat{\delta}(\theta-\theta_\textrm{I}),
\end{align}
\mylinenum{where $\hat{\delta}$ is the Dirac delta function. The total field, $\Phi$, is decomposed into incident and scattered parts $\Phi=\Phi_\textrm{I}+\Phi_\textrm{S}$ where the incident wave is given by}
\begin{align}
\label{Mac-inc-linesource}\Phi_\textrm{I}=\frac{\mathcal{A}}{4i}H^{(1)}_0\left(k\sqrt{r^2+r_\textrm{I}^2-2rr_\textrm{I}\cos(\theta-\theta_\textrm{I})}\right),
\end{align}
\mylinenum{and is subjected to BCs, \eqref{Intro-DBC} or \eqref{Intro-NBC}. Considering the ansatz $\Phi=R(r)\Theta(\theta)$, using separation of variables and applying the BCs, we obtain the following series solutions:}
\begin{align}
\label{Mac-preDir}\textrm{Dirichlet case:}\quad\Phi(r,\theta)&=\sum_{n=1}^\infty A_nR_n(r)\sin((\theta-\theta_\textrm{w})\delta n),\\
\label{Mac-preNeu}\textrm{Neumann case:}\quad\Phi(r,\theta)&=\sum_{n=0}^\infty B_nR_n(r)\cos((\theta-\theta_\textrm{w})\delta n).
\end{align}
\mylinenum{Because of the source location at $r=r_\textrm{I}$, and the need to satisfy both the edge and the radiation conditions (satisfied by the Bessel and Hankel functions respectively), we pose}
\begin{align}
\label{Mac-R-sol}R_n(r)=\begin{cases}
C_nJ_{\delta n}(kr)& r<r_\textrm{I},\\
D_nH^{(1)}_{\delta n}(kr)& r>r_\textrm{I}.
\end{cases}
\end{align}
\mylinenum{To ensure continuity across $r=r_\textrm{I}$, we require $C_n=H^{(1)}_{\delta n}(kr_\textrm{I})$ and $D_n=J_{\delta n}(kr_\textrm{I})$. We can determine the coefficients $A_n$ and $B_n$ by deriving and applying a jump condition across $r=r_\textrm{I}$. 

In the Dirichlet case, substitute \eqref{Mac-preDir} into \eqref{Mac-goveqn}, and multiply the resulting equation by $r\sin((\theta-\theta_\textrm{w})\delta m)$. Integrating w.r.t. $\theta$ from $-\theta_\text{w}$ to $\theta_\text{w}$, and using the orthogonality of sine, we obtain}
\begin{align}
\label{Mac-Dir-jump(1)}A_m\ParDer{}{r}(rR_m'(r))+A_m\left(k^2-\frac{\delta^2m^2}{r^2}\right)rR_m(r)=\frac{\mathcal{A}}{\theta_\textrm{w}}\hat{\delta}(r-r_\textrm{I})\sin((\theta_\textrm{I}-\theta_\textrm{w})\delta m).
\end{align}
\mylinenum{Now integrating \eqref{Mac-Dir-jump(1)} from $r=r_\textrm{I}-\epsilon$ to $r_\textrm{I}+\epsilon$ and taking the limit $\epsilon\rightarrow 0$ leads to the jump condition}
\begin{align}
\label{Mac-Dir-jump(2)}A_mr_\textrm{I}\left[R_m'(r)\right]_{r_\textrm{I}^-}^{r_\textrm{I}^+}=\frac{\mathcal{A}}{\theta_\textrm{w}}\sin((\theta_\textrm{I}-\theta_\textrm{w})\delta m).
\end{align}
\mylinenum{Lastly, we use \eqref{Mac-R-sol} and the Wronskian result $J_{\nu}(z)H^{(1)'}_{\nu}(z)-J_{\nu}'(z)H^{(1)}_{\nu}(z)=\frac{2i}{\pi z}$
to determine that $A_n=-i\delta\mathcal{A}\sin((\theta_\textrm{I}-\theta_\textrm{w})\delta n)$. Hence, the series solution with line source incidence and Dirichlet BCs is}
\begin{align}
\label{Mac-Dir-preseries}\Phi=\sum_{n=1}^{\infty}i\delta\mathcal{A}\sin((\theta_\textrm{w}-\theta_\textrm{I})\delta n)\sin((\theta-\theta_\textrm{w})\delta n)J_{\delta n}(kr_{<})H^{(1)}_{\delta n}(kr_{>}),
\end{align}
\mylinenum{where $r_{<}=\min(r,r_\text{I})$ and $r_{>}=\max(r,r_\text{I})$. This agrees with Macdonald's solution\footnote{Note that \citet{Macdonald1902} uses the alternate time factor $e^{i\omega t}$.}.

For the Neumann case, the coefficients $B_n$ are found by the same method using the orthogonality relation for cosine, leading to}
\begin{align}
 \Phi&=-\sum_{n=0}^{\infty}i\varepsilon_n\delta\mathcal{A}\cos((\theta_\textrm{w}-\theta_\textrm{I})\delta n)\cos((\theta-\theta_\textrm{w})\delta n)J_{\delta n}(kr_{<})H^{(1)}_{\delta n}(kr_{>}), \label{Mac-Neu-preseries}
\end{align}
\mylinenum{where $\varepsilon_0=1/2$ and $\varepsilon_n=1$ for $n\geq 1$.
 
To recover the plane wave solution, we send the source and its strength to infinity in a way that ensures that $\Phi_\text{I}$ (as defined in \eqref{Mac-inc-linesource}) behaves like $e^{-ikr\cos(\theta-\theta_\textrm{I})}$ as $r_\text{I}\rightarrow\infty$. This can be done by choosing $\mathcal{A}=\sqrt{8\pi k r_\textrm{I}}e^{-ikr_\textrm{I}+\frac{3\pi i}{4}}$ and leads to}
\begin{align}
\label{Mac-Dir-coeff-limit}\lim_{r_\text{I}\rightarrow\infty}A_nH^{(1)}_{\delta n}(kr_\textrm{I})&=4\delta(-i)^{\delta n}\sin((\theta_\textrm{I}-\theta_\textrm{w})\delta n),\\
\label{Mac-Neu-coeff-limit}\lim_{r_\text{I}\rightarrow\infty}B_nH^{(1)}_{\delta n}(kr_\textrm{I})&=4\varepsilon_n\delta(-i)^{\delta n}\cos((\theta_\textrm{w}-\theta_\textrm{I})\delta n).
\end{align}
\mylinenum{Hence, for plane wave forcing with Dirichlet or Neumann BCs respectively, the series solutions are}
\begin{align}
\nonumber\Phi(r,\theta)&=4\delta\sum_{n=1}^\infty(-i)^{\delta n}J_{\delta n}(kr)\sin((\theta-\theta_\textrm{w})\delta n)\sin((\theta_\textrm{I}-\theta_\textrm{w})\delta n),\\
\nonumber\Phi(r,\theta)&=2\delta J_0(kr)+4\delta\sum_{n=1}^\infty(-i)^{\delta n}J_{\delta n}(kr)\cos((\theta-\theta_\textrm{w})\delta n)\cos((\theta_\textrm{w}-\theta_\textrm{I})\delta n),
\end{align}
\mylinenum{which matches perfectly with \eqref{MKLT-SeriesSol-D} and \eqref{MKLT-SeriesSol-N} as required. Note that these exact series solutions have a natural embedding structure (see Section \ref{Embedding}) in the sense that they are simply sums of products of functions of one variable only.


\section{A link between the spectral function $s(z)$ and Green's integral operator}\label{app:app1} 

\subsection{Preliminary definitions and Green's function representation}

Let us introduce the generic plane wave function $w_z$ by}
\begin{eqnarray}
w_z(r,\theta)=\exp(ikr\cos(z-\theta))
\end{eqnarray}
\mylinenum{It is important to note that for any $\Theta\in[-\pi,\pi]$, $w_z(r,\Theta)$ is exponentially decaying as $r \rightarrow \infty$ as long as $z\in\Omega_{\Theta} = \Omega_0 + \Theta$, where both $\Omega_0$ and $\Omega_{\Theta}$ are understood as open sets (do not contain their boundaries) and are illustrated in Figure \ref{fig:domainsOmegastrips}.

\begin{figure}[h]\centering
	\includegraphics[width=0.4\textwidth]{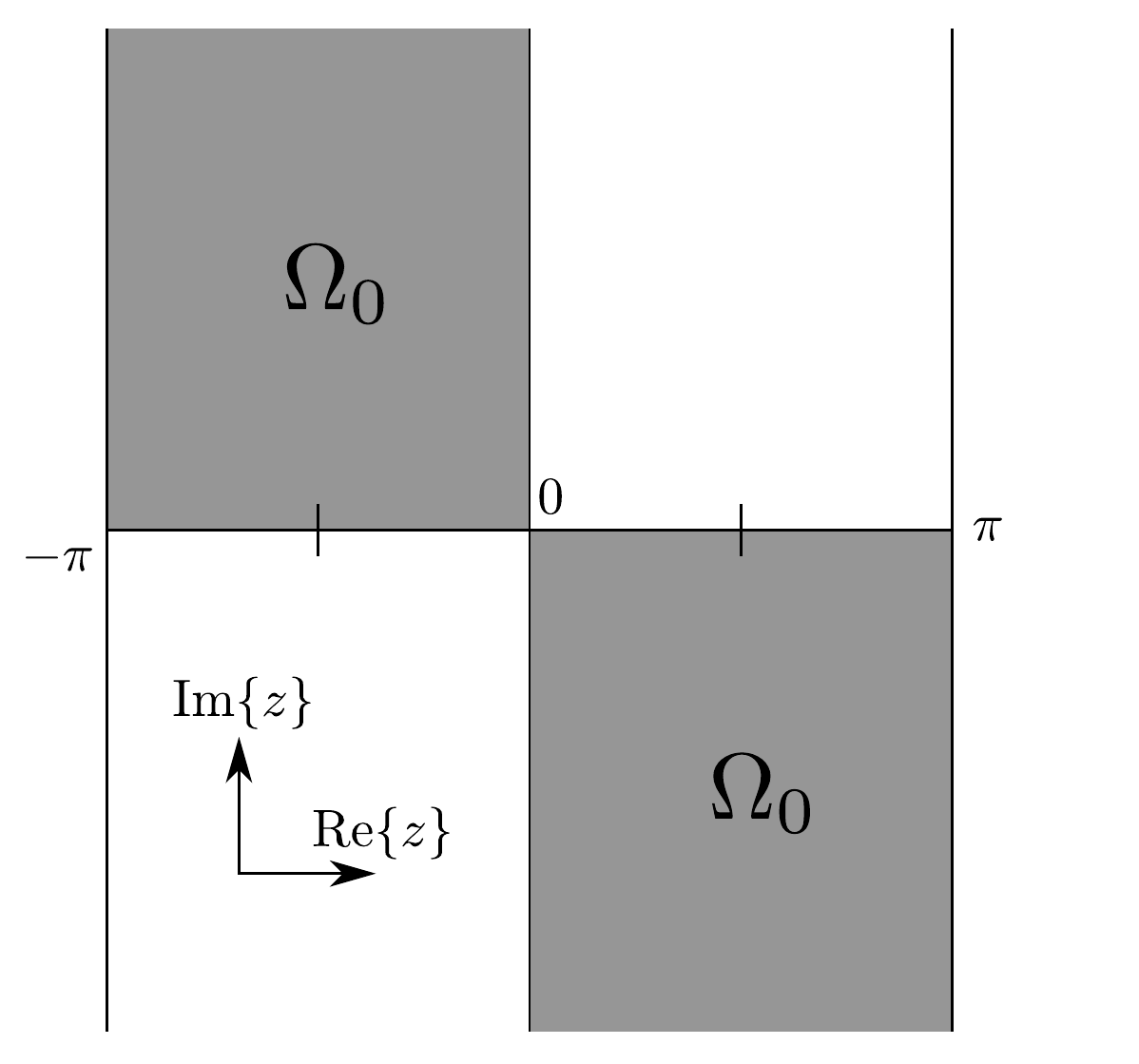}\quad
	\includegraphics[width=0.4\textwidth]{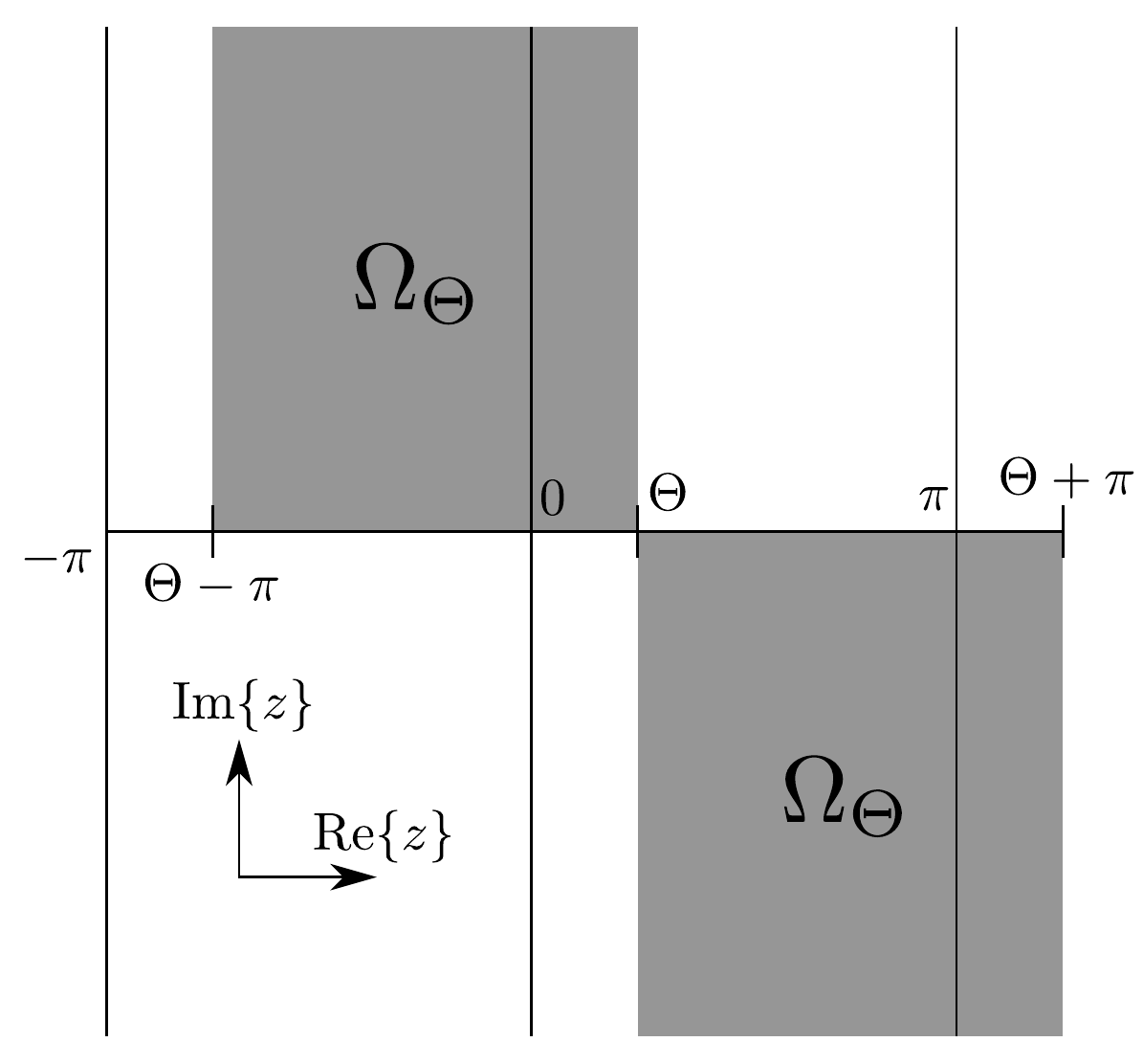}
	\caption{The domains $\Omega_0$ and $\Omega_{\Theta}$}
	\label{fig:domainsOmegastrips}
\end{figure}

Using the notations of Section \ref{SolAnalysis}, the total, diffracted and geometrical optics fields are denoted $\Phi$, $\Phi_{\tmop{Diff}}$ and $\Phi_{\tmop{GO}}$. For the exterior wedge, $\Phi_{\tmop{GO}}$ consists of an incident wave and one or two reflected waves and can hence be written in the form $\Phi_{\tmop{GO}}(r,\theta)=\sum_i a_i(\theta)w_{z_i}(r,\theta),$ where $a_i (\theta)$ is either zero or a given constant. Typically in our problem the incident wave corresponds to $z_i = \theta_{\text{I}} + \pi$, and the reflected waves to either $z_i=2\theta_{\text{w}}-\theta_{\text{I}}+\pi$ or $z_i=-2\theta_{\text{w}}-\theta_{\text{I}}+\pi$ or both depending on how many reflections we have.

For a given function $\Psi$, let us now introduce the Green integral operator $\mathcal{S}_{\Theta} (z) [\Psi]$ defined by}
\begin{eqnarray}
\mathcal{S}_{\Theta}(z)[\Psi]=\int_0^{\infty}\left[\Psi\frac{\partial w_z}{\partial \theta}-\frac{\partial \Psi}{\partial \theta}w_z\right]_{\theta=\Theta}\frac{\mathd r}{r}.
\label{eq:GreenIntegralOperator}\end{eqnarray}
\mylinenum{In the contest of this review, $\Theta\in\left[-\theta_{\text{w}},\theta_{\text{w}}\right]$. Moreover, using standard integration, one can show that $\mathcal{S}_{\Theta}(z)[w_{z_i}]=\tan\left(\frac{z_i-z}{2}\right)$, and hence we can write}
\begin{eqnarray}
\mathcal{S}_{\Theta}(z)[\Phi]=\mathcal{S}_{\Theta}(z)[\Phi_{\tmop{Diff}}]+\sum_i a_i(\Theta)\tan\left(\frac{z_i-z}{2}\right), 
\end{eqnarray}
\mylinenum{which implies in particular that each of the $z_i+\pi$ are simple poles of $\mathcal{S}_{\Theta}(z)[\Phi]$ with residue $- 2$. The objective of this appendix is to find a connection between the spectral function $s (z)$ and $\mathcal{S}_{\Theta} (z) [\Phi]$. In order to do that, we shall make use of the theory of Green's functions as follows.

Let us pick a point $(r^{\star}, \theta^{\star})$, and pick two angles $\varphi_a$ and $\varphi_b$ (the subscripts $a$ and $b$ stand for {\tmem{above}} and {\tmem{below}}) chosen such that $-\pi/2<\varphi_b<\theta^{\star}<\varphi_a<\pi/2$, and a radius $R_{\mathcal{A}}>r^{\star}$. Now consider the domain $\Omega^{\star}(\varphi_b,\varphi_a,R_{\mathcal{A}})$ to be the corresponding sector described in Figure \ref{fig:Angular_domain}. Let us further assume that $\partial \Omega^{\star}$
is oriented anti-clockwise, and that the normals $\tmmathbf{n}$ to $\partial\Omega^{\star}$ are chosen to be outgoing.

\begin{figure}[h]
\centering
  \includegraphics[width=0.25\textwidth]{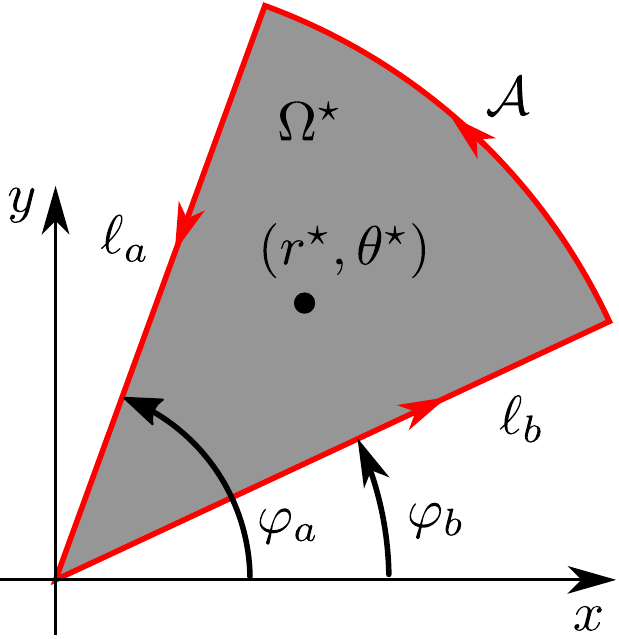}
  \caption{The angular domain $\Omega^{\star}$}
\label{fig:Angular_domain}
\end{figure}

Let $G_{\star} (r, \theta)$ be a short notation for $G(r,\theta;r^{\star},\theta^{\star})$, the free-space Green's function for the Helmholtz equation resulting from a point source at $(r^{\star}, \theta^{\star})$. Using the respective governing equations of $\Phi$ and $G_{\star}$, and the divergence theorem, we have}
\begin{align}\nonumber
\iint_{\Omega^{\star}}(\Phi\Delta G_{\star}-G_{\star}\Delta\Phi)\mathd A=\Phi(r^{\star},\theta^{\star})=\oint_{\partial \Omega^{\star}}(\Phi\nabla G_{\star}-G_{\star}\nabla\Phi)\cdot\tmmathbf{n}\mathd s
\end{align}
\mylinenum{Hence, using that on $\ell_b$, $\tmmathbf{n}=-\tmmathbf{e}_{\theta}$, on $\mathcal{A}$, $\tmmathbf{n}=\tmmathbf{e}_r$ and on $\ell_a$, $\tmmathbf{n}=\tmmathbf{e}_{\theta}$, we get}
\begin{align}
\Phi(r^{\star},\theta^{\star})=&\ 
\underbrace{-\int_0^{R_{\mathcal{A}}}\left(\Phi\frac{\partial G_{\star}}{\partial \theta}-G_{\star}\frac{\partial\Phi}{\partial\theta}\right)_{\theta=\varphi_b}\frac{\mathd r}{r}}_{\ell_b\text{ component: }I_{\ell_b}[\Phi]}+
\underbrace{\int_{\varphi_b}^{\varphi_a}\left(\Phi\frac{\partial G_{\star}}{\partial r}-G_{\star}\frac{\partial\Phi}{\partial r}\right)_{r=R_{\mathcal{A}}}R_{\mathcal{A}}\mathd\theta}_{\text{Arc $\mathcal{A}$ component: }I_{\mathcal{A}}[\Phi]}\nonumber\\
&\underbrace{+\int_0^{R_{\mathcal{A}}}\left(\Phi\frac{\partial G_{\star}}{\partial\theta}-G_{\star}\frac{\partial\Phi}{\partial\theta}\right)_{\theta=\varphi_a}\frac{\mathd r}{r}}_{\ell_a\text{ component: $I_{\ell_a}[\Phi]$}}.
\label{eq:finitesliceofcake}\end{align}
\mylinenum{Using the Hankel representation of $G_{\star}$, its far-field asymptotics, and the method of steepest descent, we can show that the only part of the far-field leading to any contribution of the arc integral as $R_{\mathcal{A}}\rightarrow\infty$ is an incident plane wave coming from within the sector. More precisely, if $\theta_{\text{I}} \in (\varphi_b, \varphi_a)$,}
\begin{align}
\nonumber\lim_{R_A\rightarrow\infty} I_{\mathcal{A}}\left[e^{-ikr\cos\left(\theta-\theta_{\text{I}}\right)}\right]
=e^{-ikr^{\star}\cos\left(\theta^{\star}-\theta_{\text{I}}\right)}=\Phi_{\text{I}}(r^{\star},\theta^{\star}). 
\end{align}
\mylinenum{All other components (reflected waves, diffracted field, incident waves from outside the sector) can be shown to have zero contribution. Hence, taking the limit as $R_{\mathcal{A}} \rightarrow \infty$ in \eqref{eq:finitesliceofcake}, we get}
\begin{align}
\Phi(r^{\star},\theta^{\star})=-\int_0^{\infty}\left(\Phi\frac{\partial G_{\star}}{\partial \theta}-G_{\star}\frac{\partial\Phi}{\partial\theta}\right)_{\theta=\varphi_b}\frac{\mathd r}{r}
+\int_0^{\infty}\left(\Phi\frac{\partial G_{\star}}{\partial\theta}-G_{\star}\frac{\partial\Phi}{\partial\theta}\right)_{\theta=\varphi_a}\frac{\mathd r}{r}+\Phi^{ab}_{\text{I}}(r^{\star},\theta^{\star}),
\label{eq:phidiffgreen}\end{align}
\mylinenum{where $\Phi_{\text{I}}^{ab}=\Phi_{\text{I}}$ if $\theta_{\text{I}}\in(\varphi_a,\varphi_b)$ and zero otherwise. Hence the knowledge of $G_{\star}$ and $\frac{\partial G^{\star}}{\partial \theta}$ on oblique lines of constant $\theta$ is important. At this stage, it is important to realise, at least informally, that if we could write them in terms of $w_z$ somehow, then we have a chance to link $\Phi$ and the Green integral operator.

\subsection{Green's functions on oblique lines}

Before finding formulae for $G_{\star}$, we will focus on the Green's function $G_0$ corresponding to a point source at the origin. First of all, it is well known (see e.g. {\cite{greensfunctionbook}}) that $G_0(r,\theta)=\frac{-i}{4}H_0^{(1)}(kr)$. Moreover, the Hankel function has the following integral representation\footnote{See, e.g. {\cite{Sommerfeld2003}} eq (6) p19, together with translators' note 4 on p78, here we use $\beta=\frac{\pi}{2}$.}}
\begin{align}
H_0^{(1)}(r)=\frac{1}{\pi}\int_{\Gamma}e^{ir\cos(z)}\mathd z\ \ \text{leading to}\ \ G_0(r,\theta)=\frac{+1}{4\pi i}\int_{\Gamma}e^{ikr\cos(z)}\mathd z,
\label{eq:intitiaschafli}\end{align}
\mylinenum{where $\Gamma$ is described in Figure \ref{SA-SDCs}. We will now endeavour to find formulae for $G_0$ valid on an oblique half-space and hence on any line that crosses the $x$ axis with an angle $\varphi\in(-\pi/2,\pi/2)$ say, and lies above (see Figure \ref{fig:lines} (left)) or below (see Figure \ref{fig:lines} (right)) the origin.

\begin{figure}[h]\centering
	\includegraphics[width=0.4\textwidth]{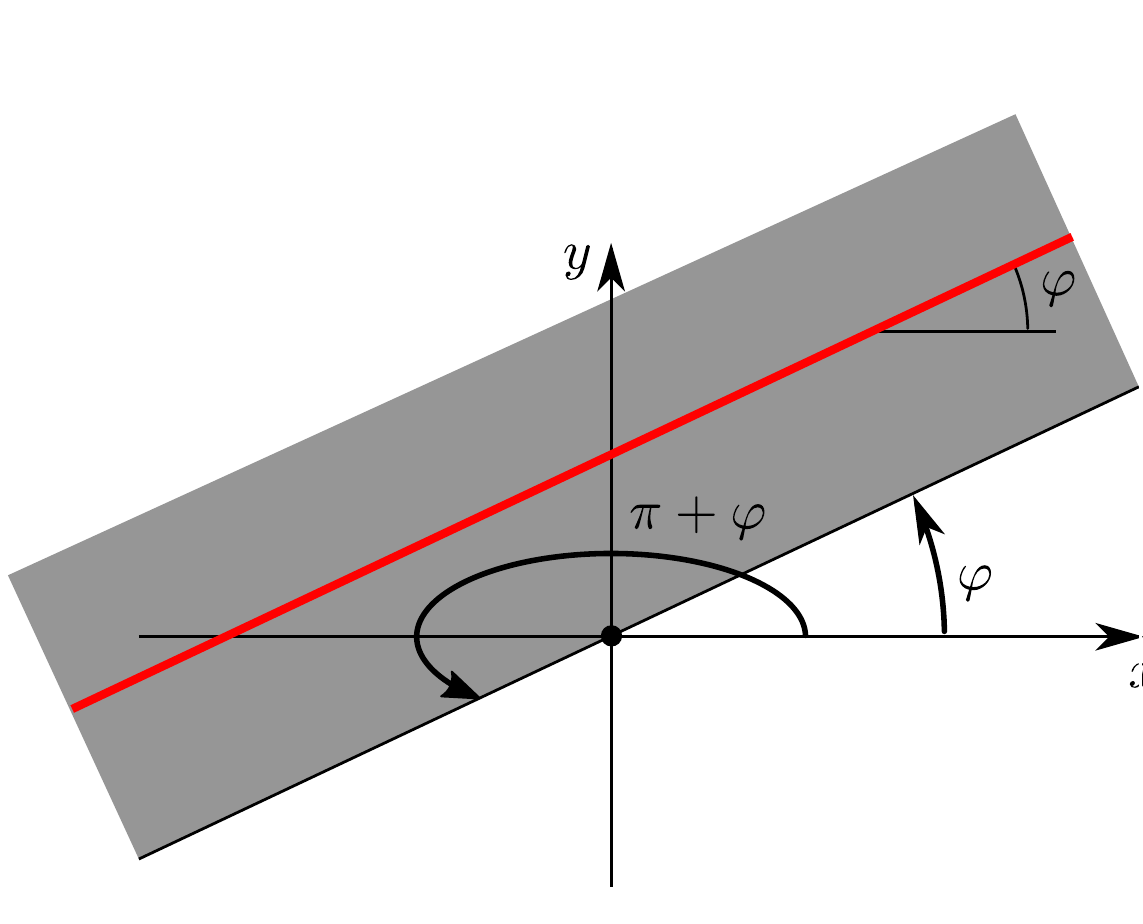}\quad
	\includegraphics[width=0.4\textwidth]{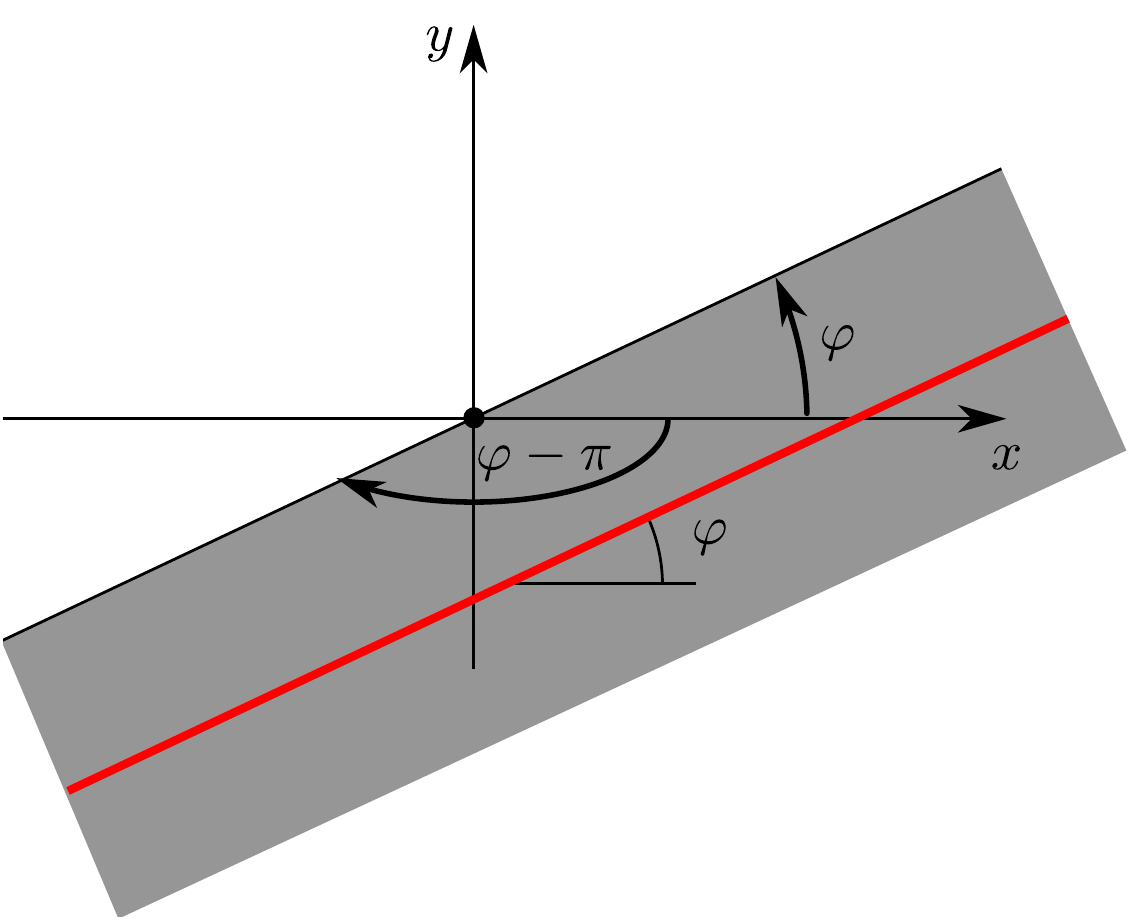}
	\caption{The two half-spaces under consideration for a given $\varphi$: above (left) and below (right)}
	\label{fig:lines}
\end{figure}

\paragraph{Oblique line above the origin}Let us consider the half-space $\varphi<\theta<\pi+\varphi$, the grey area of Figure \ref{fig:lines} (left) . Let us start from \eqref{eq:intitiaschafli} and shift the contour $\Gamma$ to the contour $\Gamma+\frac{\pi}{2}-(\theta-\varphi)$, where the new contour height is adjusted so that it goes through the origin. Because of the restriction on $\theta$, we can do that without leaving $\Omega_0$, where our integrand is analytic and exponentially decaying, and so the value of the integral and its convergence property remain unchanged. We can now perform the substitution $z'\leftrightarrow z+\theta$ to get}
\begin{eqnarray}
G_0(r,\theta)=\frac{1}{4\pi i}\int_{\gamma_a(\varphi;\theta)}w_{z'}(r,\theta)\mathd z', 
\end{eqnarray}
\mylinenum{where the contour $\gamma_a(\varphi;\theta)=\Gamma+\frac{\pi}{2}+\varphi$ goes through the point $z=\theta$ of the real axis, as shown in Figure \ref{fig:obliqueaboveadjusted} (left). This formula is valid (and the integral converges exponentially) on any oblique line with angle $\varphi$ that lies above the origin.

\begin{figure}[h]\centering
	\includegraphics[width=0.4\textwidth]{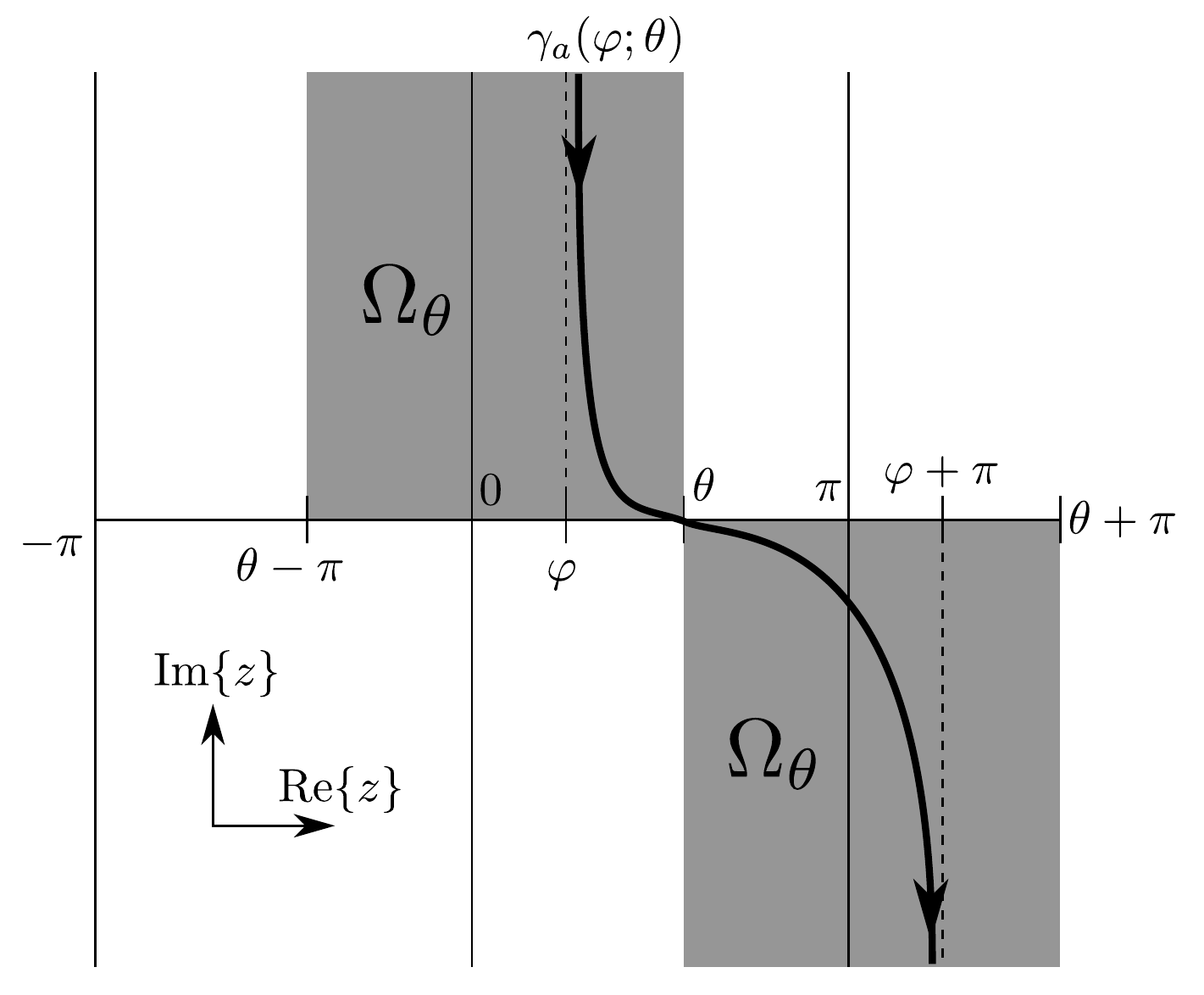}\quad
	\includegraphics[width=0.4\textwidth]{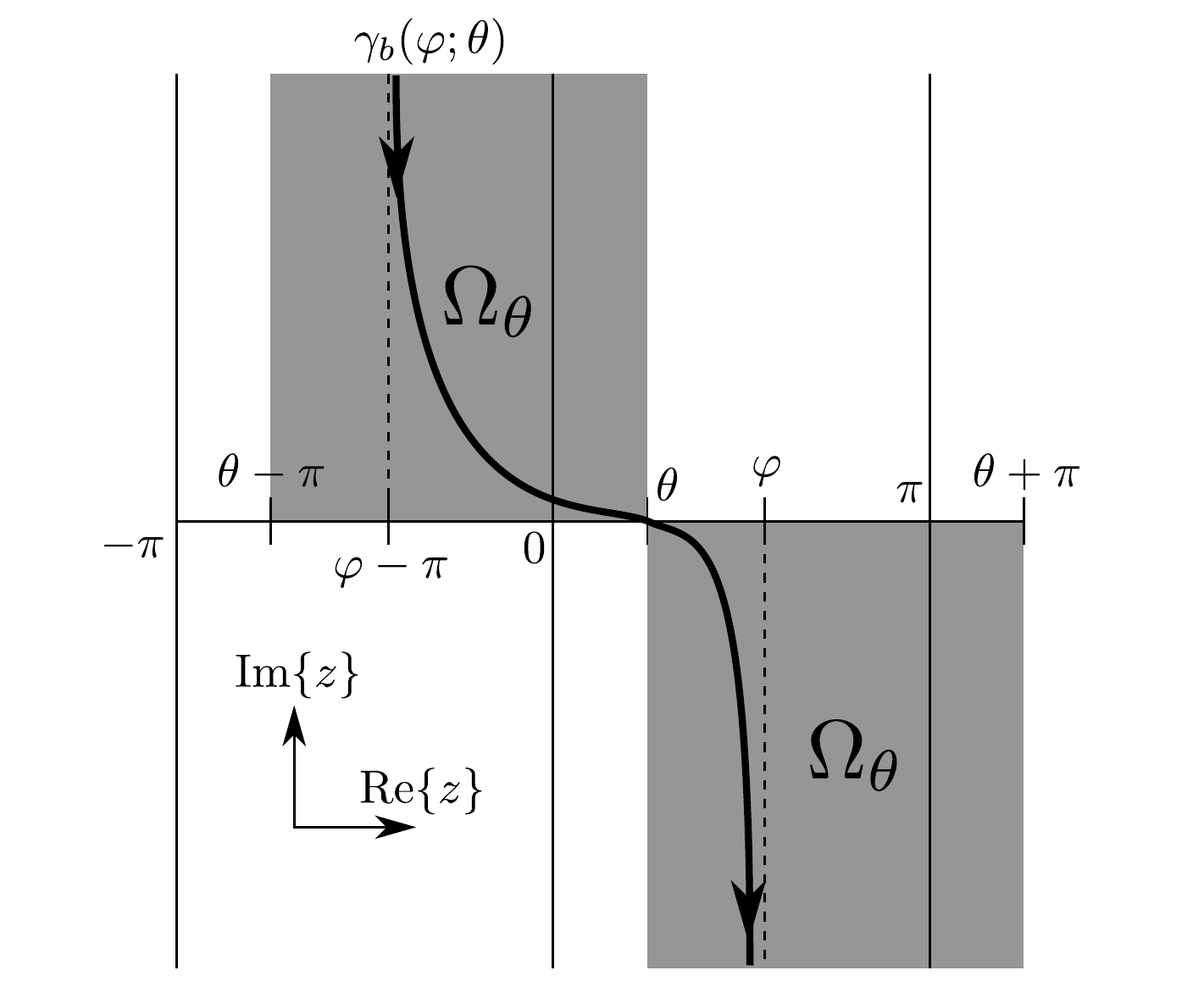}
	\caption{The contours $\gamma_a(\varphi;\theta)$ and $\gamma_b(\varphi;\theta)$}
	\label{fig:obliqueaboveadjusted}
\end{figure}

Note that here $r$ is finite, and the integrand is analytic, so we can in principle deform all the contours $\gamma_a(\varphi;\theta)$ to $\gamma_a(\varphi;\varphi)$, it is important to note that the latter crosses the real axis at $z=\varphi$ and is included (only just!) in $\Omega_{\varphi}$. This contour will just be referred to as $\gamma_a (\varphi)$ thereafter, and we get}
\begin{eqnarray}
G_0(r,\theta)=\frac{1}{4\pi i}\int_{\gamma_a(\varphi)}w_z(r,\theta)\mathd z
\label{eq:obliqueabove}\end{eqnarray}

\mylinenum{\paragraph{Oblique lines below the origin}In a very similar way, we can consider the half-space $\varphi-\pi<\theta<\varphi$ and shift the contour $\Gamma$ to a height adjusted $\Gamma-\frac{\pi}{2}+(\varphi-\theta)$ passing through the origin. Upon performing the substitution $z\leftrightarrow z+\theta$, we obtain an integral over a contour $\gamma_b(\varphi;\theta)$ illustrated on Figure \ref{fig:obliqueaboveadjusted} (right). Again by analyticity of the integrand, such integral can safely be deformed to the contour $\varphi_b(\varphi)\equiv\varphi_b(\varphi;\varphi)$ that crosses the real axis at $z = \varphi$ and lies within $\Omega_{\varphi}$, to get}
\begin{eqnarray}
G_0(r,\theta)=\frac{1}{4\pi i}\int_{\gamma_b(\varphi)}w_z(r,\theta)\mathd z,\label{eq:obliquebelow}
\end{eqnarray}
\mylinenum{\paragraph{Back to $G_{\star}$}In order to get back to $G_{\star}$, we just need to replace $r$ by $r'$ and $\theta$ by $\theta'$ in \eqref{eq:obliqueabove} and \eqref{eq:obliquebelow}, where $r'$ and $\theta'$ are the polar coordinates centred at $(r^{\star}, \theta^{\star})$. Upon noting that $r'e^{i \theta'}=re^{i\theta}-r^{\star}e^{i\theta^{\star}}$, we find that}
\begin{align}\nonumber
w_z(r',\theta')=w_z(r,\theta)e^{-ikr^{\star}\cos(z-\theta^{\star})}\ \ \text{and}\ \ G_{\star}(r,\theta)=\frac{1}{4\pi i}\int_{\gamma_s(\varphi_s)}w_z(r,\theta)e^{-ikr^{\star}\cos(z-\theta^{\star})}\mathd z,
\end{align}
\mylinenum{where from now on, the subscript $s$ is either $a$ or $b$. Since $\gamma_s(\varphi_s)$ is independent of $\theta$, we get a similar formula for $\frac{\partial G_{\star}}{\partial \theta}$. In particular, in the configuration of Figure \ref{fig:Angular_domain}, since the oblique line $\ell_a$ (resp. $\ell_b$) lies above (resp. below) the source $(r^{\star}, \theta^{\star})$ and make an angle $\varphi_a$ (resp. $\varphi_b$) with the real axis, we have}

\begin{eqnarray}
G_{\star}|_{\ell_s}=G_{\star}(r,\varphi_s) \ \text{ and } \
\frac{\partial G_{\star}}{\partial\theta}\bigg|_{\ell_s}=\frac{1}{4\pi i}\int_{\gamma_s(\varphi_s)}\frac{\partial w_z}{\partial \theta}(r,\varphi_s)e^{-ikr^{\star}\cos(z-\theta^{\star})}\mathd z.
\label{eq:restrictionGstar}\end{eqnarray}

\mylinenum{
\subsection{Connection formula between $s(z)$ and $S_0(z)$}

Before making use of our results \eqref{eq:phidiffgreen} and \eqref{eq:restrictionGstar}, we need to make use of some properties of the Green integral operator:

\begin{proposition}
  \label{prop:greenoperatorstuff}{\tmdummy}
 
  \begin{enumerate}
    \item Apart from eventual poles on the real line, as a function of $z$, $\mathcal{S}_{\Theta}(z)[\Phi]$ is analytic for $z\in\Omega_{\Theta}=\Omega_0+\Theta$.
    
    \item If $z\in\Omega_{\Theta_1}\cap\Omega_{\Theta_2}$, then $\mathcal{S}_{\Theta_1}(z)[\Phi]=\mathcal{S}_{\Theta_2}(z)[\Phi]$. Note that by analytic continuation, this allows to extend the natural domain of analyticity of $\mathcal{S}_{\Theta_{1,2}}(z)$ to $\Omega_{\Theta_1}\cup\Omega_{\Theta_2}$.
  \end{enumerate}
\end{proposition}

Now, we can input \eqref{eq:restrictionGstar} into \eqref{eq:phidiffgreen}, and, since we made sure that $\gamma_s(\varphi_s)\subset\Omega_{\varphi_s}$, we can exchange the order of integration. Let us furthermore assume that $\varphi_b<0<\varphi_a$, then the formula can be evaluated at $\theta^{\star} = 0$ to get}
\begin{eqnarray}
\Phi(r^{\star},0)=\frac{-1}{4\pi i}\int_{\gamma_b(\varphi_b)}\!\!\!\!\!\!e^{-ikr^{\star}\cos(z)}\mathcal{S}_{\varphi_b}(z)[\Phi]\mathd z
+\frac{1}{4\pi i}\int_{\gamma_a(\varphi_a)}\!\!\!\!\!\!e^{-ikr^{\star}\cos(z)}\mathcal{S}_{\varphi_a}(z)[\Phi]\mathd z+\Phi^{ab}_{\text{I}}(r^{\star},0),
\label{eq:PHIEVALON0}\end{eqnarray}
\mylinenum{where an illustration of the contour configuration is displayed in Figure \ref{fig:bothcontours} (left). Making use of point $2$ of Proposition \ref{prop:greenoperatorstuff}, the integrands of both integrals in \eqref{eq:PHIEVALON0} are actually analytical continuations of each other, and hence we can write}
\begin{align}\nonumber
\Phi(r^{\star},0)=\frac{1}{4\pi i}\int_{(\gamma_b(\varphi_b))^c+\gamma_a(\varphi_a)}e^{-ikr^{\star}\cos(z)}\mathcal{S}_{\varphi_b}(z)[\Phi]\mathd z+\Phi^{ab}_{\text{I}}(r^{\star},0),
\end{align}
\mylinenum{where $(\gamma_b (\varphi_b))^c$ is a notation for $\gamma_b (\varphi_b)$ going in the other direction. Let us consider the contours as angular (we can do that by analytic deformation), as depicted in Figure \ref{fig:bothcontours} (right). Let us also introduce a new contour $\mathcal{R}$, that is rectangular, with its centre at the origin and oriented anticlockwise, such that its left (resp. right) lateral side coincides with a part of $(\gamma_b(\varphi_b))^c$ (resp. $\gamma_a(\varphi_a)$), but in the opposite direction.

\begin{figure}[h]\centering
  \includegraphics[width=0.45\textwidth]{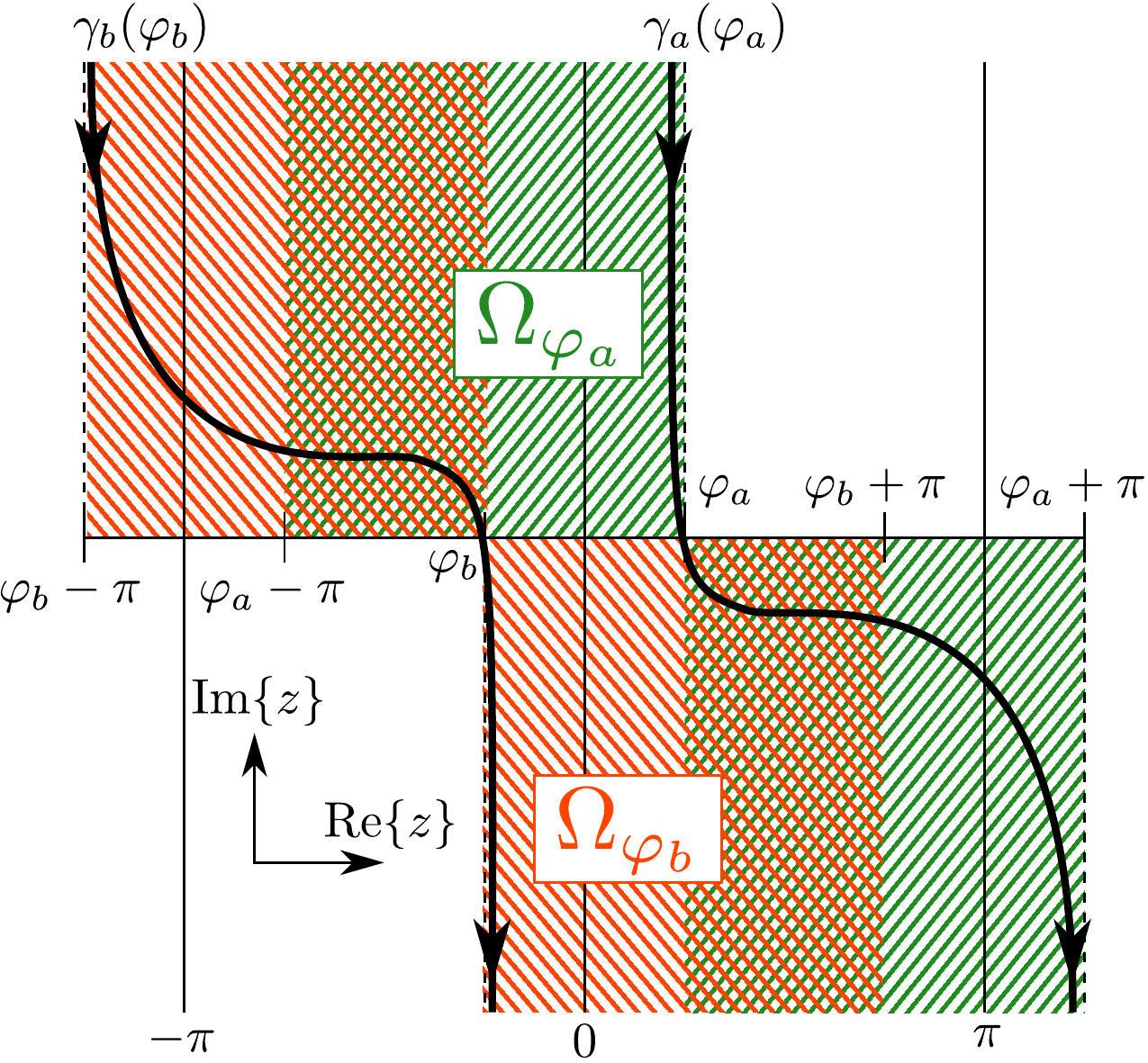}
  \includegraphics[width=0.45\textwidth]{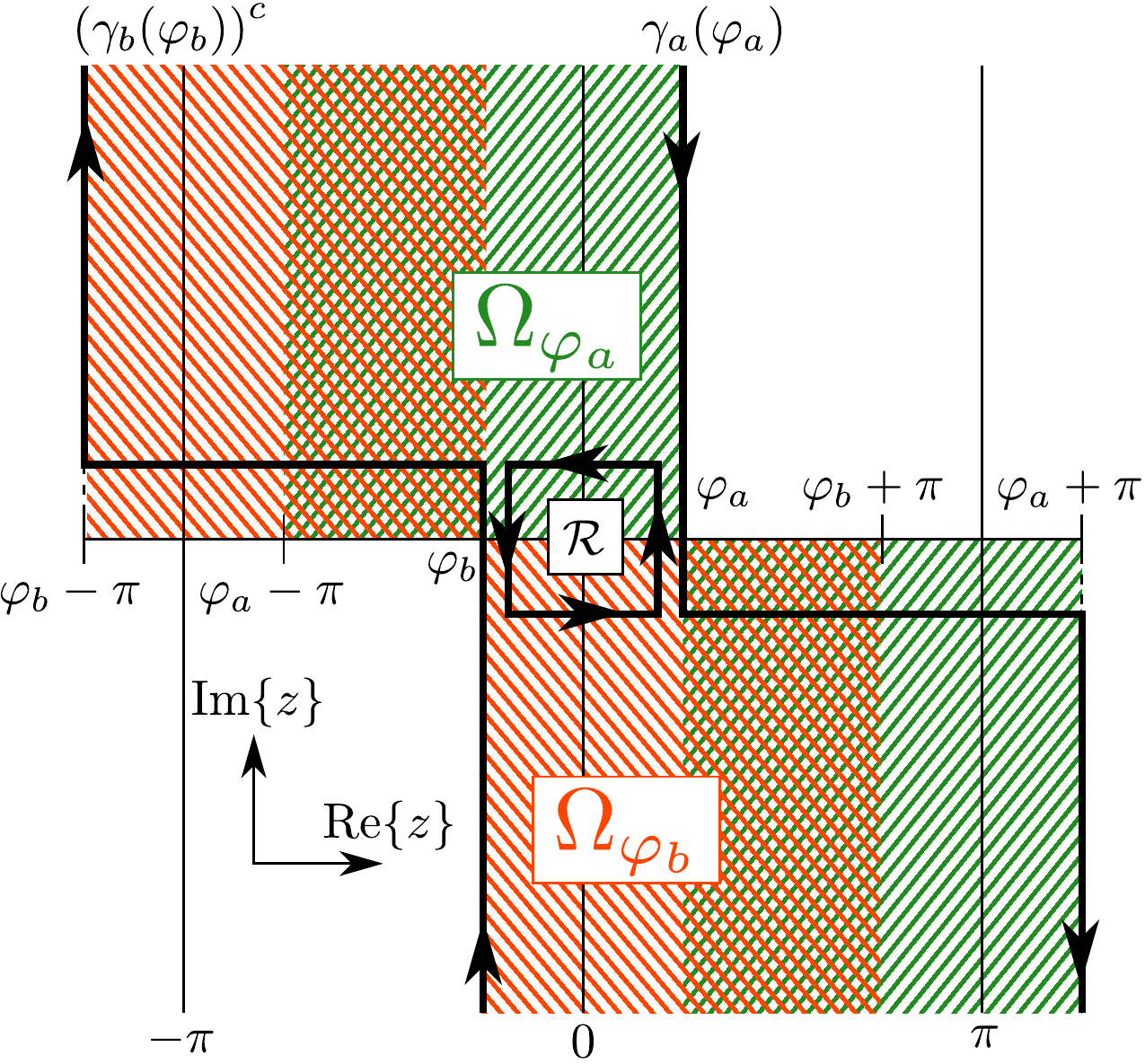}
  \caption{The two contours $\gamma_b(\varphi_b)$ and $\gamma_a(\varphi_a)$ for $-\pi<\varphi_b<0<\varphi_a<\pi$ (left), their angular counterparts and the contour $\mathcal{R}$ (right)}
\label{fig:bothcontours}
\end{figure}

We can always choose $\varphi_a$ and $\varphi_b$ close enough to zero, such that no poles related to reflected waves exist within $\mathcal{R}$. In this case, one can show that the only possible singularity is a pole corresponding to the incident wave, and we have}
\begin{align}
\frac{1}{4\pi i}\oint_{\mathcal{R}}e^{-ikr^{\star}\cos(z)}\mathcal{S}_{\varphi_b}(z)[\Phi]\mathd z=\Phi^{ab}_{\text{I}}(r^{\star},0).
\end{align}
\mylinenum{This ensures that we can write}
\begin{align}
\Phi(r^{\star},0)=\frac{1}{4\pi i}\int_{(\gamma_b(\varphi_b))^c+\gamma_a(\varphi_a)+\mathcal{R}}e^{-ikr^{\star}\cos(z)}\mathcal{S}_{\varphi_b}(z)[\Phi]\mathd z.
\end{align}
\mylinenum{Now, the coinciding lateral parts cancel each other, and the remaining contour is simply $\gamma_+ + \gamma_-$ (see Figure \ref{SMTSommcon} (left)). Now taking the limit as $\varphi_b \rightarrow 0$, or using the fact that $\mathcal{S}_{\varphi_b} (z)$ is an analytic continuation of $\mathcal{S}_0 (z)$ by Proposition \ref{prop:greenoperatorstuff}, we get}
\begin{eqnarray}
\Phi(r^{\star},0)&=&\frac{1}{4\pi i}\int_{\gamma_++\gamma_-}e^{-ikr^{\star}\cos (z)}\mathcal{S}_0(z)[\Phi]\mathd z \nonumber\\
&=&\frac{1}{2\pi i}\int_{\gamma_+}e^{-ikr^{\star}\cos(z)}\left(\frac{\mathcal{S}_0(z)[\Phi]-\mathcal{S}_0(-z)[\Phi]}{2}\right)\mathd z
\label{eq:niceformula}\end{eqnarray}
\mylinenum{Everything that has been done in this subsection can be used to get a similar formula for $\frac{\partial \Phi}{\partial \theta}$ to get}
\begin{eqnarray}
-\frac{1}{ikr^{\star}}\frac{\partial\Phi}{\partial\theta}(r^{\star},0)=\frac{1}{2\pi i}\int_{\gamma_+}e^{-ikr^{\star}\cos(z)}\sin(z)\left(\frac{\mathcal{S}_0(z)[\Phi]+\mathcal{S}_0(-z)[\Phi]}{2}\right)\mathd z
\label{eq:niceformuladerivative}\end{eqnarray}
\mylinenum{Now, comparing \eqref{eq:niceformula} and \eqref{eq:niceformuladerivative} to equations (\ref{SMTSommInt}) and (\ref{WHM-SommInt-der}), it is clear that we can apply Theorem \ref{th:MalThm} to find that}
\begin{align}
\frac{1}{2}(\mathcal{S}_0(z)[\Phi] \mp \mathcal{S}_0(-z)[\Phi])=s(z)\mp s(-z)
\end{align}
\mylinenum{leading to the sought-after formula}
\begin{eqnarray}
  s(z)=\frac{1}{2}\mathcal{S}_0(z)[\Phi], 
  \label{eq:finalconnectionformula}
\end{eqnarray}
\mylinenum{linking the spectral function $s (z)$ to the Green's operator $\mathcal{S}_0(z)[\Phi]$. Note that this formula could have been recovered from what was done at the end of the Wiener-Hopf section (Section \ref{WHT}) in particular it is a consequence of (\ref{WHM-Somm-Sol}), but this appendix is showing this link from Green's identity only. This can also be seen as a constructive way of getting to the form of the Sommerfeld integral. We can also follow the paper \citep{Malyuzhinets1958-2} to directly link the spectral function $s(z)$ with the Kontorovich-Lebedev transform of the scattered wave $\Psi$,}
\begin{align}
\nonumber\Psi(\nu,\theta)=&\ \frac{1}{\pi i\nu}\int_{-i\infty}^{i\infty}e^{-i\nu(z+\pi/2)}\left[s(\theta+z)-s(\theta-z)\right]\textrm{d}z\\
\label{Somm-KL-link}&+\frac{2(-i)^{1+\nu}}{\nu\sin(\pi\nu)}\cos\left(\left(|\theta-\theta_{\text{I}}|-\pi\right)\nu\right)
\end{align}
\mylinenum{Note that the work done in this appendix is very general and can possibly be applied to geometries other than the wedge.}

\end{document}